\documentclass[12pt]{article} 
\usepackage{amsmath}
\usepackage{amsfonts,amssymb}
\usepackage{fullpage}
\usepackage{graphicx}
\usepackage[usenames,dvipsnames]{color}
\usepackage{epsfig}
\usepackage{subfig}
\usepackage{amsmath}
\usepackage{amsfonts,amssymb}
\usepackage{amssymb}
\usepackage{latexsym}
\usepackage{url}

\def\Dn{D^{(n)}}
\def\En{E^{(n)}}

\def\In{I^{(n)}}

\def\Mn{M^{(n)}}

\def\qn{q^{(n)}}
\def\Rn{R^{(n)}}
\def\Sn{S^{(n)}}
\def\Tn{T^{(n)}}

\def\Wn{W^{(n)}}
\def\Xn{X^{(n)}}
\def\Yn{Y^{(n)}}
\def\Zn{Z^{(n)}}
\def\Ine{I^{(n)}_E}
\def\Tne{T^{(n)}_E}

\def\Sdn{S^{(n)}_D}
\def\Tdn{T^{(n)}_D}
\def\Wtn{W^{(n)}_T}
\def\Wthn{\hat{W}^{(n)}_T}

\def\qtn{\tilde{q}^{(n)}}

\def\Int{\tilde{I}^{(n)}}
\def\Inhat{\hat{I}^{(n)}}
\def\Inet{\tilde{I}^{(n)}_E}
\def\Inehat{\hat{I}^{(n)}_E}

\def\Rnhat{\hat{R}^{(n)}}
\def\Snt{\tilde{S}^{(n)}}
\def\Snhat{\hat{S}^{(n)}}
\def\Tnt{\tilde{T}^{(n)}}
\def\Tnhat{\hat{T}^{(n)}}
\def\Tnbar{\bar{T}^{(n)}}

\def\Tnt{\tilde{T}^{(n)}}
\def\Wnt{\tilde{W}^{(n)}}
\def\Wnhat{\hat{W}^{(n)}}
\def\Xnt{\tilde{X}^{(n)}}
\def\Xntprime{\tilde{X}^{(n)'}}
\def\Xnhat{\hat{X}^{(n)}}
\def\Xnprime{{X}^{(n)'}}

\def\Ynhat{\hat{Y}^{(n)}}
\def\Yncheck{\check{Y}^{(n)}}

\def\bXnbar{\bar{\boldsymbol{X}}^{(n)}}

\def\bXn{\boldsymbol{X}^{(n)}}

\def\bXnt{\tilde{\boldsymbol{X}}^{(n)}}
\def\bXnhat{\hat{\boldsymbol{X}}^{(n)}}

\def\alphan{\alpha^{(n)}}
\def\alphatn{\tilde{\alpha}^{(n)}}
\def\alphahn{\hat{\alpha}^{(n)}}
\def\betan{\beta^{(n)}}
\def\betatn{\tilde{\beta}^{(n)}}

\def\zetan{\zeta^{(n)}}
\def\zetatn{\tilde{\zeta}^{(n)}}
\def\zetahn{\hat{\zeta}^{(n)}}

\def\Estar{\mathcal{E}^*}
\def\EEn{\mathcal{E}^{(n)}}

\def\SSn{\mathcal{S}^{(n)}}

\def\E{{\rm E}}
\def\P{{\rm P}}

\newcommand\re{{\rm e}}

\def\B{\mathcal{B}}
\def\EEn{\mathcal{E}^{(n)}}

\newcommand\bx{\boldsymbol{x}}

\newcommand\bxt{\tilde{\boldsymbol{x}}}

\newcommand\rhot{\tilde{\rho}}
\newcommand\taut{\tilde{\tau}}

\newcommand\zetat{\tilde{\zeta}}
\newcommand\zetah{\hat{\zeta}}
\newcommand\At{\tilde{A}}
\newcommand\bt{\tilde{b}}

\newcommand\ct{\tilde{c}}
\newcommand\Ct{\tilde{C}}

\newcommand\gt{\tilde{g}}
\newcommand\itt{\tilde{i}}
\newcommand\ih{\hat{i}}
\newcommand\iet{\tilde{i}_E}
\newcommand\ieh{\hat{i}_E}
\newcommand\Kt{\tilde{K}}

\newcommand\st{\tilde{s}}
\newcommand\sh{\hat{s}}
\newcommand\ttild{\tilde{t}}
\newcommand\xt{\tilde{x}}
\newcommand\xh{\hat{x}}

\newcommand\Ut{\tilde{U}}

\newcommand\wt{\tilde{w}}
\newcommand\wh{\hat{w}}

\newcommand\fe{f_{\epsilon}}

\newcommand{\convas}{\stackrel{a.s.}{\longrightarrow}}

\newcommand{\convp}{\stackrel{p}{\longrightarrow}}
\newcommand{\sge}{\stackrel{st}{\ge}}
\newcommand{\sle}{\stackrel{st}{\le}}

\newtheorem{theorem}{Theorem}[section]
\newtheorem{lem}{Lemma}[section]
\newtheorem{cor}{Corollary}[section]
\newtheorem{rmk}{Remark}[section]
\newtheorem{conj}{Conjecture}[section]

\numberwithin{equation}{section}

\begin{document}


\title{Epidemics on networks with preventive rewiring}


\author{Frank Ball$^{1}$ and Tom Britton$^{2}$}
\footnotetext[1]{School of Mathematical Sciences, The University of Nottingham, University Park, Nottingham NG7 2RD, UK}
\footnotetext[2]{Department of Mathematics, Stockholm University, SE-106 91 Stockholm, Sweden}

\date{\today}
\maketitle

\begin{abstract}
A stochastic SIR (susceptible $\to$ infective $\to$ recovered) epidemic model defined on a social network is analysed. The underlying social network is described by an Erd\H{o}s-R\'{e}nyi random graph but, during the course of the epidemic, susceptible individuals connected to infectious neighbours may drop or rewire such connections. Large population limits of the model are derived giving both convergence results for the early branching process-like behaviour, and, assuming a major outbreak, the main phase of the epidemic process which converges to a deterministic model that is equivalent to a certain pair approximation model.  Law of large numbers results are also obtained for the final size (i.e.~total number of individuals infected) of a major outbreak.
Two results stand out (valid for a range of parameter set-ups): (i) the limiting final fraction infected may be discontinuous in the infection rate $\lambda$ at its threshold $\lambda_c$ (thus making a discrete jump from $0$ to a strictly positive number) and (ii) for the situation when rewiring is necessarily to uninfected individuals, if it is discontinuous, the limiting final fraction infected jumps from $0$ to $1$ as $\lambda$ passes through $\lambda_c$.
\end{abstract}




\section{Introduction}
\label{sec:intro}
Traditional mathematical models for the spread of an infectious disease assume that individuals behave the same way all through the epidemic outbreak, and that no interventions are put in place during the outbreak. More recently models have been studied that include various public health preventive measures put in place during the outbreak, thus reducing spreading (e.g.\ \cite{F06}). Usually such measures are easy to incorporate into a model, however the resulting model
is often not susceptible to mathematical analysis and the effects of measures are qualitatively understood and quantitatively obtained through large simulations studies (e.g.\ \cite{C20}).

A different type of preventive measure is those performed by separate individuals. For example, it is well known for diseases like HIV, Ebola and Covid-19 that individuals at risk also, individually, increase their preventive behaviour once they become aware of an ongoing epidemic outbreak. Recently, mathematical models containing such individual preventive measures have been defined and studied numerically (e.g.\ \cite{FSJ10}). These individually-based preventive measures, triggered by an ongoing epidemic, are however harder to analyse and few analytical results are available for epidemic models incorporating such \emph{adaptive dynamics}.

In the current paper we continue the study of one such SIR (susceptible $\to$ infective $\to$ recovered) epidemic model, where transmission occurs along links of an underlying social network, and in which uninfected individuals that are connected to infectious individuals may choose to drop or rewire connections to infectious neighbors in order to avoid getting infected.  \cite{LBSB18} study such a model and prove that for some (but not all!) networks and parameter set-ups, the model with preventive rewiring may actually give a \emph{larger} epidemic outbreak than the corresponding network and parameter set-up but without rewiring.
The surprising conclusion is hence that such preventive rewiring at the individual level \emph{reduces} the risk for an individual to get infected, but if collectively performed it may \emph{increase} the final epidemic size.

Durrett and co-workers (\cite{JKYJD:2019}) study the same model for the special case where the underlying social network is an Erd\H{o}s-R\'{e}nyi network
(i.e.~one in which between any two distinct individuals an edge is present independently with a fixed probability)
and conjecture, by means of a convincing simulation study, that the final epidemic size may be discontinuous at criticality. By this we mean that, fixing the Erd\H{o}s-R\'{e}nyi network model and all other model parameters apart from the transmission rate $\lambda$, the limiting final epidemic outbreak size $\tau(\lambda)$ (more specifically, the fraction of the population infected by an epidemic initiated by one infective, given the occurrence of a \emph{major outbreak}, i.e.~one that takes off and becomes established)
satisfies
\begin{equation}
\tau(\lambda)
\begin{cases}
	      =0& \text{ if } \lambda \le \lambda_c, \\
	      \ge a& \text{ if } \lambda > \lambda_c,
\end{cases}
\end{equation}\label{equ:tau-disc}
$\lambda_c$ being the value for which the basic reproduction number $R_0=1$, and $a$ being a strictly positive constant. In~\cite{JKYJD:2019} the
authors provide details of various approximations to the epidemic model with rewiring (none of which displayed the discontinuity) in the hope that
someone can find an accurate approximation that explains the reason for discontinuous phase transition at $\lambda=\lambda_c$.  The main
achievement of the present paper is to provide such a deterministic approximation, together with an associated law of large numbers for the temporal trajectory of the stochastic epidemic as the population size $n \to \infty$.  We give almost identical necessary and sufficient conditions for the final
size of this deterministic model to display the discontinuity~\eqref{equ:tau-disc}. For the stochastic model, we prove conditions under which the final size of a major outbreak satisfies~\eqref{equ:tau-disc} in the limit as $n \to \infty$ and conditions under which it does not satisfy~\eqref{equ:tau-disc}.

Our deterministic approximation is derived from a novel construction of the stochastic epidemic with rewiring, in which the Erd\H{o}s-R\'{e}nyi
network and epidemic
are constructed simultaneously.  A key observation underlying this construction is that owing to symmetries in the Erd\H{o}s-R\'{e}nyi
network and the rewiring process it is sufficient to keep track of just the number (and not the end points) of susceptible-susceptible rewired edges.  The construction yields
an infinite-dimensional continuous-time Markov chain. However the epidemic dynamics are encapsulated by the four-dimensional process which records the evolution of the numbers of susceptibles, infectives, infective-susceptible edges and susceptible-susceptible rewired edges with time.  We prove
weak law of large numbers convergence of this four-dimensional process to the above-mentioned deterministic approximation using a theorem in~\cite{DN:08}.

The construction yields a process whose law is very close to but not exactly the same as that of the original model.  However, the difference is very small and we prove that convergence in probability results transfer from the constructed to the original process. The approximating deterministic model is
equivalent to the pair-approximation model of the epidemic, so the above law of large numbers provides a rigorous justification of the pair-approximation  model.

A common approach to analysing the final size of an SIR epidemic is via a suitable random time-scale transformation which leaves the final size unchanged.
For the present model this yields a deterministic limiting process that is not Lipschitz in the vicinity of disease-free states, which means that
results concerning convergence of terminal values, in for example~\cite{DN:08}, cannot be applied.  This explains why only partial results, derived by considering appropriate bounding processes, are proved concerning discontinuity of the final size of the stochastic model.  For the special case of the
SI model, in which there is no recovery from infection so infectious individuals remain so forever, this difficulty disappears.  Moreover, for the SI model, the deterministic
limiting process can be analysed to yield an explicit non-linear equation for the final size $\tau$, which enables us to prove a necessary and sufficient
condition for there to be a discontinuity at the threshold $\lambda=\lambda_c$.

In the above models, when a susceptible individual rewires an edge from an infective neighbour it is to an individual chosen uniformly at random from all the other $n-2$ individuals in the population.  We analyse also an alternative model in which such rewiring is to an individual chosen uniformly at random
from all other \emph{susceptible} individuals.  This model also does not suffer from the Lipschitz problem and a full rigorous analysis of the stochastic model is available.  In particular, a necessary and sufficient condition for the (in probability) limiting final size of a major outbreak $\tau(\lambda)$ to have a discontinuity at $\lambda_c)$ is proved.  Moreover, when it exists, the discontinuity at $\lambda_c$ is more striking in that
\begin{equation}
\label{equ:tau-disc-susonly}
\tau(\lambda)
\begin{cases}
	      =0& \text{ if } \lambda \le \lambda_c, \\
	      =1& \text{ if } \lambda_c<\lambda \le \omega-\gamma,\\
           \in (0,1)& \text{ if } \lambda > \omega-\gamma,
\end{cases}
\end{equation}
where $\omega$ and $\gamma$ are the rewiring and recovery rates, respectively.  (A necessary but not sufficient condition for $\tau(\lambda)$ to be discontinuous at $\lambda_c$ is $\omega > \gamma$.)  Thus the fraction of the population infected by a major outbreak jumps from $0$ to $1$
when $\lambda$ passes through the threshold $\lambda_c$.  When $\lambda$ is very large susceptibles become infected before they can rewire away from
neighbouring infectives and in the limit $n \to \infty$ an epidemic that takes off infects all members of the giant component of the underlying
Erd\H{o}s-R\'{e}nyi network.

While we were completing our research, Durrett and another co-worker (\cite{YD:20}) published a second paper concerned with epidemics with rewiring on
random graphs generated by the configuration model, in which conditions (that are almost complementary) for the final size of the SIR
epidemic to have and not have a discontinuity at $\lambda_c$ are proved.  They also prove similar results for the SI model.  When the degree distribution in the configuration model graph is chosen to be Poisson, the conditions in~\cite{YD:20} for $\tau(\lambda)$ to have a discontinuity at $\lambda_c$
coincide with those we give for epidemics on an Erd\H{o}s-R\'{e}nyi network but (at least for the SI model) the function $\tau(\lambda)$ in~\cite{YD:20}
is different from ours.  In our work an explicit equation satisfied by $\tau(\lambda)$ is not available for the SIR model.  As we discuss later,
the results in~\cite{YD:20} are obtained by analysing a construction of the network and epidemic which we believe is not probabilistically equivalent
to the original model, and which yields an incorrect deterministic limit as $n \to \infty$ (even if numerically quite close).

The remainder of the paper is organised as follows.  The main SIR model with rewiring is described in Section~\ref{sec:model}.  It is slightly more general than that indicated above in that susceptible individuals may also drop edges to infective neighbours.  The main results of the paper are given in Section~\ref{sec:mainresults}.  The SIR model is treated in Section~\ref{sec:SIR}, where theorems concerning branching process approximation of the early stages of an epidemic with one initial infective (Theorem~\ref{thm:BPapprox}) and a weak law of large numbers for the temporal behaviour of an epidemic
in which a strictly positive fraction of the population are initially infective (Theorem~\ref{thm:temporalWLLN}) are stated.  A weak law of large numbers for the final size of a major outbreak is conjectured (Conjecture~\ref{conj:SIRfinal}).  Theorems giving conditions for the final size of the
approximating deterministic model and stochastic model to have a discontinuity at the threshold $\lambda_c$ are stated (Theorems~\ref{thm:detdisc} and~\ref{thm:SIRdiscon}, respectively).

The SI model is treated in Section~\ref{sec:SI}.  A weak law of large numbers for the final size of a major outbreak is given in Theorem~\ref{thm:SIfinalsize} and a necessary and sufficient condition for the final size $\tau(\lambda)$ to have a discontinuity at $\lambda_c$
is presented in Theorem~\ref{thm:SIdisc}.  Note that a major outbreak is not guaranteed even though infected individuals remain so forever, as the initial infective may not belong to the giant component of the original network, and even if it does, susceptibles can rewire their connections away from
infected neighbours.  A detailed analysis of the final size of a major outbreak is given in Corollary~\ref{cor:SIfinal}.  The SIR model in which rewiring is necessarily to susceptibles is treated in Section~\ref{sec:SIRsusonly}.  Weak laws of large numbers for the temporal behaviour of an epidemic with a positive fraction initially infected and for the final size of a major outbreak for an epidemic with one initial infective are given in Theorems~\ref{thm:temporalWLLNsusonly} and~\ref{thm:SIfinalsizesusonly}, respectively.  The latter includes the behaviour described by~\eqref{equ:tau-disc-susonly}.  Two other rewiring models are treated briefly in Section~\ref{sec:altrewire}.  Throughout Section~\ref{sec:mainresults} results are illustrated by simulations and numerical studies.

The construction of the epidemic model which underpins its analysis is described in Section~\ref{sec:construction}.  The approximating deterministic model is shown to be equivalent to a pair-approximation model of the epidemic in Section~\ref{sec:pairapprox}.  A comparison of our results with those in~\cite{YD:20} is given in Section~\ref{sec:YaoDurrett}, where a numerical study shows that the corresponding asymptotic final sizes differ and a simulation study gives very strong support to the correctness of our asymptotic final size.  A heuristic derivation of the final size of the SI model is given in Section~\ref{sec:heuristic}.  The proofs are given in Section~\ref{sec:proofs} and some very brief concluding comments are given in Section~\ref{sec:conc}.


\section{The model, results and construction}
\label{sec:modelmain}
\subsection{The network epidemic model  with social distancing}
\label{sec:model}
Consider a community of size $n$ socially structured by the Erdös-Renyi random graph $G(n,\mu/n)$, i.e.\ where there is an edge, independently, between each pair of individuals with probability $\mu/n$.

Given the social graph/network a Markovian SIR epidemic process is defined as follows. The epidemic starts by one or more individuals, selected uniformly at random from the population, being infectious and the rest being susceptible. During the course of the epidemic, each infectious individual infects each susceptible neighbour at rate $\lambda$ (i.e.~at the points of independent Poisson processes each having rate $\lambda$, and each infectious individual recovers (and becomes immune) at rate $\gamma$ (implying that the duration of the infectious period follows an exponential distribution having mean $\gamma^{-1}$). The preventive feature of the model is that each susceptible individual rewires or drops its edges to any infectious neighbour independently at rate $\omega$
for each edge. More precisely, the susceptible individual rewires such an edge with probability $\alpha$, otherwise the susceptible individual drops the  edge. When a susceptible rewires an edge, the rewiring is to an individual chosen uniformly at random from the other $n-2$ individuals in the population. (In  Section~\ref{sec:SIRsusonly} we treat also the case where the rewiring is necessarily to a susceptible individual, and in Section~\ref{sec:altrewire} we outline the case where the rewiring is to anyone except those currently infectious.) All events described above are defined to occur mutually independently.  The epidemic stops when there is no infective individual remaining in the population.

The parameters of the model are hence: $\mu$: the mean degree of the underlying graph, $\lambda$ the transmission rate, $\gamma$: the recovery rate, $\omega$: the rewiring/dropping rate, and $\alpha$: the probability of rewiring (rather than dropping).

Throughout the paper we assume implicitly that $\omega>0$, unless explicitly stated otherwise. If $\omega=0$ the model reduces to an SIR epidemic on an Erd\H{o}s-R\'{e}nyi random graph having known large-$n$ behaviour (\cite{AB00} and \cite{Neal:2003}).

\subsection{Results}
\label{sec:mainresults}
\subsubsection{SIR model}
\label{sec:SIR}

We first introduce some notation, the precise meaning of which is made clear when we describe our construction of the epidemic model in Section~\ref{sec:construction}.  Let $\bXn=\{\bXn(t):t \ge 0\}$, where $\bXn(t)=(\Sn(t), \In(t), \Ine(t), \Wn(t))$ and $\Sn(t)$, $\In(t)$, $\Ine(t)$ and $\Wn(t)$ are respectively the numbers of susceptibles, infectives, infectious (i.e.~infective-susceptible) edges and susceptible-susceptible rewired edges at time $t$.

Our first result concerns the initial phase of an epidemic starting with one initial infective, which we approximate by a branching process.  Let $\B$ be a continuous-time branching process, with one ancestor, in which the lifetime of an individual follows an exponential distribution with mean $\gamma^{-1}$. At birth an individual is assigned ${\rm Po}(\mu)$ infectious edges, where ${\rm Po}(\mu)$ denotes a Poisson random variable with mean $\mu$.  An individual drops each of its infectious edges independently at rate $\omega$ and infects down them independently at rate $\lambda$.  When an individual infects down an infectious edge a new individual is born and the edge is dropped. When an individual dies all of its remaining infectious edges are dropped.  For $t \ge 0$, let $I(t)$, $I_E(t)$ and $T(t)$ be respectively the number of individuals, number of infectious edges and the total progeny (including the initial ancestors) in $\B$ at time $t$.  For $n=1,2,\dots$, let $t_n=\inf\{t \ge 0: T(t)\ge \log n\}$, where $t_n=\infty$ if $T(t) <  \log n$ for all $t \ge 0$.  Let $\EEn$ be the epidemic described in Section~\ref{sec:model}, assuming a single initial infective.
Let $\Tn(t)=n-\Sn(t)$ be the total number of infections in $\EEn$ during $[0,t]$.

\begin{theorem}
\label{thm:BPapprox}
There exists a probability space $(\Omega, \mathcal{F}, \P)$ on which are defined realisations of $\EEn$ $(n=1,2,\dots)$ and $\B$ satisfying
\[
\sup_{0 \le t \le t_n}\left|\left(\In(t), \Ine(t), \Tn(t)\right)-\left(I(t), I_E(t), T(t)\right)\right| \convp 0 \quad \mbox{as } n \to \infty.
\]
\end{theorem}

Let
\begin{equation}
\label{equ:Rzero}
R_0=\frac{\mu\lambda}{\lambda+\omega+\gamma}
\end{equation}
be the offspring mean of $\B$. The quantity $R_0$ is known as the basic reproduction number for the epidemic $\EEn$.  We say that a \emph{major outbreak} occurs in $\EEn$ if at least $\log n$ individuals are infected.  It follows from Theorem~\ref{thm:BPapprox} that as $n \to \infty$ the probability of a major outbreak converges to the probability that $\B$ does not go extinct.  Thus with high probability a major outbreak is possible if and only if $R_0>1$.
(An event, $A_n$ say, is said to hold with high probability if $\P(A_n) \to 1$ as $n \to \infty$.)
In particular, if all parameters other than $\lambda$ are held fixed, with high probability a major outbreak is possible if and only if $\lambda>\lambda_c$, where
\begin{equation}
\label{equ:lambdac}
\lambda_c=\frac{\gamma+\omega}{\mu-1}.
\end{equation}

Let $\Tn=n-\Sn(\infty)$ be the final size of the epidemic $\EEn$, i.e.~the total number of individuals infected during the epidemic, and $\Tnbar=n^{-1}\Tn$ be the fraction of the population that become infected.  The above criticality of $\lambda_c$ and the following lemma follow from~\cite{JKYJD:2019}, Theorem 2,
though the proof there is quite different from the current proof, given in Section~\ref{sec:proofBP}.
\begin{lem}
\label{lem:majbound}
Suppose that $R_0>1$.  Then there exists $\tau'=\tau'(\mu, \lambda, \gamma, \omega)>0$ such that
\[
\lim_{n \to \infty} \P(\Tnbar \ge \tau'|\Tn \ge \log n)=1.
\]
\end{lem}

\begin{rmk}
\label{rmk:dropping}
Note that $\lambda_c$ and $\tau'$ are independent of the value of $\alpha$.  It is clear that the model with $\alpha=0$, so there is dropping but no rewiring of edges, provides a lower bound for any corresponding model with $\alpha>0$.  For the model with $\alpha=0$, the epidemic and random graph
can be constructed simultaneously on a generation basis using an obvious extension of the construction in~\cite{Neal:2003} to yield a process that
can be represented as a randomised Reed-Frost process (\cite{ML:1986}).  The central limit theorem in~\cite{ML:1986} then yields a central limit theorem for
the final size of a major outbreak.  We omit the details.
\end{rmk}

Our next result is concerned with the deterministic approximation of the main body of an epidemic started with many initial infects.  For $t \ge 0$, let $x(t)=(s(t), i(t), i_E(t), w(t))$ be the solution of the system of ordinary differential equations
\begin{align}
\dfrac{ds}{dt}&=-\lambda i_E, \label{equ:dsdt}\\
\dfrac{di}{dt}&=-\gamma i + \lambda i_E\label{equ:didt}\\
\dfrac{di_E}{dt}&= -\lambda i_E -\gamma i_E +\lambda \mu i_E s -\lambda \frac{i_E^2}{s} +2 \lambda i_E \frac{w}{s}-\omega i_E(1-\alpha+\alpha(1-i)),\label{equ:diEdt}\\
\dfrac{dw}{dt}&= \omega \alpha i_E s-2 \lambda i_E \frac{w}{s},\label{equ:dwdt}
\end{align}
having initial condition $x(0)=(s(0), i(0), i_E(0), w(0))$.  (The vector $x(t)$ is not in bold to help link with the theory in~\cite{DN:08}, where $\bx(t)$ has a separate meaning, that is used in the proofs.)  Let $\bXnbar(t)=n^{-1}\bXn(t)$.

\begin{theorem}
\label{thm:temporalWLLN}
Suppose $\bXnbar(0) \convp x(0)$ as $n \to \infty$, where $i(0)>0$ and $i_E(0)>0$.  Then, for any $t_0>0$,
\[
\sup_{0 \le t \le t_0} \left|\bXnbar(t)-x(t)\right| \convp 0 \quad \mbox{as } n \to \infty.
\]
\end{theorem}

\begin{rmk}
We show in Section~\ref{sec:pairapprox} that~\eqref{equ:dsdt}-\eqref{equ:dwdt} are equivalent to the pair-approximation model of the epidemic, so
Theorem~\ref{thm:temporalWLLN} provides a rigorous justification of the pair-approximation model of the epidemic $\EEn$ when a positive fraction of individuals are initially infective; cf.~\cite{Altmann98} and~\cite{JBTR:2018} who prove similar weak law of large numbers for other dynamic network epidemic models.
\end{rmk}

\begin{figure}
\begin{center}
\resizebox{0.49\textwidth}{!}{\includegraphics{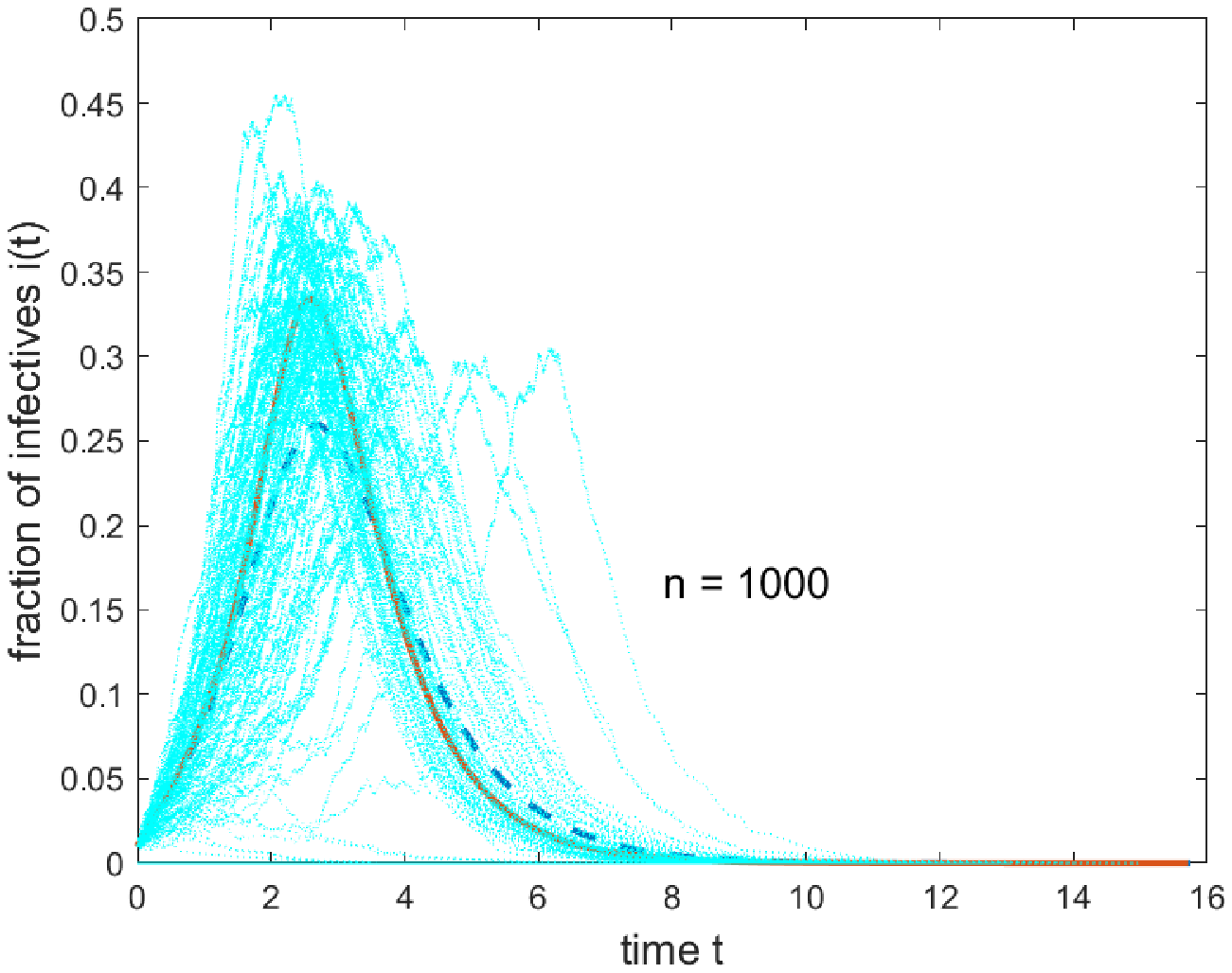}}
\resizebox{0.49\textwidth}{!}{\includegraphics{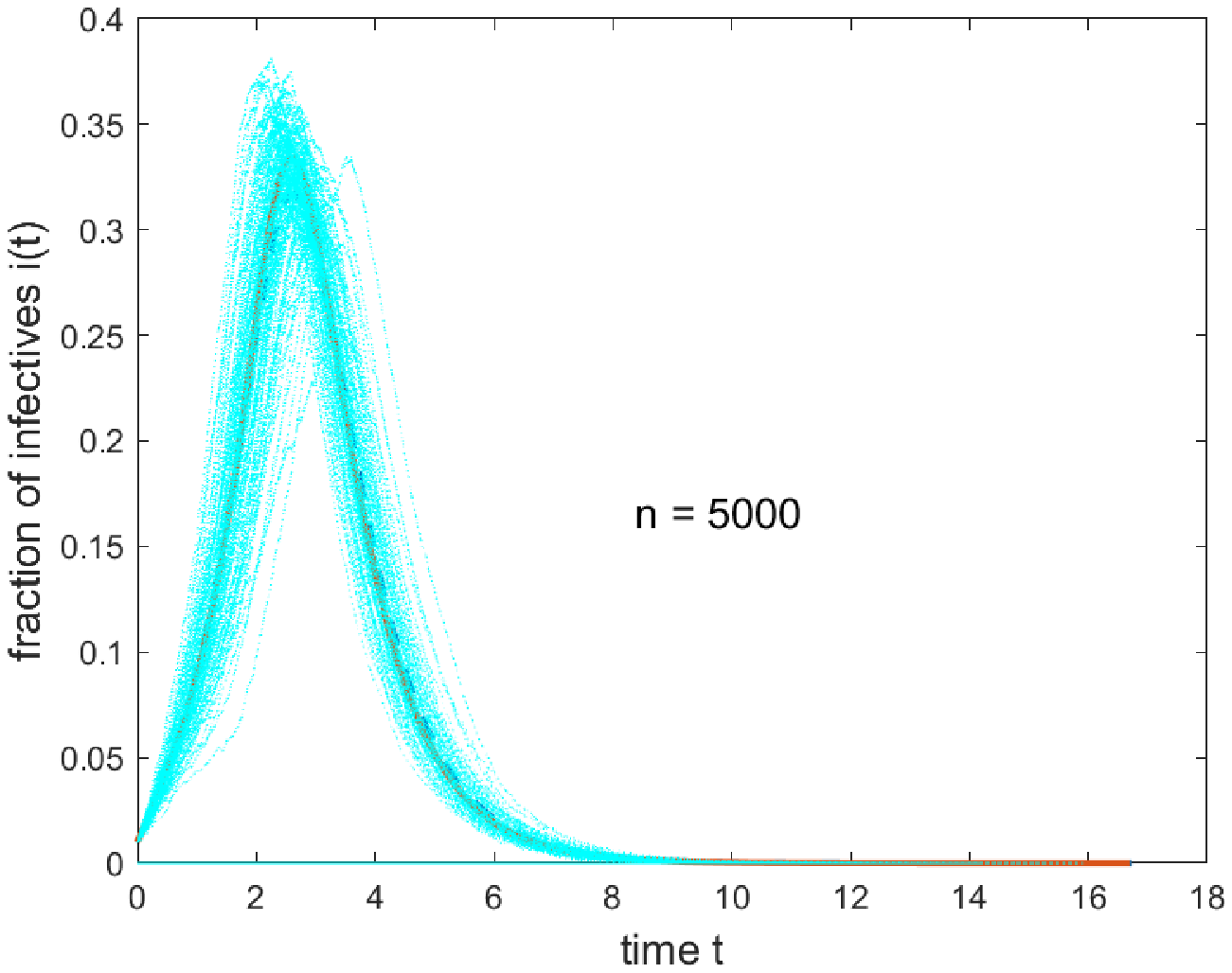}}\\
\end{center}
\caption{$100$ simulated realisations of trajectories of fraction infected in SIR epidemics in populations of size $n=1,000$ and $n=5,000$ with $1\%$ initially infective, when
$\mu=5$, $\lambda=1.5$, $\gamma=1$, $\omega=4$ and $\alpha=1$.  Also shown is the deterministic fraction $i(t)$ (solid curve) and the mean of the stochastic
trajectories (dashed curve).  See text for details.}
\label{fig:wllnpic}
\end{figure}

Theorem~\ref{thm:temporalWLLN} is illustrated in Figure~\ref{fig:wllnpic}, which is based on simulations of epidemics with $\lambda=1.5$, $\gamma=1$, $\omega=4$ and $\alpha=1$ (so no dropping of edges) on networks with $\mu=5$.  For each of $n=1,000$ and $n=5,000$, a realisation of the  Erdös-Renyi random graph $G(n,\mu/n)$ was simulated and then an epidemic, in which initially $1\%$ of the population were infected and the remaining $99\%$ susceptible,  was simulated on that graph, with the whole process being repeated $100$ times.  The trajectories of the fraction of the population infected in the $100$ simulations are shown, together with the trajectory of their mean (dashed curve) and the deterministic fraction $i(t)$ (solid curve).  When $n=1,000$, there are only $10$ initial infectives.  Consequently a few of the epidemics failed to take off and for those that took off there was noticeable variation in the time of the peak of the process of infectives.  Both of these phenomena have the effect of reducing the mean trajectory and can be explained by considering the branching process $\B$ which approximates the early stages of an epidemic.  See~\cite{BR:2013} for a formal proof of the latter phenomenon for a range of epidemic models.  When $n=5,000$, there are $50$ initial infectives, all $100$ simulated epidemics took off and there was appreciably less variation in the time of the peak.  The initial conditions are now closer to those of Theorem~\ref{thm:temporalWLLN} and the deterministic model provides a good approximation to the corresponding stochastic model.

We now consider the final outcome of an epidemic.  Recall that $\Tn=n-\Sn(\infty)$ is the final size of the epidemic and let $\Tnbar=n^{-1}\Tn$ be the fraction of the population that become infected.

\begin{conj}
\label{conj:SIRfinal}
\begin{enumerate}
\item[(a)] Suppose $\bXnbar(0) \convp x(0)$ as $n \to \infty$, where $i(0)>0$ and $i_E(0)>0$.  Then
\[
\Tnbar \convp 1-x(\infty) \quad \mbox{as } n \to \infty.
\]
\item[(b)]
Suppose that $R_0>1$ and for each $n$ the epidemic is started by $1$ infective, with the
rest of the population susceptible.  For $ \epsilon \in (0,1)$, let $x^{\epsilon}(t)=(s^{\epsilon}(t), i^{\epsilon}(t), i_E^{\epsilon}(t), w^{\epsilon}(t))$ be the solution
of~\eqref{equ:dsdt}-\eqref{equ:dwdt} with $x^{\epsilon}(0)=(1-\epsilon, \epsilon, L^{-1} \epsilon,0)$, where
$L=\frac{\lambda}{\lambda(\mu-1)-\omega}$, and $\tau=1-\lim_{\epsilon \downarrow 0} s^{\epsilon}(\infty)$.  Then, conditional upon a major outbreak,
\[
\Tnbar \convp \tau \quad \mbox{as } n \to \infty.
\]
\end{enumerate}
\end{conj}

\begin{rmk}
\label{rmk:proofSIRfinal}
We present arguments in support of Conjecture~\ref{conj:SIRfinal} in Section~\ref{sec:proofthmSIRfinal}, although we do not have a complete proof.
The usual approach to proving limit theorems for the final size of an SIR epidemic models is via a random time-scale transformation (see the start of Section~\ref{sec:proofthmSIfinal}).  This leads to a deterministic model in which the right-hand sides of~\eqref{equ:dsdt}-\eqref{equ:dwdt} are divided
by $\lambda i_E$ (see~\eqref{equ:didtB}-\eqref{equ:dwdtB} in Section~\ref{sec:conjSIRsupp}).  Note that the vector field corresponding to the time-transformed deterministic model is not Lipschitz in the neighbourhood of $(i, i_E)=(0,0)$ owing to the term $i/i_E$ in the transformed version of~\eqref{equ:didt}, so standard results for the approximation of Markov chains by differential equations cannot be applied.  Note that for the SI epidemic
(i.e.~when $\gamma=0$), $i(t)=1-s(t)$ for all $t \ge 0$. The time-transformed differential equations for $(s,i_E,w)$ form a closed system that is Lipschitz,  so the above difficulty disappears and rigorous results are readily available.
\end{rmk}

Write $\tau$ in Conjecture~\ref{conj:SIRfinal}(b) as $\tau_{SIR}(\mu, \lambda, \gamma, \omega, \alpha)$ to show explicitly its dependence on the parameters of the epidemic and let $\tau_{SIR}(\mu, \lambda, \gamma, \omega, \alpha)=0$ if $R_0 \le 1$. Our next result gives a sufficient (and almost necessary) condition for $\tau_{SIR}(\mu, \lambda, \gamma, \omega, \alpha)$ to be discontinuous at $\lambda=\lambda_c$ when all other parameters are held fixed.
\begin{theorem}
\label{thm:detdisc}
\[
\lim_{\lambda \downarrow \lambda_c} \tau_{SIR}(\mu, \lambda, \gamma, \omega, \alpha)
\begin{cases}
	      =0& \text{ if } \gamma >\omega(2\alpha-1) \text{ or } \mu < \frac{2\omega\alpha}{\omega(2\alpha-1)-\gamma}, \\
	      >0& \text{ if } \gamma <\omega(2\alpha-1) \text{ and } \mu > \frac{2\omega\alpha}{\omega(2\alpha-1)-\gamma}.
\end{cases}
\]
\end{theorem}

Theorem~\ref{thm:detdisc} would immediately yield rigorous results concerning a discontinuity at $\lambda=\lambda_c$ for the final size of stochastic epidemics that take off if Conjecture~\ref{conj:SIRfinal}(b) is true.
By considering suitable bounding processes we can however obtain conditions, together with a proof, under which the final size of an epidemic started with a single infective
that take off is discontinuous at the threshold and other conditions when it is continuous at the threshold $\lambda=\lambda_c$.  Our main result is that
the final size may be discontinuous at the threshold $\lambda=\lambda_c$.
\begin{theorem}
\label{thm:SIRdiscon}
Suppose that for each $n$ the epidemic is started by $1$ infective, with the
rest of the population susceptible.
\begin{enumerate}
\item[(a)]
Suppose that $\omega(2\alpha-1)>\gamma$ and $\mu > \frac{2\omega\alpha}{\omega(2\alpha-1)-\gamma}$.  Then there exists $\tau_0=\tau_0(\mu, \lambda, \omega, \alpha)>0$ such that, conditional upon a major outbreak,
\[
\lim_{\lambda \downarrow \lambda_c} \lim_{n \to \infty} \P(\Tnbar > \tau_0)=1.
\]
\item[(b)]
Suppose that $\omega(3\alpha-1) \le \gamma$ or $\mu \le \frac{3\omega\alpha}{\omega(3\alpha-1)-\gamma}$.
Then, conditional upon a major outbreak, for all $c>0$,
\[
\lim_{\lambda \downarrow \lambda_c} \lim_{n \to \infty} \P(\Tnbar < c)=1.
\]
\end{enumerate}
\end{theorem}

\begin{rmk}
\label{rmk:disc}
Theorem~\ref{thm:SIRdiscon} is not fully satisfactory as there is a gap between the conditions in parts (a) and (b).   Note that Conjecture~\ref{conj:SIRfinal}(b) and Theorem~\ref{thm:detdisc} imply the conjecture that the condition in Theorem~\ref{thm:SIRdiscon}(b) can be replaced by
\begin{equation}
\label{equ:disccond}
\omega(2\alpha-1) < \gamma \quad\mbox{or}\quad  \mu < \frac{2\omega\alpha}{\omega(2\alpha-1)-\gamma}.
\end{equation}
See Section~\ref{sec:proofthmSIRfinal} for further evidence in support of this conjecture is provided, which suggests also that the strict inequalities
in~\eqref{equ:disccond} can be replaced by weak inequalities.
\end{rmk}

\begin{figure}
\begin{center}
\resizebox{0.49\textwidth}{!}{\includegraphics{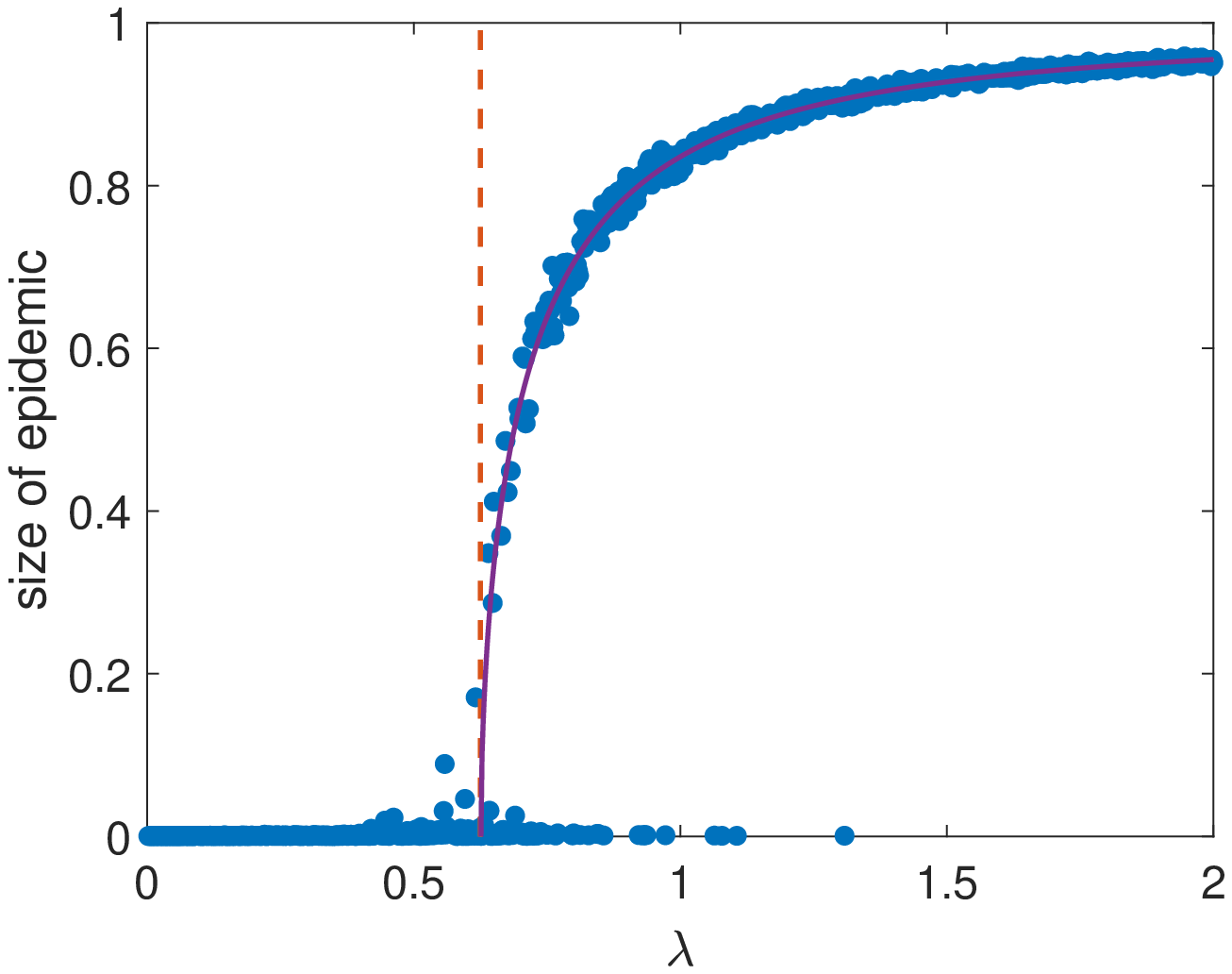}}
\resizebox{0.49\textwidth}{!}{\includegraphics{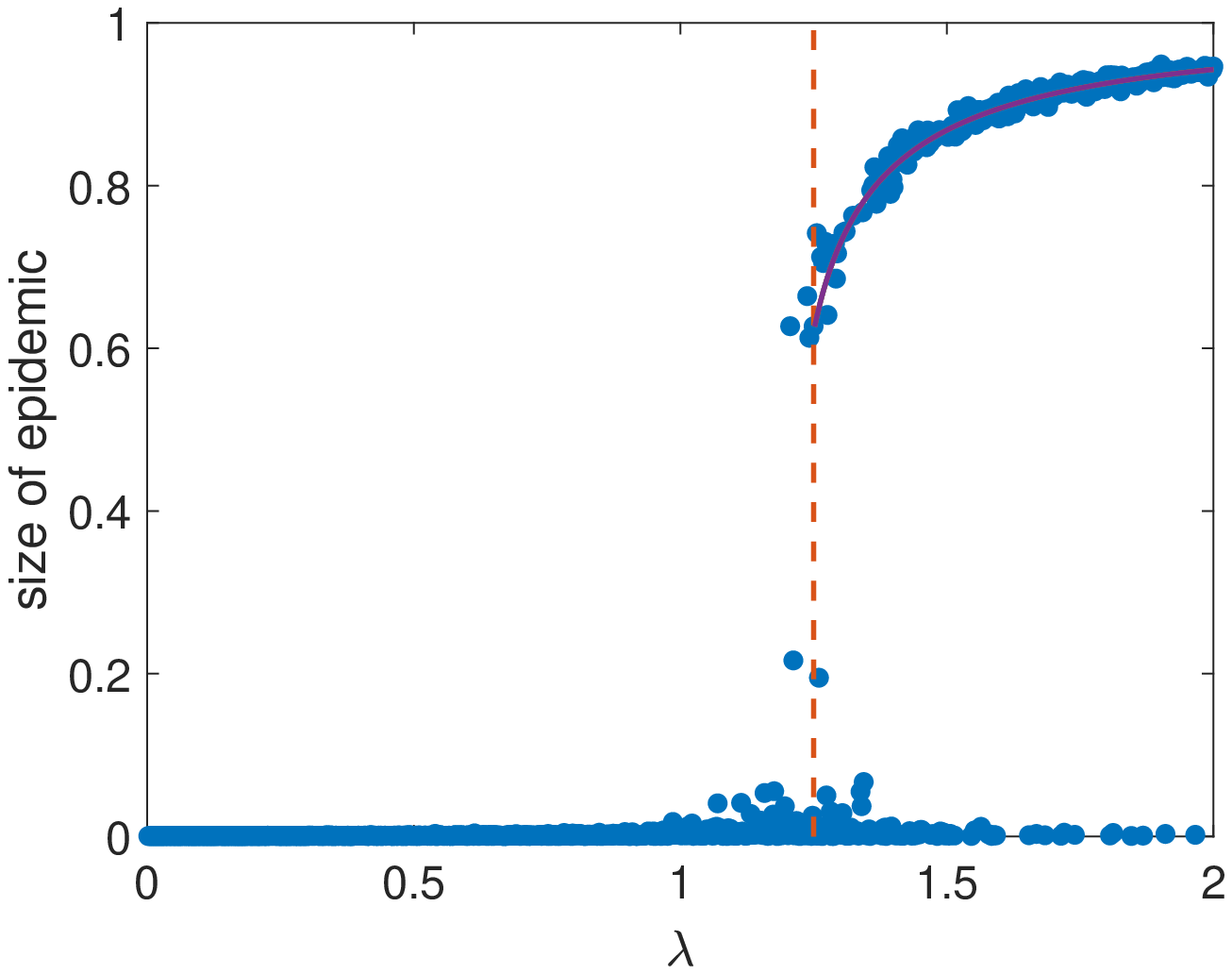}}
\end{center}
\caption{1,000 simulations of final size of SIR epidemic when $n=10,000, \mu=5, \gamma=1, \alpha =1$ and varying $\lambda$; $\omega=\frac{3}{2}$ in the left panel and $\omega=4$ in the right panel. See text for details.}
\label{fig:sirfinal}
\end{figure}

Figure~\ref{fig:sirfinal}, which is inspired by Figure 4 of~\cite{JKYJD:2019}, shows the results of simulations of the final outcome of the SIR epidemic with rewiring.  In each of the two plots the final size (expressed as a fraction of the population) of $1,000$ epidemics in a population of size $n=10,000$,
with values of $\lambda$ equally spaced in $[0, 2]$, $\gamma=1$, $\alpha=1$ (so there is no dropping of edges) and $\omega=\frac{3}{2}$ (left panel) and
$\omega=4$ (right panel).  Each simulation was initiated by $5$ infectives, chosen uniformly at random from the population, with the rest of the population being susceptible.  The value of $\omega$ in the right panel is that used in~\cite{JKYJD:2019}, Figure 4, for which Theorem~\ref{thm:SIRdiscon}(a) indicates that the final size is discontinuous at the threshold $\lambda=\lambda_c$.  The value of $\omega$ in the left panel is chosen so that our conjecture predicts the final size is continuous at $\lambda=\lambda_c$.  The solid curves in the two panels show the limiting fraction infected predicted by Conjecture~\ref{conj:SIRfinal}(b); they are calculated numerically solving the differential equation for $x^{\epsilon}(t)$ with $\epsilon=10^{-10}$ and deeming the deterministic epidemic finished when $i^{\epsilon}(t)$ crosses $10^{-8}$ from above.  Note the excellent agreement between the stochastic simulations and the asymptotic limits in both panels, lending credence to Conjecture~\ref{conj:SIRfinal}(b), and the marked discontinuity at $\lambda=\lambda_c$ in the right panel.  In both panels some large epidemics occur when $\lambda \le \lambda_c$.  This is a finite-population effect
even though $n$ is relatively large.

\subsubsection{SI model}
\label{sec:SI}

In this section we assume that $\gamma=0$, so there is no recovery and the model becomes SI.  We focus on the final size of the epidemic; results concerning
temporal behaviour follow by setting $\gamma=0$ in corresponding results for the SIR model.  Note that as infectives remain so forever the epidemic stops
when there is no edge between a susceptible and infective individuals.  For this model we obtain an explicit equation satisfied by
the limiting fraction ultimately infected as $n \to \infty$ (see Theorem~\ref{thm:SIfinalsize}), which enables analysis of the dependence of that fraction on model parameters and proof of a necessary and sufficient condition for its discontinuity (Corollary~\ref{cor:SIfinal}).  For $\epsilon \in [0,1)$, define $F_{\epsilon}:[\epsilon,1] \to \mathbb{R}$ by
\begin{equation}
\label{equ:Feps}
F_{\epsilon}(x)=1-x-(1-\epsilon)\exp\left(-\frac{(\lambda\mu+\omega\alpha)x-\omega\alpha\epsilon}{\lambda+\omega(1-\alpha)+2\omega\alpha(1-x)}\right).
\end{equation}

\begin{theorem}
\label{thm:SIfinalsize}
\begin{enumerate}
\item[(a)] Suppose that $\bXnbar(0) \convp (1-\epsilon, \epsilon, \mu \epsilon(1-\epsilon),0)$ as $n \to \infty$, where $\epsilon \in (0, 1)$ and let $\tau$
be the smallest solution in $(\epsilon,1)$ of $F_{\epsilon}(x)=0$.
Then, provided $F_{\epsilon}'(\tau)<0$,
\[
\Tnbar \convp \tau \quad \mbox{as } n \to \infty.
\]

\item[(b)] Suppose that $R_0>1$ and for each $n$ the epidemic is started by one infective, with the
rest of the population susceptible.  Then, conditional upon a major outbreak,
\[
\Tnbar \convp \tau \quad \mbox{as } n \to \infty,
\]
where $\tau$ is the unique solution in $(0,1)$ of $F_0(x)=0$.
\end{enumerate}
\end{theorem}

For $R_0>1$, let $\tau_{SI}(\mu, \lambda, \omega, \alpha)$ denote the  solution of $F_0(x)=0$ in $(0,1)$, i.e.~the fraction of the population that is infected by a major outbreak in the limit as $n \to \infty$, assuming $1$ initial infective.  The fact that $\tau_{SI}(\mu, \lambda, \omega, \alpha)$ satisfies an explicit equation facilitates a detailed analysis of $\tau_{SI}(\mu, \lambda, \omega, \alpha)$, which is given in Corollary~\ref{cor:SIfinal} below.  We first highlight a key result that is an immediate consequence of Corollary~\ref{cor:SIfinal}(b).  Note that setting $\gamma=0$ in~\eqref{equ:Rzero} and~\eqref{equ:lambdac} yields $R_0=\frac{\mu\lambda}{\lambda+\omega}$ and $\lambda_c=\frac{\omega}{\mu-1}$. Let $\tau_{SI}(\mu, \lambda, \omega, \alpha)=0$ if $R_0 \le 1$.
\begin{theorem}
\label{thm:SIdisc}
Suppose that $\mu>1, \omega$ and $\alpha$ are held fixed. Then $\tau_{SI}(\mu, \lambda, \omega, \alpha)$ is discontinuous at the
threshold $\lambda=\lambda_c$ if and only if $\alpha > \frac{1}{3}$ and $\mu>\frac{3\alpha}{3\alpha-1}$.
\end{theorem}

Before stating Corrolary~\ref{cor:SIfinal} we need some more notation.  For $(\mu, \alpha)\in (1,\infty) \times [0,1]$, let $\theta(\mu, \alpha)=\frac{2\alpha(\mu-1)}{\mu+\alpha(\mu-1)}$ and
\begin{equation}
\label{equ:f0def}
f_0(x)=\log(1-x)+\frac{x}{1-\theta(\mu, \alpha)x} \qquad (0 < x < 1).
\end{equation}
Let
\begin{equation}
\label{equ:h}
h(\mu, \alpha)=\log \left(\frac{\alpha-(1-\alpha)\mu}{2\alpha \mu}\right)+\frac{\mu+\alpha(\mu-1)}{2\alpha}\qquad (\frac{1}{2} < \alpha \le 1, 1<\mu<\frac{\alpha}{1-\alpha})
\end{equation}
and, for $\alpha \ge \frac{7}{9}$, let
\[
\hat{\mu}(\alpha)=\frac{\alpha}{2(1-\alpha)}\left(1+\sqrt{\frac{9\alpha-7}{1+\alpha}}\right).
\]
Let
\begin{equation}
\label{equ:x0}
x_0(\mu,\alpha)=\frac{1+\alpha}{2\alpha}-\frac{1}{2 \mu}.
\end{equation}

The following facts, required for the statement of Corollary~\ref{cor:SIfinal} below, are proved in Section~\ref{sec:proofcorSIfinal}.   There exists a unique $\alpha^* \in (\frac{7}{9},1)$ such that
$h(\hat{\mu}(\alpha),\alpha)=0$. Let $\tau^*=x_0(\alpha^*,\hat{\mu}(\alpha^*))$.   (Numerical calculation yields $\alpha^* \approx 0.8209$, $\hat{\mu}(\alpha^*) \approx 3.3482$ and $\tau^*\approx 0.8764$.)  For fixed $\alpha \in (\alpha^*,1)$, the equation
$h(\mu, \alpha)=0$ has two solutions for $\mu \in (1,\frac{\alpha}{1-\alpha})$, which we denote by
$\mu_L^*(\alpha)$ and $\mu_U^*(\alpha)$, where $\mu_L^*(\alpha)<\mu_U^*(\alpha)$.  Finally,
$h(\mu, 1)=0$ has a unique solution $\hat{\mu}^*(1)\approx 1.7564$ in $(1, \infty)$, which satisfies $2\mu=\re^{\mu-\frac{1}{2}}$.

\begin{cor}
\label{cor:SIfinal}
\begin{enumerate}
\item[(a)] For fixed $\mu >1, \omega > 0$ and $\alpha \in [0,1]$,
\[
\lim_{\lambda \downarrow \lambda_c} \tau_{SI}(\mu, \lambda, \omega, \alpha)= \tau_0(\mu,\alpha),
\]
where $\tau_0(\mu,\alpha)$ is the largest solution in $[0, 1)$ of $f_0(x)=0$.  Note that the limit
$\tau_0(\mu,\alpha)$ is independent of $\omega$.

\item[(b)] If $\alpha \le \frac{1}{3}$ then $\tau_0(\mu,\alpha)=0$ for all $\mu > 1$.  If
$\alpha >\frac{1}{3}$ then
\[
\tau_0(\mu,\alpha)>0 \mbox{ if and only if } \mu > \frac{3 \alpha}{3 \alpha -1}.
\]

\item[(c)] For fixed $\omega > 0$ and $\alpha \in [0,1]$ we have the following, where implicitly $\lambda>\lambda_c$.
\begin{enumerate}
\item[(i)]
If $\alpha < \alpha^*$ then $\tau_{SI}(\mu, \lambda, \omega, \alpha)$ is strictly increasing in $\lambda$ for all $\mu$.
\item[(ii)] If $\alpha = \alpha^*$ then $\tau_{SI}(\mu, \lambda, \omega, \alpha)$ is strictly
    increasing in $\lambda$, unless $\mu=\hat{\mu}(\alpha^*)$ when $\tau_{SI}(\mu, \lambda, \omega, \alpha)=\tau^*$ for all $\lambda>\lambda_c$.
\item[(iii)] If $\alpha^* < \alpha < 1$ then $\tau_{SI}(\mu, \lambda, \omega, \alpha)$ is strictly increasing in
$\lambda$ if $\mu<\mu_L^*(\alpha)$ or $\mu>\mu_U^*(\alpha)$, strictly decreasing in $\lambda$ if
$\mu_L^*(\alpha)<\mu<\mu_U^*(\alpha)$, and independent of $\lambda$ if $\mu=\mu_L^*(\alpha)$ or $\mu_U^*(\alpha)$; $\tau_{SI}(\mu, \lambda, \omega, \alpha)=x_0(\mu_L^*(\alpha), \alpha)$ for all $\lambda$ if $\mu=\mu_L^*(\alpha)$ and $x_0(\mu_U^*(\alpha), \alpha)$ for all $\lambda$ if $\mu=\mu_U^*(\alpha)$.
\item[(iv)] If $\alpha=1$ then $\tau_{SI}(\mu, \lambda, \omega, \alpha)$ is strictly increasing in $\lambda$ if
$\mu<\hat{\mu}^*(1)$, independent of $\lambda$ if $\mu=\hat{\mu}^*(1)$ and strictly decreasing in $\lambda$ if
$\mu>\hat{\mu}^*(1)$.
\end{enumerate}
\end{enumerate}

\end{cor}

\begin{rmk}
\label{rmk:SIRSI}
Assuming that the conjecture in Remark~\ref{rmk:disc} is correct, as $\gamma \downarrow 0$ in the SIR model the final size is discontinuous at $\lambda_c$ if and only if $\alpha>\frac{1}{2}$ and $\mu> \frac{2\alpha}{2\alpha-1}$,
which is different from that for the $SI$ model, viz.~$\alpha>\frac{1}{3}$ and $\mu> \frac{3\alpha}{3\alpha-1}$  (see Corollary~\ref{cor:SIfinal}(b)).
\end{rmk}

\begin{rmk}
\label{rmk:varyomega}
By linearly rescaling time it is immediate that $\tau_{SI}(\mu, \lambda, \omega, \alpha)$ depends on
$\lambda$ and $\omega$ only through $\lambda/\omega$.  For fixed $\lambda>0$, let $\omega_c=(\mu-1)\lambda$, so, in the limit $n \to \infty$, a major outbreak is possible if and only if
$\omega \in [0, \omega_c)$.  It follows immediately from Corollary~\ref{cor:SIfinal}(a) that, for fixed
$\mu >1, \lambda \ge 0$ and $\alpha \in [0,1]$,
\[
\lim_{\omega \uparrow \omega_c}=\tau_{SI}(\mu, \lambda, \omega, \alpha)= \tau_0(\mu,\alpha).
\]
Hence, by Corollary~\ref{cor:SIfinal}(b), $\tau_{SI}(\mu, \lambda, \omega, \alpha)$ is discontinuous at the threshold $\omega_c$ if and only if $\alpha>\frac{1}{3}$ and $\mu > \frac{3 \alpha}{3 \alpha -1}$.
The obvious analogue of  Corollary~\ref{cor:SIfinal}(c) holds.
\end{rmk}

\begin{rmk}
\label{rmk:giant}
Note that if $\omega=0$ and $\lambda>0$ then the epidemic ultimately spreads to all individuals in the components of the underlying graph $G(n, \mu/n)$ that have initial infectives,
so the distribution of the final size $\Tn$ is independent of $\lambda$.  In particular, for $\mu>1$ and $\lambda>\mu^{-1}$ (so $R_0>1$), $\rho(\mu)=\tau_{SI}(\mu, \lambda, 0, 0)$ gives the fraction of individuals in the giant component of $G(n, \mu/n)$ in the limit as $n \to \infty$.  Setting $\omega=0$ in Theorem~\ref{thm:SIfinalsize}(b) yields that $\rho(\mu)$ is the unique solution in $(0,1)$ of $1-x=\re^{-\mu x}$ (see, for example, \cite{Durrett:2007}, Theorem 2.3.2).
\end{rmk}

\begin{figure}
\begin{center}
\resizebox{0.6\textwidth}{!}{\includegraphics{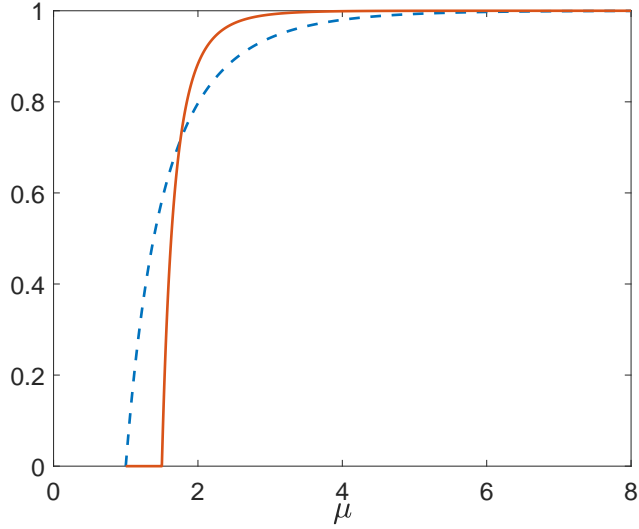}}
\end{center}
\caption{Plots of $\tau_0(\mu,1)$ (solid curve)and $\rho(\mu)$ (dashed curve).}
\label{fig:tau0GC}
\end{figure}

\begin{rmk}
Suppose that $\omega>0$ and $\alpha=1$, so there is no dropping of edges.
Plots of $\tau_0(\mu,1)$, which gives the size of the discontinuity in $\tau$ at the threshold $\lambda_c$, and $\rho(\mu)$ are shown in Figure~\ref{fig:tau0GC}.  Note that $\tau_0(\mu,1)< \rho(\mu)$ for $\mu < \hat{\mu}^*(1)$  and $\tau_0(\mu,1)> \rho(\mu)$ for $\mu > \hat{\mu}^*(1)$, which is consistent with Corollary~\ref{cor:SIfinal}(c)(iv).  Thus, when $\alpha=1$, rewiring reduces the size of a major epidemic if $\mu < \hat{\mu}^*(1)$ and increases it if  $\mu > \hat{\mu}^*(1)$.

It is interesting to note that~\cite{BT:2012} find the identical threshold $\hat{\mu}^*(1)\approx 1.7564$ when studying the size of the giant component of a Poissonian random graph where nodes have mixed-Poisson degrees with mean $\mu$.  They found that whenever $\mu > \hat{\mu}^*(1)$ the giant is maximized when the Poisson degree is not mixed and instead is ${\rm Po}(\mu)$. The giant component is hence \emph{maximized} with minimal heterogeneity in degree distribution whereas the present result shows that the fraction infected by a major outbreak is \emph{minimized} when there is no rewiring whenever $\mu > \hat{\mu}^*(1)$.
\end{rmk}

\subsubsection{SIR model with rewiring only to susceptibles}
\label{sec:SIRsusonly}

Suppose now that when a susceptible individual rewires an edge away from an infective, instead of rewiring to an individual chosen uniformly at random from the other $n-2$ individuals in the population it rewires
to an individual chosen uniformly at random from the other \emph{susceptible} individuals in the population.  If
there is no other susceptible individual then it does not rewire.  The deterministic approximation becomes
\begin{align}
\dfrac{ds}{dt}&=-\lambda i_E, \label{equ:dsdt1}\\
\dfrac{di}{dt}&=-\gamma i + \lambda i_E\label{equ:didt1}\\
\dfrac{di_E}{dt}&= -\lambda i_E -\gamma i_E +\lambda \mu i_E s -\lambda \frac{i_E^2}{s} +2 \lambda i_E \frac{w}{s}-\omega i_E,\label{equ:diEdt1}\\
\dfrac{dw}{dt}&= \omega \alpha i_E -2 \lambda i_E \frac{w}{s}.\label{equ:dwdt1}
\end{align}
Note that the only changes from the ordinary differential equations~\eqref{equ:dsdt}-\eqref{equ:dwdt}
are in the last term in~\eqref{equ:diEdt1}, where the factor $1-\alpha+\alpha(1-i)$ has been dropped as now rewiring necessarily leads to the loss of one infectious edge, and in the first term in~\eqref{equ:dwdt1}, where the factor $s$ has been dropped as now rewiring necessarily leads to the
gain of one susceptible-susceptible rewired edge.  (The difference when there is only one susceptible remaining does not affect the deterministic limit.)

\begin{theorem}
\label{thm:temporalWLLNsusonly}
Suppose $\bXnbar(0) \convp x(0)$ as $n \to \infty$, where $i(0)>0$ and $i_E(0)>0$.  Then, for any $t_0>0$,
\[
\sup_{0 \le t \le t_0} \left|\bXnbar(t)-x(t)\right| \convp 0 \quad \mbox{as } n \to \infty,
\]
where $x(t)=(s(t), i(t), i_E(t), w(t))$ is given by the solution of~\eqref{equ:dsdt1}-\eqref{equ:dwdt1}
with initial state $x(0)$.
\end{theorem}

Note that the differential equations~\eqref{equ:dsdt1}, \eqref{equ:diEdt1} and~\eqref{equ:dwdt1} form a closed system, whose time-transformed version is
Lipschitz (see Section~\ref{sec:proofsusonly}).  Consequently the final outcome of the deterministic model is readily susceptible to analysis and corresponding laws of large numbers are readily available for the stochastic model.
For brevity, we restrict attention to the case of the stochastic model with one initial infective.  Theorem~\ref{thm:BPapprox} holds also for the present model, see Remark~\ref{rmk:susonly} at the end of Section~\ref{sec:proofBP}, so $R_0$ and $\lambda_c$ are still given by~\eqref{equ:Rzero} and~\eqref{equ:lambdac}, respectively.   Let
\[
g(x)=\left(1+\frac{\gamma+\omega(1-2\alpha)}{\lambda}\right)\log(1-x)
+\left(\mu-\frac{2\alpha\omega}{\lambda}\right)x \quad (0 \le x<1)
\]
and 
\begin{equation}
\label{equ:rdef}
r(\mu, \gamma, \omega, \alpha)=\mu(\gamma+\omega-2\alpha\omega)+2\alpha\omega.
\end{equation}

\begin{theorem}
\label{thm:SIfinalsizesusonly}

Suppose that $R_0>1$ and for each $n$ the epidemic is started by $1$ infective, with the
rest of the population susceptible.  Then, conditional upon a major outbreak,
\[
\Tnbar \convp \tilde{\tau}=\tilde{\tau}(\mu, \lambda, \gamma, \omega, \alpha) \quad \mbox{as } n \to \infty,
\]
where
\begin{enumerate}
\item[(a)]
if $r(\mu, \gamma, \omega, \alpha)\ge 0$ then, for all $\lambda > \lambda_c$, $\tilde{\tau}$ is
given by the unique solution in $(0,1)$ of g(x)=0;
\item[(b)]
if $r(\mu, \gamma, \omega,\alpha)< 0$ then $\tilde{\tau}=1$, for $\lambda_c <\lambda \le \omega(2\alpha-1)-\gamma$, and $\tilde{\tau}$ is given by the unique solution in $(0,1)$ of $g(x)=0$, for $\lambda> \omega(2\alpha-1)-\gamma$.
\end{enumerate}
\end{theorem}

\begin{rmk}
Note that a necessary condition for $r(\mu, \gamma, \omega, \alpha)< 0$ is $\alpha>\frac{\gamma+\omega}{2\omega}$.
Moreover, if this condition is satisfied then, in the limit $n \to \infty$, the fraction of the population that is infected by a major outbreak is one!
\end{rmk}

\begin{rmk}
Note from Theorem~\ref{thm:SIfinalsizesusonly} that if $\lambda$ is increased, with all other parameters held fixed, then for all sufficiently large $\lambda$ the final size of the epidemic is strictly less than one and is given by the root of $g$ in $(0,1)$.  Thus, in the case of Theorem~\ref{thm:SIfinalsizesusonly}(b) (i.e.~$r(\mu,\gamma, \omega, \alpha)<0$), as $\lambda$ is increased from $0$, the final size jumps from $0$ to $1$ at the threshold $\lambda_c$, stays at $1$ until $\lambda=\omega(2\alpha-1)-\gamma$ and then decreases.  This is because when $\lambda$ is very large susceptible individuals with infectious neighbours become infected before they are able to rewire away from them. Indeed, in the limit $\lambda \to \infty$, all individuals in components of the graph that contain initial infectives become infected.  In particular, the final size of a major outbreak coincides with
that of the giant component, as is easily verified by letting $ \lambda \to \infty$ in the equation $g(x)=0$; cf.~Remark~\ref{rmk:giant}.
\end{rmk}

\begin{figure}
\begin{center}
\resizebox{0.49\textwidth}{!}{\includegraphics{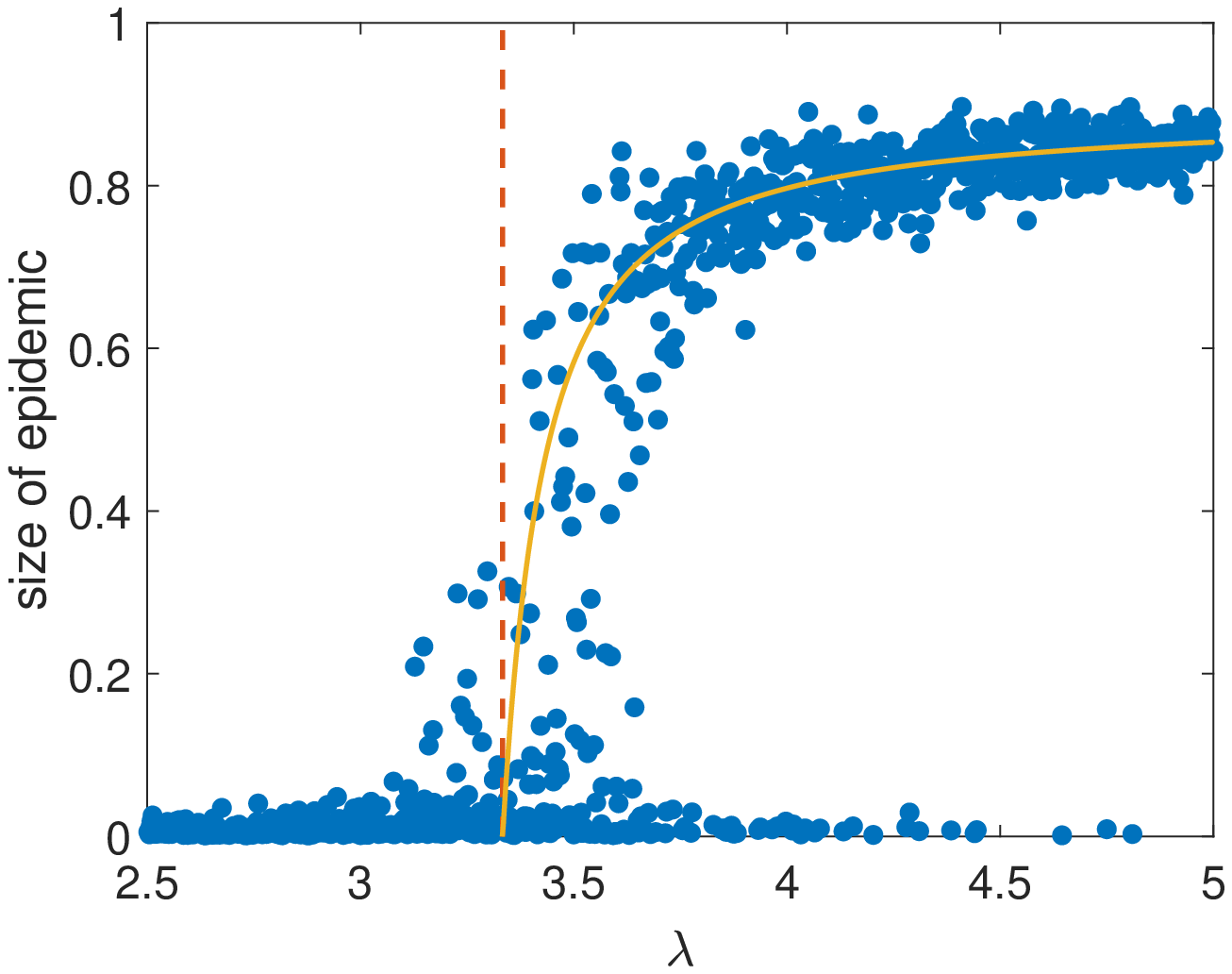}}
\resizebox{0.49\textwidth}{!}{\includegraphics{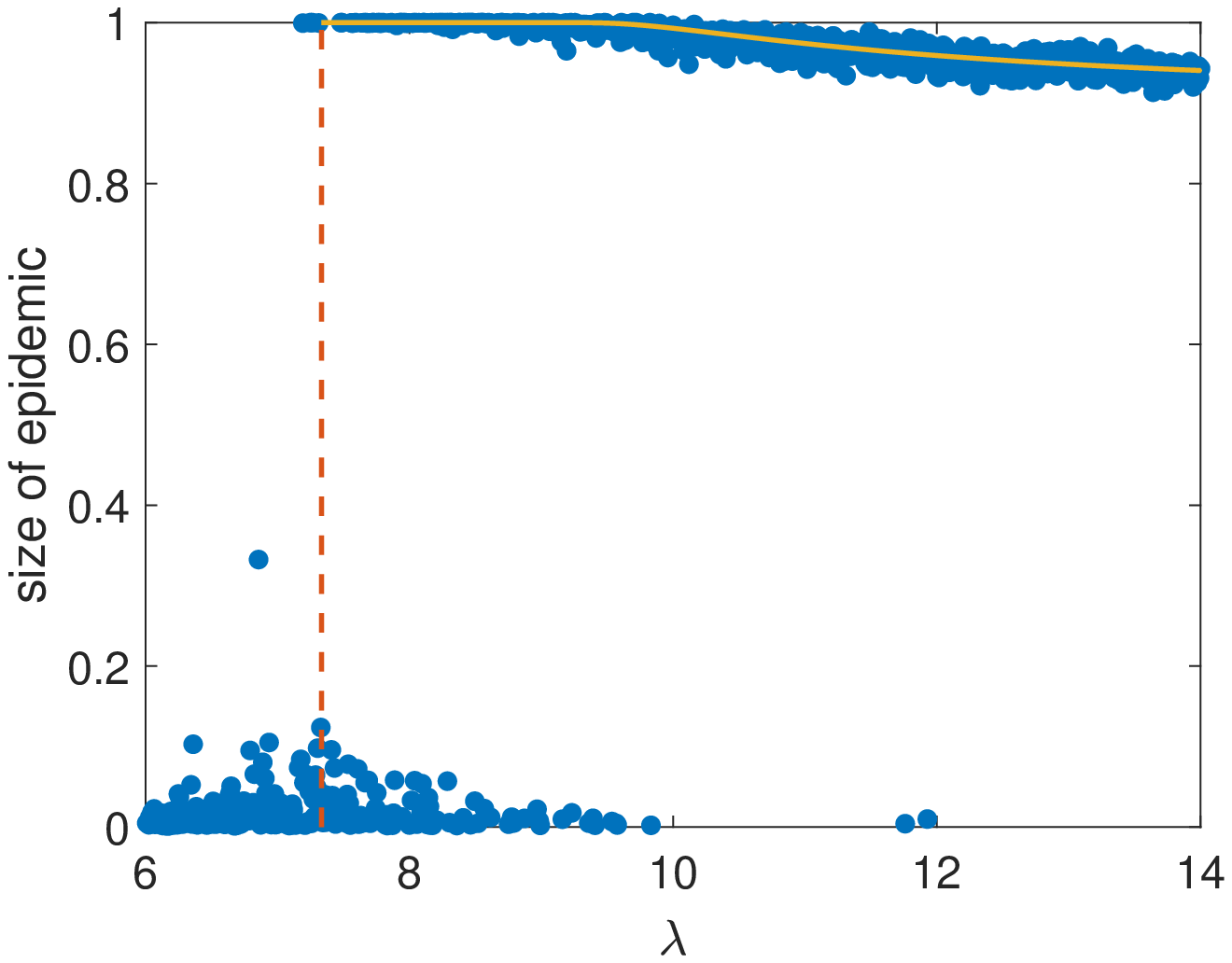}}
\end{center}
\caption{1,000 simulations of final size of SIR epidemic with rewiring only to susceptibles when $n=10,000, \mu=2.5, \gamma=1, \alpha=1$ and varying $\lambda$; $\omega=4$ in the left panel and $\omega=10$ in the right panel. See text for details.}
\label{fig:susonly}
\end{figure}

Figure~\ref{fig:susonly} shows the results of simulations of the SIR model with rewiring
only to susceptibles.  In each of the two plots the final size of $1,000$ epidemics in a population of size $n=10,000$, with values of $\lambda$ equally spaced in the given range, $\gamma=1$, $\alpha=1$ (so there is no dropping of edges) and $\omega=4$ (left panel) and $\omega=10$ (right panel).  Each simulation was started with $10$ infectives, chosen uniformly at random from the population.  The two values of $\omega$ were chosen so that $r(\mu, \gamma, \omega, \alpha)< 0$ (left panel) and $r(\mu, \gamma, \omega, \alpha)> 0$ (right panel).  The vertical dashed line shows the critical value $\lambda_c$ of $\lambda$.
The solid curve shows the final size given by Theorem~\ref{thm:SIfinalsizesusonly}(a).
The figure demonstrates clearly Theorem~\ref{thm:SIfinalsizesusonly}.  Note in the right panel that the final size is less than $1$ for $\lambda>\omega-\gamma=9$.  There are a few simulations in the right panel in which the whole population is infected when $\lambda<\lambda_c$.  This is a finite-population effect, which would be reduced if $n$ was increased and enhanced if $n$ was decreased.

\subsubsection{SIR models with alternative rewiring}
\label{sec:altrewire}

Two other rewiring models are when a susceptible rewires from an infective it (i) rewires to a randomly chosen \emph{recovered} individual or (ii) it rewires to an individual chosen uniformly at random from \emph{all} those individuals that are \emph{not} currently infectious.  Note that in model (i) the edge is effectively dropped as far as disease transmission is concerned, so the model is equivalent to the dropping model (i.e.~the corresponding model with $\alpha=0$).
The analysis of model (ii) is similar to, and has the same difficulties as, that of the SIR model considered in Section~\ref{sec:SIR}.  Models (i) and (ii) each have $R_0$ given by~\eqref{equ:Rzero} and critical infection rate $\lambda_c$ given by~\eqref{equ:lambdac}.  We omit the details but
Theorem~\ref{thm:detdisc} holds also for model (ii) and we conjecture that the conjecture in Remark~\ref{rmk:disc} holds as well.

\subsection{Construction}
\label{sec:construction}

The main idea of the construction is to explicitly consider only infectious individuals and connections emanating from these individuals and to keep track also of the pool of rewired edge-stubs. Only when a new individual is infected are its degree and whether it has received any rewired edges determined.

To help introduce the construction, consider first the SIR model without rewiring (or dropping).  We construct a realisation of the epidemic and
Erd\H{o}s-R\'{e}nyi random graph simultaneously as follows (cf.~\cite{BO'N:2002} and~\cite{Neal:2003}).  We start with $I_0$ infectives, $n-I_0$ susceptibles and no edges in the graph.
Attach to each infective independently ${\rm Po}(\mu_n)$ edges, where $\mu_n=-n \log(1-\frac{\mu}{n})$ so $\mu_n$ satisfies
$1-\frac{\mu}{n}={\rm e}^{-\frac{\mu_n}{n}}$. Such unconnected edges are called live infectious edges. Infectives make infectious contacts along each of their live infectious edges independently
at rate $\lambda$.  When an infective, individual $i_0$ say, makes infectious contact along an edge, that edge emanating from $i_0$ becomes dead and an individual, $j_0$ say, is sampled independently and uniformly from the whole population.  If individual $j_0$ is susceptible then an edge between $i_0$ and $j_0$ is formed, individual $j_0$ becomes infected and (independently) has ${\rm Po}(\mu_n)$ live infectious edges attached to it.  If individual $j_0$ is not susceptible then
nothing happens (apart from the edge from $i_0$ becoming dead).  Infectives recover independently at rate $\gamma$.  When an infective recovers, all live edges attached to it become dead.  The epidemic stops as soon as there is no infective remaining. (Infection spread stops as soon as there is no live infectious edge remaining.  In the SI model the epidemic stops when there is no live infectious edge remaining.) It is easily verified that this construction gives an epidemic that is probabilistically equivalent to the SIR epidemic on an Erd\H{o}s-R\'{e}nyi graph. (Note that $\mu_n$ is chosen so that if
$\gamma=0$, so an infective infects down \emph{all} of the edges that are initially attached  to it, then it tries to infect the other $n-1$ individuals independently,
each with probability $\frac{\mu}{n}$.)

An alternative approach, which is more easily generalised to incorporate rewiring, is to modify the above construction so that only live infectious edges that will be paired with susceptibles are included,
and each time an infection occurs along such an edge it is with an individual chosen independently and uniformly
from all susceptible individuals at that time.  More specifically, suppose an individual $i_0$ say, is infected at time $t$.  Then ${\rm Po}(\mu_n \frac{\Sn(t)}{n})$ live infectious edges are attached to it.  However, we now have to delete each of the live infectious edges at time $t-$
(apart from the one that led to the infection of $i_0$, which is necessarily deleted) independently with probability $\frac{1}{\Sn(t)}$, to take account of the fact that all live infectious edges are necessarily to susceptible individuals.

Consider now the SIR model $\EEn$ \emph{with} rewiring and dropping of edges described in Section~\ref{sec:modelmain}.  We adopt the above alternative approach and treat rewiring/dropping as follows.  Let $\Wn(t)$ be the number of
rewired susceptible-susceptible edges at time $t$, so $\Wn(0)=0$.  Each infective independently sends warnings
down each of its unconnected live infectious edges independently at rate $\omega$.  Suppose that an individual, $i_0$ say, sends
such a warning at time $t$.  Then $i_0$ drops that edge with probability $1-\alpha$; otherwise the edge is rewired to a susceptible individual, an infective individual or a recovered individual with probabilities $(\Sn(t)-1)/(n-2)$, $(\In(t)-1)/(n-2)$ and $1-(\Sn(t)+\In(t)-2)/(n-2)$,
respectively, where $\In(t)$ is the number of infectives at time $t$ (a rewiring is always to a new individual other than the individual itself).
If the edge is rewired to a susceptible then $\Wn(t)$ is increased by $1$.  (We do not decide the two susceptible
individuals involved at this stage.)  If the edge is rewired to an infective, then $i_0$ loses an infectious
edge and an individual chosen uniformly at random from the other $\In(t)-1$ infectives at time $t$ gains
an infectious edge.  If the edge is rewired to a recovered individual then $i_0$ loses an infective
edge but nothing else happens.

As before infectives infect along their live infectious edges at rate $\lambda$.
If an infective, $i_0$ say, infects down a live infectious edge at time $t$.  Then that live edge is connected with an individual, $j_0$ say, chosen uniformly from $\SSn(t)$, where, for $t \ge 0$, $\SSn(t)$ is the set comprising of the $\Sn(t)$ susceptibles at time $t$. The live edge then becomes dead. Individual $j_0$ is infected
at time $t$ and acquires ${\rm Po}(\mu_n \frac{\Sn(t)}{n})$ live infectious edges (corresponding to edges in the original network), the fraction $ \Sn(t)/n$ comes from now only considering edges to be connected to susceptibles.  Individual $j_0$ may also acquire rewired infectious edges as follows.  For each $i=1,2,\dots,\Wn(t- )$,
independently sample two individuals, $k_i$ and $l_i$ say, without replacement from $\SSn(t-)$. If $j_0 \in \{k_i, l_i\}$ then $j_0$ acquires a further live infectious edge and $\Wn(t)$ is reduced by one.  Finally, each of the live infectious
edges at time $t-$, apart from the one that infected $j_0$, is dropped independently with probability $1/\Sn(t-)$ (thus reducing the number of live edges to be connected to susceptibles correctly). The live infectious edge that infected $j_0$ is necessarily dropped. The construction now continues in the obvious manner.
The epidemic stops when there is no infective remaining (SIR model) or no infectious edge remaining (SI model).

The above construction is not fully faithful to the model $\EEn$, as it allows the possibility
of a susceptible, $i_0$ say, to rewire from an infective, $j_0$ say, and then be subsequently infected by
$j_0$ without having first rewired back to $j_0$.  Let $F_n$ be the probability that no such imperfection occurs.  We show at the start of Section~\ref{sec:proofs} that $\liminf_{n \to \infty} \P(F_n)>0$, and consequently
any statement which holds with high probability for the construction, holds also with high probability when the construction is conditioned on being faithful to the original model (cf.~\cite{Janson:2009}).

For $t \ge 0$, let $\Ine(t)$ be the total number of live infectious edges at time $t$.  We show in Section~\ref{thm:temporalWLLN} that the $\{(\Sn(t), \In(t), \Ine(t), \Wn(t)):t \ge 0\}$ converges to the deterministic model given by~\eqref{equ:dsdt}-\eqref{equ:dwdt} as $n \to \infty$, in
the sense made clear in Theorem~\ref{thm:temporalWLLN}.  The terms in~\eqref{equ:diEdt} are explained as follows.  The first two terms correspond to a live infectious edge being lost if the associated infective infects along it or recovers.  The third term is owing to the creation of new live infectious edges (corresponding to edges in the original network) when
an infective transmits infection along an infectious edge; on average each such infection creates $\mu s$ new live infection edges as only edges that will be paired with a susceptible are included.  The fourth term arises from live infectious edges that are lost when a susceptible is infected.  The fifth term
comes from the gain of live infectious edges from the pool of rewired susceptible-susceptible edges when a susceptible is infected; the factor $2$ is present as each such rewired susceptible-susceptible edge has probability $2/\Sn(t-)$ of having an end at the newly infected susceptible.  The final term comes from rewiring/dropping of live infectious edges; note that any edge down which an infective sends a warning is lost unless the edge is rewired to another infective.  The terms in the other $3$ equations are explained similarly.

\section{Pair approximation}
\label{sec:pairapprox}
In this section we describe a deterministic pair-approximation model for the SIR epidemic of Section~\ref{sec:SIR} (see~\cite{gross06} for a similar pair approximation to the corresponding SIS epidemic) and show that it is equivalent to
the deterministic model given by~\eqref{equ:dsdt}-\eqref{equ:dwdt}.  It then follows using Theorem~\ref{thm:temporalWLLN} that the pair-approximation SIR
model is exact in the limit as $n \to \infty$.  For simplicity we treat the case $\alpha=1$, i.e.~no dropping of edges.  We use the pair-approximation paradigm described in~\cite{AB00}, Chapter 7.

Consider the stochastic SIR model with rewiring described in Section~\ref{sec:model}.  Label the nodes (individuals) in the population $1,2,\dots,n$.  Let $G=[g_{ij}]$ be the $n \times n$ matrix having elements $g_{ij}$, where $g_{ij}$ is the number of edges joining nodes $i$ and $j$, if $i \neq j$, and $g_{ii}=0$.  (In many applications $G$ is an adjacency matrix but here the elements of $G$ can be $>1$ owing to rewiring.)  Each node is of type $S, I$ or $R$, according to its disease status.  For $A \in \{S, I, R\}$ and $i=1,2,\dots, n$, let $A_i=1$ if node $i$ is of type $A$ and $A_i=0$ otherwise.
For $A , B, C \in \{S, I, R\}$, let
\[
[A]_n=\sum_{i=1}^nA_i, [AB]_n=\sum_{i,j=1}^nA_i g_{ij} B_j\quad \text{and}\quad
[ABC]_n=\sum_{\substack{i,j,k=1\\i \ne k}}^nA_i g_{ij} B_j g_{jk}C_k.
\]
Note that $[A]_n$, $[AB]_n$ and $[ABC]_n$ are all functions of time $t$, as indeed is $G$ owing to rewiring.  For example, $[I]_n(t)$ is the number of infectives at time $t$ and $[SI]_n(t)$ is the number of edges between infective and susceptible nodes at time $t$.  Note that $[SS]_n(t)$ is \emph{twice}
the number of edges between susceptible nodes at time $t$.

Suppose that $n^{-1}[A]_n(t), n^{-1}[AB]_n(t)$ and $n^{-1}[ABC]_n(t)$ each tend to deterministic limits as $n \to \infty$. For example, $n^{-1}[S]_n(t) \to s(t)$ and $n^{-1}[SI]_n(t)\to x_{SI}(t)$.  A heuristic argument then yields the following system of ordinary differential equations:
\begin{align}
\dfrac{ds}{dt} &= -\lambda x_{SI}, \label{equ:pairs}\\
\dfrac{di}{dt}&= \lambda x_{SI} -\gamma i,\label{equ:pairi}\\
\dfrac{dx_{SS}}{dt} &= 2\omega s x_{SI} -2\lambda x_{SSI},\label{equ:pairss}\\
\dfrac{dx_{SI}}{dt} &= -\lambda x_{SI} + \lambda x_{SSI} - \lambda x_{ISI} -[\omega(s+r)+\gamma ]x_{SI}.\label{equ:pairsi}
\end{align}
The second and third terms in~\eqref{equ:pairsi} are explained as follows.  The second term is added to the first since an infection along an
$SSI$ link does not change the number of $SI$ links.  The third term is present because not only is the rate doubled, although this is already counted since it consists of two $SI$ links, but when an infection occurs the number of $SI$ links is reduced by $2$.  The factor $2$ in~\eqref{equ:pairss} arises from the double counting in $[SS]_n(t)$.  There are other variables, such as $x_{II}$, for which we could write down a differential equation but the above system is sufficient for our purposes.

To close the above system we make the pair approximation
\[
x_{ABC}=\frac{x_{AB} x_{BC}}{x_B} \qquad \mbox{for all } A, B, C \in  \{S, I, R\},
\]
where, for example, $x_S=s$.  Equations~\eqref{equ:pairss} and~\eqref{equ:pairsi} then become
\begin{align}
\dfrac{dx_{SS}}{dt} &= 2\omega s x_{SI} -2\lambda \frac{x_{SI}x_{SS}}{s},\label{equ:pairss1}\\
\dfrac{dx_{SI}}{dt} &= -\lambda x_{SI}\left(1-\frac{x_{SS}}{s}+\frac{x_{SI}}{s}\right) -\left(\omega(s+r)+\gamma\right)x_{SI}.\label{equ:pairsi1}
\end{align}

To connect informally with the deterministic model given by~\eqref{equ:dsdt}-\eqref{equ:dwdt}, note that $[SS]_n$ includes $SS$ links in the original
network and $SS$ links arising from rewiring.  Taking account of the double counting in $[SS]_n$, it follows that
\begin{equation}
\label{equ:wxssconnect}
w(t)=\frac{1}{2}\left(x_{SS}(t)-\mu s(t)^2\right)\qquad\mbox{for all }t \ge 0.
\end{equation}
For a formal proof, substitute from~\eqref{equ:wxssconnect} into~\eqref{equ:pairsi1} to obtain
\[
\dfrac{dx_{SI}}{dt} = -\lambda x_{SI}\left(1-2\frac{w}{s}-\mu s+\frac{x_{SI}}{s}\right) -\left(\omega(s+r)+\gamma\right)x_{SI},
\]
which, noting that $s+r=1-i$, coincides with~\eqref{equ:diEdt}, with $\alpha=1$ and $i_E$ replaced by $x_{SI}$.  Also, differentiating~\eqref{equ:wxssconnect} and using~\eqref{equ:pairs}, \eqref{equ:pairss1} and~\eqref{equ:wxssconnect}, yields
\begin{align*}
\dfrac{dw}{dt}&=\frac{1}{2}\left(\dfrac{dx_{SS}}{dt}-2\mu s \dfrac{ds}{dt} \right)\\
&=\frac{1}{2}\left(2\omega s x_{SI} -2\lambda \frac{x_{SI}x_{SS}}{s}+2\mu s \lambda x_{SI}\right)\\
&=\omega s x_{SI}-2\lambda \frac{x_{SI}w}{s},
\end{align*}
which agrees with~\eqref{equ:dwdt}.  Finally, \eqref{equ:pairs} and~\eqref{equ:pairi} clearly agree with~\eqref{equ:dsdt} and~\eqref{equ:didt}.

\section{Comparison with~\cite{YD:20}}
\label{sec:YaoDurrett}
\cite{YD:20} consider SI and SIR epidemics with rewiring on graphs generated by the configuration model (see, for example, \cite{Durrett:2007}, Chapter 3, and~\cite{vdHofstad:2016}, Chapter 7) or its close relatives, which include Erd\H{o}s-R\'{e}nyi~random graphs.  Briefly, the configuration model random graph on $n$ individuals, labelled $1,2,\dots, n$, is constructed as follows.  Let $D_1, D_2, \dots, D_n$ be i.i.d.~copies of a random variable $D$, which takes values in $\mathbb{Z}_+$.  Attach $D_i$ half-edges to individual $i$ ($i=1,2,\dots,n$) and then pair up the $D_1+D_2+\dots+D_n$ half-edges uniformly at random to form the graph.  (If $D_1+D_2+\dots+D_n$ is odd then there is a left-over half-edge, which is ignored.)  This yields the version of the model introduced in~\cite{NSW:2001}.  In an alternative version (\cite{MR:1995}) the degrees are prescribed deterministically.  The configuration model was first introduced by~\cite{Bollobas:1980}.

We outline the main result obtained by~\cite{YD:20} for the SI model.  In our notation, $\alpha=1$, so there is no dropping of edges.  Let $\mu=\E[D]$ and $G(z)=\E[z^D]$ be respectively the mean and probability-generating function of $D$.  Suppose that $\E[D^3]<\infty$.  Let $\beta=\frac{\omega\mu}{\lambda}$ ($\omega$ is denoted by $\rho$ in~\cite{YD:20}, and $\beta$ by $\alpha$).  Let $f:(0,1) \to \mathbb{R}$ be defined by
\[
f(x)=\log\left(\frac{\mu x}{G'(x)+\beta(1-x)G(x)}\right)+\frac{\beta}{2}(x-1)^2.
\]
Suppose that $\lambda>\lambda_c$ and let
\begin{align*}
\sigma&=\sup\{x: 0 < x < 1, f(x)=0\} \quad\mbox{with } \sup(\emptyset)=0,\\
v&=1-\exp\left(-\frac{\beta}{2}(\sigma-1)^2\right)G(\sigma).
\end{align*}
Suppose that either $\sigma=0$, or $\sigma \in (0,1)$ and there exists $\delta>0$ such that $f<0$ on
$(\sigma-\delta,\sigma)$. Suppose that there is initially one infective, chosen uniformly at random from the population of $n$ individuals. Then (\cite{YD:20}, Theorem 2), for any $\epsilon>0$,
\begin{equation}
\label{equ:YDSIfinal}
\lim_{n \to \infty}\P(\Tnbar<\nu+\epsilon)=\lim_{\eta \downarrow 0}\liminf_{n \to \infty}
\P(\Tnbar>\nu-\epsilon|\Tnbar>\eta)=1.
\end{equation}

\begin{figure}
\begin{center}
\resizebox{0.6\textwidth}{!}{\includegraphics{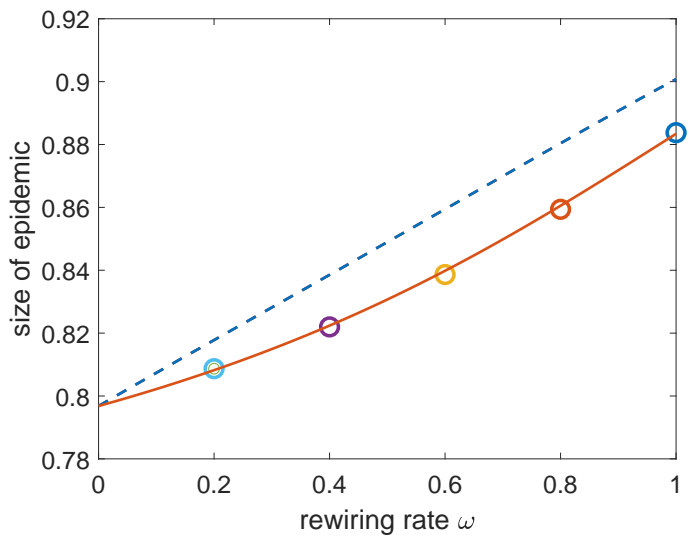}}
\end{center}
\caption{Plots of $\tau$ (solid curve) and $\nu$ (dashed curve) as functions of $\omega$ for SI epidemics with $\mu=2$, $\lambda=1$ and $\alpha=1$.  The circles show the mean sizes of $1,000$ simulated major epidemics in a population of $n=5,000$.  See text for details.}
\label{fig:YDcomp}
\end{figure}

Note that $\nu$ in~\eqref{equ:YDSIfinal} has essentially the same interpretation as $\tau$ in Theorem~\ref{thm:SIfinalsize}(b) (when $\alpha=1$), so one would expect them to be equal when $D \sim {\rm Po}(\mu)$.  Calculation shows that is not the case.  Plots of $\nu$ and $\tau$ as functions of the rewiring rate $\omega$  when $\mu=2$ and $\lambda=1$ are shown in Figure~\ref{fig:YDcomp}, together with mean fraction infected based on $1,000$ simulations of major outbreaks in a population with $n=5,000$ for each of $\omega=0.2, 0.4, 0.6, 0.8$ and $1.0$.  Each simulated epidemic was started by one infective, with the rest of the population being susceptible.  A cut-off of $3,000$ was used to determine whether a major epidemic had occurred.  (For each $\omega$, epidemics were simulated until $1,000$ major epidemics had occurred.  The smallest major epidemic had size $3,840$ and the largest non-major epidemic size $1,164$.)
Note that $R_0=1$ when $\omega=1$, so major outbreaks have zero probability of occurrence in the limit $n \to \infty$.  However, they do occur when $n=5,000$.  The value of $\tau$ when $\omega=1$ is $\tau_0(2,1)$ (see Remark~\ref{rmk:varyomega}).
The simulations are consistent with Theorem~\ref{thm:SIfinalsize}(b) and suggest strongly that~\eqref{equ:YDSIfinal} is incorrect.

The results in~\cite{YD:20} are based on an effective-degree type construction (\cite{BN:2008}) in which the random graph and epidemic are constructed simultaneously.  Briefly, half-edges are labelled susceptible, infective or recovered according to the disease status of the individual to which they are attached.  Initially no half-edge is paired.  Infectives transmit infection down their infective half-edges independently at rate $\lambda$.  When such transmission occurs, the half-edge is paired with a half-edge chosen uniformly at random from all available half-edges, to form an edge in the graph.  If the chosen half-edge is attached to a susceptible individual then that individual becomes infected and can then transmit infection down its remaining half-edges.  If it is attached to a recovered or infective individual then nothing happens, apart from the loss of the two edges.  Such an approach has been used by~\cite{JLW:2014} to prove law of large numbers results for SIR epidemics on Molloy--Reed configuration model random graphs.

In~\cite{YD:20}, rewiring is added to the construction by assuming that each infective half-edge is removed from the vertex to which it is attached independently at rate $\omega$ and is immediately re-attached to a vertex chosen uniformly at random from all other vertices.  It is claimed that the construction then yields a model whose final size has the same distribution as the original model.
However, it is not clear that is the case.  Suppose that an infective half-edge is moved from individual $i$ to individual $j$ at time $t_0$ and individual $j$ becomes infected and transmits infection down that half-edge at time $t_1>t_0$.  First, the half-edge needs to be paired to determine whether or not it should be moved; it should be moved only if it is paired with a susceptible individual.  Second, if it is moved and infection is transmitted along it at time $t_1$ then that half-edge should be paired with an individual chosen from susceptible individuals at time $t_0$, weighted by their degree at that time, and  not from all half-edges at time $t_1$.  Third, if a half-edge is moved then one of the susceptibles ($k$ say) at time $t_0$ should have a half-edge such that if $k$ is subsequently infected (and the corresponding edge has not been rewired) and transmits infection down that half-edge at time $t_2$, then that half-edge should be paired with an individual chosen (approximately) uniformly at random from the rest of the population and not with a half-edge chosen uniformly at random from the other half-edges at time $t_2$.
In view of these observations and the simulations, it would appear that the model analysed in~\cite{YD:20}
is an approximation of the epidemic with rewiring. However, it is interesting to note that in both the SI and SIR models the condition in~\cite{YD:20} for the final size to be discontinuous at the threshold $\lambda=\lambda_c$ (see their Examples 3 and 7, respectively) coincides with that obtained by setting $\alpha=1$ in Theorems~\ref{thm:detdisc} and~\ref{thm:SIdisc}.

\section{Heuristic argument for final size of SI epidemic}
\label{sec:heuristic}
Consider first the SI model with $\alpha=0$, so there is dropping of edges but no rewiring.  The fraction
infected by a major outbreak can be obtained heuristically as follows.  Suppose that a fraction $\tau$ get infected.  Fix an initial susceptible, $i$ say.  The probability that $i$ gets infected can be calculated by conditioning on the spread among all members of the population excluding $i$.  For $n$ large, in the event of a major outbreak, the number of infected neighbours of $i$ has a ${\rm Po}(\mu \tau)$ distribution, since $i$ has ${\rm Po}(\mu)$ neighbours, each of which is independently infected by the epidemic with probability $\tau$.  (These are approximations which become exact in the limit $n \to \infty$.) Each infected neighbour of $i$ infects $i$ with probability $\frac{\lambda}{\lambda+\omega}$, so the probability that $i$ avoids infection is $\exp(-\frac{\mu \tau \lambda}{\lambda+\omega})$.  This probability is also given by $1-\tau$, so $\tau$
satisfies
\[
1-\tau=\exp\left(-\frac{\mu \lambda \tau}{\lambda+\omega}\right),
\]
in agreement with Theorem~\ref{thm:SIfinalsize}(b) with $\alpha=0$.  This argument is straightforward to extend to the SIR model with $\alpha=0$.

We now extend the argument to the model with rewiring and for ease of exposition assume that $\alpha=1$, i.e.~that there is no dropping of edges.   This extension involves two new features.
\begin{enumerate}
\item[(1)] It is no longer the case that the probability that a given infective neighbour, $i^*$ say, of $i$ leads to
$i$ being infected by that neighbour is given by $\lambda/(\lambda+\omega)$ because if $i$ rewires away from $i^*$ it may be to another infective.
\item[(2)] Individual $i$ may acquire infective neighbours through rewiring, which were not neighbours of $i$ in the original network.  This happens
if a susceptible rewires to $i$ and that susceptible is subsequently infected.
\end{enumerate}
We treat these in turn, again assuming that a fraction $\tau$ of the population gets infected.

Consider an infective neighbour, $j$ say, of $i$.  Individual $i$ rewires from $j$ at rate $\omega$ and is infected by $j$ at rate
$\lambda$.  However, if $i$ rewires to an infective its position as far as the disease is concerned is unchanged.  Now $i$ rewires to a susceptible
(in which case $i$ avoids infection owing to the original link to $j$) at rate $\omega(1-\tau)$.  It follows that the probability, $p_I$ say,
that $i$ is infected owing to the original link to $j$ is given by
\begin{equation}
\label{equ:pI}
p_I=\frac{\lambda}{\lambda+\omega(1-\tau)}.
\end{equation}

We now determine heuristically the (approximate) distribution of the number of infective neighbours, $N_I^{(RW)}(i)$ say, that $i$ acquires through rewiring.
Let $\Estar$ denote the epidemic in which $i$ is excluded and assume that $n$ is large.  For any susceptible--infective edge, the probability that the susceptible rewires that edge before being infected along it is $\omega/(\lambda+\omega)$.  Thus if the epidemic $\Estar$
infects $(n-1) \tau$ individuals, by the law of large numbers there will be $\approx (n-1)\tau \omega /\lambda$ rewiring events in it.  Each such rewiring has probability $1/(n-2)$ of being to individual $i$, so if $N^{(RW)}(i)$ denotes
the number of edges that rewire to individual $i$, then $N^{(RW)}(i)\sim {\rm Po}(\tau \omega/\lambda)$ approximately.

Consider a typical edge that rewires to $i$.  For large $n$,
the number of infections in $\Estar$ that occur prior to that rewiring is distributed approximately as $\left \lceil (n-1) \tau U \right \rceil$, where $U \sim {\rm U}(0,1)$.  The number of infections
that occur after that rewiring is $\approx \left \lfloor (n-1) \tau (1-U) \right \rfloor$.  By symmetry, these $\left \lfloor (n-1) \tau (1-U) \right \rfloor$ infections are with individuals drawn uniformly at random without replacement
from the $\approx n-1-\left \lceil (n-1) \tau U \right \rceil$ suscpetibles remaining at the time of this typical rewire to $i$.  Thus, given $U$, the probability that the susceptible who rewires to $i$ is subsequently infected is approximately
\begin{equation*}
\frac{\tau(1-U)}{1-\tau U}=1-\frac{1-\tau}{1-\tau U}.
\end{equation*}
Taking expectations with respect to $U$, it follows that for a typical edge that rewires to $i$, the probability that the susceptible is subsequently infected is given approximately by
\begin{equation*}
\int_0^1 1-\frac{1-\tau}{1-\tau u} \,{\rm d}u=1+\frac{(1-\tau)}{\tau}\log(1-\tau).
\end{equation*}
Thus in the event of a major outbreak, recalling that $N^{(RW)}(i)\sim {\rm Po}(\tau \omega/\lambda)$ approximately, we have approximately for large $n$ that
\begin{equation}
\label{equ:NIRWidist}
N_I^{(RW)}(i) \sim {\rm Po}\left(\frac{\tau \omega}{\lambda}\left[1+\frac{(1-\tau)}{\tau}\log(1-\tau)\right] \right).
\end{equation}

Let $N_I(i)$ be the number of neighbours of $i$ in the original network that are infected by the epidemic.  Then, as before,
\begin{equation}
\label{equ:NIidist}
N_I(i) \sim {\rm Po}(\mu \tau),
\end{equation}
approximately.  Now $N_I(i)$ and $N_I^{(RW)}(i)$ are approximately independent, so the total number of edges from $i$ to an infective is approximately Poisson distributed with mean $\E[N_I(i)]+\E[N_I^{(RW)}(i)]$.  Each such edges leads to the infection of $i$
independently with probability $p_I$ so, using~\eqref{equ:pI}-\eqref{equ:NIidist}, $\tau=\tau_{SI}(\mu, \lambda, \omega, \alpha=1)$  satisfies
\begin{eqnarray}
1-\tau&=&\exp\left(-\left\{\mu\tau+\frac{\tau \omega}{\lambda}\left[1+\frac{(1-\tau)}{\tau}\log(1-\tau)\right]\right\}\frac{\lambda}{\lambda+\omega(1-\tau)}\right)\nonumber \\
&=&\exp\left(-\frac{\tau(\mu \lambda+\omega)+\omega(1-\tau)\log(1-\tau)}{\lambda+\omega(1-\tau)} \right). \label{equ:taur}
\end{eqnarray}

Taking logarithms of both sides of~\eqref{equ:taur} and rearranging yields
\[
\log(1-\tau)=-\frac{\tau(\mu \lambda+\omega)}{\lambda+2\omega(1-\tau)},
\]
so $F_0(\tau)=0$, in agreement with Theorem~\ref{thm:SIfinalsize}(b) when $\alpha=1$.

The arguments leading to~\eqref{equ:NIRWidist} and~\eqref{equ:NIidist} still hold for the SIR model. However, $p_I$ is no longer given by~\eqref{equ:pI}
since (i) the infective, $j$ say, at the end of an edge from $i$ may recover before it infects $i$ and (ii) if $i$ rewires from $j$ to an individual, $k$ say, that was infected by  epidemic the rewiring needs to occur \emph{before} $k$ recovers from infection.  There does not seem to be an easy way to approximate $p_I$.

\section{Proofs}
\label{sec:proofs}

We prove the results described in Section~\ref{sec:mainresults} by showing that corresponding results hold for the construction described in Section~\ref{sec:construction}.  As noted there, the construction is not fully faithful to the rewiring model described in Section~\ref{sec:model} as it allows for the possibility of a susceptible, $i_0$ say, to rewire an edge in the ``original network" from an infective, $j_0$ say, and then be subsequently infected by $j_0$ along another edge in the ``original network".  (By an edge in the ``original network" we mean one of the un-rewired infectious edges allocated to $i_0$ when it was infected.)  A realisation of a process that is probabilistically equivalent to the original model is obtained by conditioning on there being no such imperfections. The allocation of ``original network" infectious edges to infectives in the construction can be achieved by labelling the individuals in the population $1,2,\dots, n$, letting $X_1, X_2, \dots, X_n$ be
i.i.d.~${\rm Po}(\mu_n)$ random variables and $\chi_{i,j}$ $(i=1,2,\dots,n, j=1,2,\dots)$ be i.i.d. discrete
uniform random variables on $\{1,2,\dots,n\}$.  The $i$th individual to be infected is allocated $X_i$ ``original network" infectious edges, that are linked to individuals $\chi_j$ $(j=1,2,\dots, X_i)$.  Those that are not linked to susceptibles are dropped immediately.  Let $F_n$ be the event that the construction remains faithful to the model of Section~\ref{sec:model} and $D_n$ be the event that $\chi_j$ $(j=1,2,\dots, X_i)$ are distinct for each $i=1,2,\dots,n$.  Then $D_n \subseteq F_n$.  Now
\begin{align*}
\P(\chi_1, \chi_2, \dots, \chi_{X_i}\mbox{ are distinct}|X_i) &= \prod_{i=1}^{X_i-1} \left(1-\frac{i}{n}\right)\\
&\ge 1-\frac{X_i(X_i-1)}{2n},
\end{align*}
so
\begin{align*}
\P(D_n) &\ge \left(1-\frac{\E[X_1(X_1-1)]}{2n}\right)^n\\
&=\left(1-\frac{\mu_n^2}{2n}\right)^n\\
&\to \re^{-\frac{\mu^2}{2}}\quad \mbox{as } n \to \infty.
\end{align*}
Thus $\liminf_{n \to \infty} \P(F_n)>0$ and it follows that any statement which holds with high probability for the construction, holds also with high probability when the construction is conditioned on being faithful to the original model (cf.~\cite{Janson:2009}).

In the proofs that follow, we use the same notation for the construction as that defined for the original model in Section~\ref{sec:model}.  We treat the final size of the SI model before that of the SIR model as the former is simpler.

\subsection{Proofs of Theorem~\ref{thm:BPapprox} and Lemma~\ref{lem:majbound}}
\label{sec:proofBP}
Let $(\Omega, \mathcal{F}, \P)$ be a probability space on which are defined independent unit-rate Poisson processes $\eta_1, \eta_2,\dots$.  Let
$Z^E_i=\eta_i([0,\mu])$ $(i=1,2,\dots)$.  For $n=1,2,\dots$, let $\Zn_i=\eta_i([0,\mu_n(1-n^{-1}i)])$ $(i=1,2,\dots,n)$.  A realisation of the branching process
$\B$ is constructed on $(\Omega, \mathcal{F}, \P)$ in which the $i$th individual born in $\B$ is assigned at birth $Z^E_i$ infectious edges.  For $n=1,2,\dots$, a realisation of the epidemic process $\EEn$ is constructed on $(\Omega, \mathcal{F}, \P)$, using the construction in Section~\ref{sec:construction}, in which when first infected the $i$th infective has $\Zn_i$ infectious edges.  Further details of these constructions are
indicated below.  First we show that some events associated with the epidemic processes hold with high probability.

Let $\Tne=\Zn_1+\Zn_2+\dots+\Zn_{\left \lfloor{\log n}\right \rfloor}$.  Then the law of large numbers implies that $\Tne \le 2\mu\log n$ with high probability.  Let $\mathcal{I}_n$ denote the set consisting of the first $\left \lfloor{\log n}\right \rfloor$ individuals infected in $\EEn$.  (If the final size of
$\EEn$ is less than $\left \lfloor{\log n}\right \rfloor$ then $\mathcal{I}_n$ consists of all individuals infected in $\EEn$.)  Recall that $\Tn(t)=n-\Sn(t)$ is the total number of infections in $\EEn$ during $[0,t]$ and let $\hat{t}_n=\inf\{t \ge 0: \Tn(t)\ge \log n\}$, where $\hat{t}_n=\infty$ if $\Tn(\infty)< \log n$.
Let $\Wtn$ be the total number of infectious edges rewired from individuals in $\mathcal{I}_n$ during $[0, \hat{t}_n]$.  Clearly no individual in $\mathcal{I}_n$ can
have more than $\Tne$ infectious edges when it is first infected.  If every individual in $\mathcal{I}_n$ were to have every infectious edge they had when they were first infected rewired then, given $\Tne$, the expected number that would be rewired to individuals in $\mathcal{I}_n$ is less than $(\log n)\Tne \frac{\log n}{n-2}$ and it follows that, unconditionally, the number of such rewirings is zero with high probability.  Thus $\Wtn \le 2\mu\log n$ with high probability.

Recall that in $\EEn$ a newly infected individual, at time $t$ say, acquires ${\rm Bin}(\Wn(t-), 2/\Sn(t-))$ rewired infectious edges upon infection. (As usual, ${\rm Bin}(n,p)$ denotes a binomial random variable with $n$ trials and success probability $p$.)  Let
$\Wthn$ be the total number of such edges acquired by the infectives in $\mathcal{I}_n$.  Then
\[
\E\left[\Wthn|\Wtn\right]\le (\log n) \Wtn \frac{2}{n-\log n}
\]
and it follows that $\Wthn=0$ with high probability.  Recall also that in $\EEn$ when an infection occurs, say individual $i_0$ is infected at time $t$, each infectious edge in $\EEn$ at time $t$ is dropped independently with probability $1/\Sn(t-)$, apart from the one that infected $i_0$ which is necessarily dropped.  Let $\Tdn$ be the total number of such edges dropped by the infectives in $\mathcal{I}_n$.  Then,
\[
\E\left[\Tdn|\Tne, \Wthn\right] \le (\log n)\left(\Tne+\Wthn\right)\frac{1}{n-\log n},
\]
whence $\Tdn=0$ with high probability.

Collating the above, we have shown that in $\EEn$ with high probability, during $[0, \hat{t}_n]$, (i) all rewired edges are to individuals not in
$\mathcal{I}_n$, (ii) no individual in $\mathcal{I}_n$ acquires a rewired infectious edge at infection and (iii) no individual in $\mathcal{I}_n$ drops an infectious edge at infection.  Recalling that $\Zn_i$ and $Z^E_i$ are both defined on the same Poisson process $\eta_i$, we have that $|\Zn_i-Z^E_i| \sim {\rm Po}(|\mu_n(1-\frac{i}{n})-\mu|$ $(i=1,2,\dots,n)$.  Thus, noting that $\mu_n>\mu$,
\begin{align*}
\E\left[\sum_{i=1}^{\left \lfloor{\log n}\right \rfloor} |\Zn_i-Z^E_i|\right] 
&\le \log n \left(\mu_n-\mu+\mu_n\frac{\log n}{n}\right)\\
&\to 0 \quad \mbox{as } n \to \infty,
\end{align*}
since $\mu_n-\mu=O(1/n)$ as $n \to \infty$.  Thus, by Markov's inequality, $\Zn_i=Z^E_i$ $(i=1,2,\dots, \left \lfloor{\log n}\right \rfloor)$ with high probability.
Given (i)-(iii) above, it is then straightforward to define realisations of $\EEn$ $(n=1,2,\dots)$ from a realisation of $\B$ so that with high probability the process of infectives in $\EEn$ coincides with the realisation of $\B$ over $[0, t_n]$ and Theorem~\ref{thm:BPapprox} follows.

Turning to the proof of Lemma~\ref{lem:majbound}, note by Theorem~\ref{thm:BPapprox} and Remark~\ref{rmk:dropping} it is sufficient to consider the case $\alpha=0$, i.e.~the model with only dropping of edges.  Let $r$ denote the Malthusian parameter of the branching process $\B$. Then $r$ is given by the unique real solution of the equation
\[
\int_0^{\infty} \re^{-r t} \mu \lambda \re^{-(\lambda+\omega+\gamma)t} \,{\rm d}t =1,
\]
so $r=\lambda(\mu-1)-\gamma-\omega$.  For $t \ge 0$, let $I^E_1(t)$ be the number of individuals alive in $\B$ at time $t$ that had precisely $1$ infectious edge at birth and still have that edge at time $t$.  Recall that $T(t)$ is the total progeny of $\B$ at time $t$.  Let $A_{\rm ext}=\{\omega' \in \Omega: T(\infty, \omega')<\infty\}$ denote the set on which $\B$ goes extinct.  Then it follows using~\cite{Nerman81}, Theorem 5.4 and Corollary 3.2, that there exists a random variable $W_{\infty}\ge 0$, where $W_{\infty}(\omega')=0$ if and only if $\omega' \in A_{\rm ext}$, such that, as $t \to \infty$,
\begin{equation}
\label{equ:nerman1}
\re^{-r t} T(t) \convas \frac{r+\gamma}{r} W_{\infty} \qquad \mbox{and} \qquad \re^{-r t} I^E_1(t) \convas\frac{\re^{-\mu}(r+\gamma)}{\lambda} W_{\infty}.
\end{equation}

We use the construction in the above proof of Theorem~\ref{thm:BPapprox} to define realisations of $\EEn$ $(n=1,2,\dots)$ and $\B$ on $(\Omega, \mathcal{F}, \P)$.  Recall that $t_n=\inf\{t \ge 0: T(t)\ge \log n\}$.  It follows using Theorem~\ref{thm:BPapprox} and the first limit in~\eqref{equ:nerman1} that, with high probability, $\Tn>\log n$ if and only if $\omega' \in A_{\rm ext}^C$.  Using the second limit in~\eqref{equ:nerman1}, there exists $c_1>0$ such that $I^E_1(t_n)>c_1 \log n$ with high probability.

Recall that $R_0>1$.  There exists $\epsilon>0$ such that $R_0^{\epsilon}=(1-\epsilon)R_0>1$.  Note that in the above construction of $\EEn$, $\Zn_i \sge {\rm Po}(\mu(1-\epsilon))$ for $i \le n \epsilon$, where $\sge$ denotes the usual stochastic ordering.  Thus, with high probability, if $\Tn \ge \log n$, then until time $t$ such that $\Tn(t)\ge n \epsilon$, the process of infectives in the epidemic $\EEn$ is bounded below by  the branching process $\B^{\epsilon}$, which follows
$\B$ until time $t_n$, after which each birth is aborted independently with probability $\epsilon$.  Let $T^{\epsilon}(\infty)$ be the total progeny of $\B^{\epsilon}$.  Let $\pi_{\epsilon}$ be the extinction probability
for the progeny of an individual with one infectious edge at time $t_n$ in $\B^{\epsilon}$, so $\pi^{\epsilon}<1$ as $R_0^{\epsilon}>1$.  It follows
that
\begin{align*}
\liminf_{n \to \infty} \P(\Tn \ge \epsilon n | \Tn \ge \log n) &\ge \liminf_{n \to \infty} \P(T^{\epsilon}(\infty)\ge \epsilon n|t_n < \infty)\\
&\ge \liminf_{n \to \infty} \left(1-\pi_{\epsilon}^{c_1 \log n}\right)\\
&=1,
\end{align*}
and Lemma~\ref{lem:majbound} follows on setting $\tau'=\epsilon$.

\begin{rmk}
\label{rmk:susonly}
It is clear from the above proofs that Theorem~\ref{thm:BPapprox} and Lemma~\ref{lem:majbound} hold also for the model of Section~\ref{sec:SIRsusonly}, in which rewiring is necessarily to susceptibles.
\end{rmk}

\subsection{Proof of Theorem~\ref{thm:temporalWLLN}}
\label{sec:proofthmtemporalWLLN}

The proof is an application of~\cite{DN:08}, Theorem 4.1, and much of our notation is chosen to aid that application.

For $t \ge 0 $, let $\Sn(t)$ be the number of susceptibles at time $t$, $\In_i(t)$ be the number of infectives with $i$ infective edges at time $t$ ($i=0,1,\dots$) and $\Wn(t)$ be the number of rewired
susceptible-susceptible edges at time $t$.  Let $\Xn=\{\Xn(t):t \ge 0\}$, where $\Xn(t)=(\Sn(t), \In_0(t),\In_1(t),\dots,\Wn(t))$.  Then $\Xn$ is a continuous-time Markov chain with state space
$\En$ having typical element $\xi=(n^S, n_0^I, n_1^I,\dots,n^W)$.  Thus
\[
\En=\{(n^S, n_0^I, n_1^I,\dots,n^W):n^S,  n_0^I, n_1^I,\dots,n^W \in \mathbb{Z}_+, n^S+n^I \le n\},
\]
where $n^I=\sum_{k=0}^{\infty} n_k^I$ is the number of infectives corresponding to state $\xi$. For future reference, $n^E = \sum_{k=0}^{\infty} k n_k^I$ is the number of infectious edges corresponding to state $\xi$.

For distinct states $\xi, \xi' \in \En$, let $\qn(\xi,\xi')=\lim_{t \downarrow 0}t^{-1}\P(\Xn(t)=\xi'|\Xn(0)=\xi)$ denote the jump rate of $\Xn$ from $\xi$ to $\xi'$.
We do not write down expressions for these jump rates, as there are many different types of jumps
and some jump rates are quite complex owing to binomial samplings in the construction.  Moreover, the individual jump rates are not required for our purpose; see~\eqref{equ:betan} and~\eqref{equ:alphan} below.

Let $\bx:\En \to \mathbb{R}^4$ be defined by
\begin{equation}
\label{equ:coordb}
\bx(\xi)=(x^1(\xi), x^2(\xi), x^3(\xi), x^4(\xi))=n^{-1}\left(n^S, n^I, n^E, n^W\right)
\end{equation}
and $\bXnbar=\{\bXnbar(t):t \ge 0\}$, where $\bXnbar(t)=\bx(\Xn(t))$ $(t \ge 0)$.  Thus $\bXnbar$ corresponds to the process with the same notation in Theorem~\ref{thm:temporalWLLN}.  (To ease notation we have not indexed $\bx$ by $n$.  To aid connection with~\cite{DN:08}, the 4-dimensional vector $\bXnbar(t)$ is rendered in bold and the infinite-dimensional vector $\Xn(t)$ is not.)

For each $\xi \in \En$, let
\begin{equation}
\label{equ:betan}
\betan(\xi)=\sum_{\xi' \ne \xi} (\bx(\xi')-\bx(\xi))\qn(\xi,\xi')
\end{equation}
be the drift function for $\bXnbar$ and
\begin{equation}
\label{equ:alphan}
\alphan(\xi)=\sum_{\xi' \ne \xi} |\bx(\xi')-\bx(\xi)|^2\qn(\xi,\xi').
\end{equation}

To define the limiting deterministic process, let $U=[0,1]^2\times [0, \infty)^2$ and write
$x^1=s, x^2=i, x^3=i_E$ and $x^4=w$.  Let $b:U \to \mathbb{R}^4$ be the vector field defined by
\begin{equation}
\label{equ:bxpart}
b(\bx)=b_R(\bx)+b_D(\bx)+b_{RWR}(\bx)+b_{RWS}(\bx)+b_I(\bx),
\end{equation}
where
\begin{align*}
b_R(\bx)&=\gamma(0,-i,-i_E,0),\\
b_D(\bx)&=\omega(1-\alpha) i_E(0,0,-1,0),\\
b_{RWR}(\bx)&=\omega \alpha i_E(0,0,-(1-i-s),0),\\
b_{RWS}(\bx)&=\omega \alpha i_E(0,0,-s,s),\\
b_I(\bx)&=\lambda i_E (-1,1,-1+\mu s-\frac{y_E}{s}+2\frac{w}{s},-2\frac{w}{s}).
\end{align*}
(The suffices $R, D, RWR, RWS$ and $I$ correspond to  recovery, dropping of an infectious edge, rewiring of an infectious edge to a recovered individual, rewiring of an infectious edge to a susceptible individual and infection, respectively.  Rewiring of an infectious edge to an infective has no effect on
$\bx$.)
Let $x(t)=(s(t), i(t), i_E(t), w(t))$ be the solution of the differential equation $\dfrac{dx}{dt}=b(x)$
with initial state $x(0)$.  Observe that $b(\bx)$ is such that $x(t)$ coincides with the solution of~\eqref{equ:dsdt}-\eqref{equ:dwdt} in Section~\ref{sec:mainresults}.

To apply~\cite{DN:08}, Theorem 4.1, we need bounds for $|\betan(\xi)-b(\bx(\xi))|$ and $\alphan(\xi)$, which we now derive.  We partition $\betan(\xi)$ and $\alphan(\xi)$ into sums analogous to that for $b(\bx)$ at~\eqref{equ:bxpart}. We derive first expressions for the components of $\betan(\xi)$. Recall that $\xi=(n^S, n_0^I, n_1^I,\dots,n^W)$, $n^I=\sum_{k=0}^{\infty} n_k^I$ and $n^E = \sum_{k=0}^{\infty} k n_k^I$. Then
\begin{align}
\betan_R(\xi)&=n^{-1}\gamma (0, -n^I, -n^E,0) =b_R(\bx(\xi)),\nonumber\\
\betan_D(\xi)&=\omega (1-\alpha) n^E n^{-1}(0, 0, -1, 0)=b_D(\bx(\xi)),\nonumber\\
\betan_{RWR}(\xi)&=\omega\alpha n^E \left(\frac{n-n^S-n^I-2}{n-2}\right) n^{-1}(0,0,-1,0)\nonumber \\
&=\omega\alpha x^3(\xi)\left(\frac{n(1-x^1(\xi)-x^2(\xi))-2}{n-2}\right)(0,0,-1,0),\nonumber\\
\betan_{RWS}(\xi)&=n^E \omega\alpha \frac{n^S-1}{n-2} n^{-1}(0,0,-1,0)=\omega\alpha x^3(\xi)\frac{nx^1(\xi))-1}{n-2} (0,0,-1,0),\nonumber\\
\betan_{I}(\xi)&=\lambda n^E  n^{-1}
\left(-1,1,\frac{\mu_n(n^S-1)}{n}-1+2\frac{n^W}{n^S}-\frac{n^E-1}{n^S}, -2\frac{n^W}{n^S}\right)\label{equ:betanI}\\
&=\lambda x^3(\xi)\left(-1,1,\mu_n(x^1(\xi)-n^{-1})-1+2\frac{x^4(\xi)}{x^1(\xi)}
-\frac{x^3(\xi)-n^{-1}}{x^1(\xi)}, -2\frac{x^4(\xi)}{x^1(\xi)}\right).\nonumber
\end{align}

The expression for $\betan_{I}(\xi)$ is derived as follows.  (The derivations of the expressions for the other components of $\betan(\xi)$ are much simpler and hence omitted.) When $\Xn$ is in state $\xi$, infections occur at overall rate $\lambda n^E$.  A single infection clearly decreases $n^S$ by one and increases
$n^I$ by one.  The newly infected individual, $i_0$ say, acquires $X \sim {\rm Po}(\mu_n n^{-1}(n^S-1))$ un-rewired infectious edges and $Y \sim {\rm Bin}(n^W,2/n^S)$ rewired infectious edges.  The infectious edge that infected $i_0$ is dropped, as are $Z \sim {\rm Bin}(n^E-1,1/n^S)$ other infectious edges. (The ones that would attempt to infect $i_0$ if the edge was formed.)  Thus, in view of~\eqref{equ:coordb}, the associated change in $\bx(\xi)$ is
\begin{equation}
\label{equ:betanI1}
n^{-1}(-1, 1, X+Y-Z-1,Z).
\end{equation}
Summing $(\bx(\xi')-\bx(\xi))\qn(\xi,\xi')$ over $\xi' \ne \xi$ is equivalent to multiplying~\eqref{equ:betanI1} by $\lambda n^E$ and taking expectations with respect to $(X, Y, Z)$, which yields~\eqref{equ:betanI}.  For future reference, note that $X, Y$ and $Z$ are independent.

Simple calculation then yields
\begin{align}
|\betan_R(\xi)-b_R(\bx(\xi))|&=0,\label{equ:betmb1}\\
|\betan_D(\xi)-b_D(\bx(\xi))|&=0,\label{equ:betmb2}\\
|\betan_{RWR}(\xi)-b_{RWR}(\bx(\xi))|&=2(n-2)^{-1}\omega\alpha x^{3}(\xi)(x^{1}(\xi)+x^{2}(\xi)),\label{equ:betmb3}\\
|\betan_{RWS}(\xi)-b_{RWS}(\bx(\xi))|&=\sqrt{2}(n-2)^{-1}\omega \alpha x^3(\xi)|1-2x^1(\xi)|,\label{equ:betmb4}\\
|\betan_I(\xi)-b_I(\bx(\xi))|&=\lambda x^3(\xi) |(\mu_n-\mu)x^1(\xi)-n^{-1}\mu_n|\label{equ:betmb5}
\end{align}
and, letting $n^E_2 = \sum_{k=0}^{\infty} k^2 n_k^I$,
\begin{align}
\alphan_R(\xi)&=n^{-2}\gamma(n^I+n^E_2),\label{equ:alphan1}\\
\alphan_D(\xi)&=n^{-2}\omega(1-\alpha)n^E=n^{-1}\omega(1-\alpha)x^3(\xi),\label{equ:alphan2}\\
\alphan_{RWR}(\xi)&=n^{-2}\omega\alpha n^E \left(\frac{n-n^S-n^I-2}{n-2}\right)\nonumber\\
&=(n-2)^{-1} \omega\alpha x^3(\xi)(1-x^1(\xi)-x^2(\xi)-2n^{-1}),\label{equ:alphan3}\\
\alphan_{RWS}(\xi)&=n^{-2}\omega\alpha n^E 2 (n^S-1)/(n-2),\nonumber\\
&=2(n-2)^{-1}\omega\alpha x^3(\xi)(x^1(\xi)-n^{-1}).\label{equ:alphan4}
\end{align}

Calculation of $\alphan_{I}(\xi)$ is more involved; an upper bound is sufficient, which we now derive.
Using~\eqref{equ:betanI1},
\begin{equation}
\label{equ:alphanI}
\alphan_{I}(\xi)=\lambda n^E n^{-2}\E[2+(X+Y-Z-1)^2+Z^2].
\end{equation}
Noting that, $\E[X^2]=\E[X]+\E[X]^2$ (as $X$ has a Poisson distribution), $\E[Y^2] \le \E[Y](1+\E[Y])$ and $\E[Z^2] \le \E[Z](1+\E[Z])$ (as $Y$ and $Z$ have binomial distributions), a simple calculation yields
\[
\E[2+(X+Y-Z-1)^2+Z^2]<3+(\E[X]+\E[Y])^2+2\E[Z](2+\E[Z]),
\]
so
\begin{equation}
\label{equ:alphan5}
\alphan_{I}(\xi)<\lambda  n^{-1} x^3(\xi) \left[3+\left(\mu_n(x^1(\xi)-n^{-1})+\frac{2x^4(\xi)}{x^1(\xi)}\right)^2
+2\frac{x^3(\xi)}{x^1(\xi)}\left(2+\frac{x^3(\xi)}{x^1(\xi)}\right)
\right].
\end{equation}

We need bounds for $\Mn_E=\max_{0 \le t \le t_0} \Ine(t)$, $\Mn_{RW}=\max_{0 \le t \le t_0} \Wn(t)$
and $\Mn_D=\max_{0 \le t \le t_0} \max_{k \ge 0} (k:\In_k(t)>0)$ which hold with high probability.  Clearly both $\Mn_E$ and $\Mn_{RW}$
are bounded above by the total of the original degrees of all infectives during $[0,t_0]$, which in turn is bounded above by $\Sdn=\Dn_1+\Dn_2+\dots+\Dn_n$, where $\Dn_1, \Dn_2, \dots, \Dn_n$  are i.i.d.~${\rm Po}(\mu_n)$ random variables.  Since $\mu_n \to \mu$ as $n \to \infty$, it follows immediately from the law of large numbers that $\lim_{n \to \infty} \P(\Sdn > 2\mu n)=0$, whence
$\lim_{n \to \infty} \P(\Mn_E > 2\mu n)=\lim_{n \to \infty} \P(\Mn_{RW} > 2\mu n)=0$.  Let $\Mn=\max(\Dn_1,\Dn_2,\dots,\Dn_n)$.  Then it follows from~\cite{Kimber83} that
\begin{equation}
\label{equ:MNWHP}
\lim_{n \to \infty}\P(\Mn > \log n)=0.
\end{equation}

Turning to rewirings, let $\Rn_1, \Rn_2, \dots, \Rn_n$, where $\Rn_i$ is the total number of rewirings to the $i$th infective in $\EEn$, then $\Rn_i \sle \Rnhat_i$
$(i=1,2,\dots,n)$, where $\Rnhat_1, \Rnhat_2, \dots, \Rnhat_n$ are obtained by assuming that all $n$ individuals are infectious throughout $[0,t_0]$, where $t_0$ is as in the statement of Theorem~\ref{thm:temporalWLLN}.  Then, conditional upon $\Sdn$, the quantities
$\Rnhat_1, \Rnhat_2, \dots, \Rnhat_n$ are independent Poisson random variables, each having mean $\omega\alpha
t_0 \Sdn/(n-1)$.  (Each warning independently has probability $\frac{n-2}{n-1}$ of not being to $i$ and the probability that the person who received the warning rewires to $i$ is $\frac{1}{n-2}$.)  Given that
$\lim_{n \to \infty} \P(\Sdn > 2\mu n)=0$ it follows that with high probability $\Rnhat_1, \Rnhat_2, \dots, \Rnhat_n$ are bounded above by independent ${\rm Po}(3 \omega\alpha t_0 \mu)$ random variables, whence
\begin{equation}
\label{equ:maxRnWHP}
\lim_{n \to \infty}\P(\max_{i=1,2,\dots,n}\Rn_i> \log n)=0.
\end{equation}

Finally, we consider rewired susceptible-susceptible edges acquired by infectives at infection.  As seen later in the proof, we need consider only such acquired edges while the number of susceptibles $\Sn(t) \ge s_0 n$, for fixed $s_0 \in (0,1)$.   As noted above, $\Wn(t) \le \Sdn(t)$ for all $t \ge 0$.  For $i=2,3,\dots, n$, let $\Yn_{SS,i}$ be the number of rewired susceptible-susceptible edges acquired
by the $i$th infective at infection.  Recall that if the $i$th infection occurs at time $t$ then, conditionally, $\Yn_{SS,i} \sim {\rm Bin}\left(\Wn(t-), \frac{2}{\Sn(t-)}\right)$.  It follows that with high probability $\Yn_{SS,i} \sle \Ynhat_{SS,i}$ $(i=2,3,\dots,\left \lfloor{n(1-s_0)}\right \rfloor)$, where the $\Ynhat_{SS,i}$ are
i.i.d. ${\rm Bin}(\left \lceil{2\mu n}\right \rceil, \frac{2}{ns_0})$  random variables. Now $X \sle Y$ if $X \sim {\rm Bin}(1, p)$ and $Y \sim {\rm Po}(-\log(1-p))$, so with high probability $\Yn_{SS,i} \sle \Yncheck_{SS,i}$ $(i=2,3,\dots,\left \lfloor{n(1-s_0)}\right \rfloor)$, where the $\Yncheck_{SS,i}$ are
i.i.d. ${\rm Po}(f_n)$ random variables and $f_n=\mu n \log(1-\frac{2}{ns_0})$.  Now $f_n \to \frac{2\mu}{s_0}$ as $n \to \infty$ so, for all $n$ sufficiently large, $\Yncheck_{SS,i}$
$(i=2,3,\dots, \left \lfloor{n(1-s_0)}\right \rfloor)$ are stochastically smaller than i.i.d. ${\rm Po}(\frac{3\mu}{s_0})$ random variables. It then follows that
\begin{equation}
\label{equ:maxYSSnWHP}
\lim_{n \to \infty}\P(\max_{i=2,3,\dots, \left \lfloor{n(1-s_0)}\right \rfloor} \Yn_{SS,i}>\log n)=0.
\end{equation}

Equations~\eqref{equ:MNWHP}, \eqref{equ:maxRnWHP} and~\eqref{equ:maxYSSnWHP} imply that, provided the number of susceptibles remains above $ns_0$, then $\Mn_D$, the maximum number of infectious edges attached to
any individual throughout $[0,t_0]$, is bounded above by $3\log n$ with high probability.

Now $s(t_0)>0$.  Choose $s_0 \in (0,s(t_0))$.  Let $U_*=[s_0,1]\times[0,1] \times [0,2\mu]^2$.  The partial derivatives associated with the vector field $b$ are bounded on $U^*$, so $b$ is Lipschitz on
$U^*$, with Lipschitz constant $K$ say.  Choose $\epsilon \in (0, s(t_0)-s_0)$ such that the path
$x(t) (0 \le t \le t_0)$ lies at a distance greater than $\epsilon$ from the complement of $U_*$ and
let $\delta=\epsilon \re^{-K t_0}/3$.
Let
$\Tn=\inf\{t \ge 0: \bXnbar(t) \notin U_*\}$.  Then, with high probability, $x^1(\Xn(t)) \in [s_0, 1]$, $x^2(\Xn(t)) \in [0, 1]$, $x^3(\Xn(t)) \in [0, 2\mu]$ and $x^4(\Xn(t)) \in [0, 2\mu]$ for all
$t \in [0,\Tn]$.  It then follows using~\eqref{equ:betmb1}-\eqref{equ:betmb5} that with high probability
\[
|\betan(\Xn(t))-b((\bx(\Xn(t)))|\le (8+\sqrt{2})\mu \alpha \omega (n-2)^{-1}+2\lambda \mu (|\mu_n-\mu)+n^{-1}\mu_n),
\]
for all $t \in [0,\Tn]$, which implies that the event $\Omega_1$ occurs with high probability, where
\[
\Omega_1=\left\{\int_0^{\Tn \wedge t_0} |\betan(\Xn(t))-b(\bx(\Xn(t)))|\,{\rm d}t \le \delta\right\}.
\]

Further, \eqref{equ:alphan2}-\eqref{equ:alphan5} imply that with high probability
\begin{equation}
\label{equ:alphanb1}
\alphan_D(\Xn(t))+\alphan_{RWR}(\Xn(t))+\alphan_{RWS}(\Xn(t))\le 2 \mu \omega (n-2)^{-1}(1+2\alpha)
\end{equation}
and
\begin{equation}
\label{equ:alphanb2}
\alphan_I(\Xn(t))\le \lambda n^{-1}\left[3+(\mu_n+4\mu s_0^{-1})^2+8\mu s_0^{-1}(1+\mu s_0^{-1})\right],
\end{equation}
for all $t \in [0,\Tn]$.  Also, since $\Mn_D \le 3 \log n$ implies that $\sum_{k=1}^{\infty} k^2 \In_k(t) \le 3\log n \Ine(t)$, \eqref{equ:alphan1} implies that with high probability
\begin{equation}
\label{equ:alphanb3}
\alphan_R(\Xn(t))\le\gamma n^{-1}(1+3\mu \log n),
\end{equation}
for all $t \in [0,\Tn]$.  The inequalities~\eqref{equ:alphanb1}-\eqref{equ:alphanb3} imply that there exists $C \in (0,\infty)$ such that with high probability $\alphan(\Xn(t)) \le C n^{-1} \log n$ for all
$t \in [0,\Tn]$, whence the event $\Omega_2$ occurs with high probability, where
\[
\Omega_2=\left\{\int_0^{\Tn \wedge t_0} \alphan(\Xn(t))\,{\rm d}t \le A_n t_0\right\}
\]
and $A_n=C n^{-1} \log n$.  The event $\Omega_0=\left\{\left|\bXnbar(0)-x(0)\right| \le \delta\right\}$ occurs with
high probability as $\bXnbar(0) \convp x(0)$ as $n \to \infty$.  By~\cite{DN:08}, Theorem 4.1,
\begin{equation}
\label{equ:psup}
\P\left(\sup_{0 \le t \le t_0} \left|\bXnbar(t)-x(t)\right|>\epsilon\right) \le \frac{4 A_n}{\delta^2}
+\P(\Omega_0^c\cup \Omega_1^c \cup \Omega_2^c)
\end{equation}
and Theorem~\ref{thm:temporalWLLN} follows since the right-hand side of~\eqref{equ:psup} converges to $0$ as $n \to \infty$ and $\epsilon>0$ can be chosen arbitrarily small.

\subsection{Proof of Theorem~\ref{thm:SIfinalsize}}
\label{sec:proofthmSIfinal}

Let $\zetan=\inf\{t \ge 0: \Ine(t)=0\}$, then the final size $\Tn$ of the epidemic is given by
$\Tn=n-\Sn(\zetan)$.  Let $\Xnt=\{\Xnt(t):t \ge 0\}$, where $\Xnt(t)=(\Snt(t), \Int_0(t), \Int_1(t), \dots,
\Wnt(t))$, be the process having the same state space as $\Xn$, viz.~$\En$, and jump rates $\qtn(\xi,\xi')=\qn(\xi,\xi')/(\lambda n^E)$ $(\xi,\xi' \in \En, \xi \ne \xi')$.  Recall that $n^E$ is the total number of infective edges when corresponding to state $\xi$.  The process $\Xnt$ stops when it reaches a state $\xi$ with $n^E=0$.  For $t \ge 0$, let $\Inet(t)=\sum_{k=1}^{\infty} k \Int_k(t)$. Let $\zetatn=\inf\{t \ge 0: \Inet(t)=0\}$ and $\Tnt=n-\Snt(\zetatn)$.  Note that $\Xnt$ is a random time-scale transformation of $\Xn$, in which when $\Xn$ is in state $\xi$ the clock for $\Xnt$ is slowed down by the factor $\lambda n_E(\xi)$ (i.e.~by the total force of infection), so $\Tn$ and $\Tnt$ have the same distribution.  Consequently we can use the time-transformed process $\Xnt$ and $\Tnt$ to study the asymptotic behaviour of $\Tn$ as $n \to \infty$.  The advantage of using the time-transformed process is that, as $n \to \infty$, $\zetatn \convp \zetat<\infty$ (where $\zetat$ is non-random and defined below), whereas $\zetan \convp \infty$.
Such random time-scale transformations have a long history of application to the final outcome of epidemics; see, for example, \cite{Watson80}, \cite{EK86}, page 467, \cite{DN:08}, Section 5.2, and~\cite{JLW:2014}.

Let $\bXnt=\{\bXnt(t):t \ge 0\}$, where $\bXnt(t)=\bxt(\Xnt(t))$ $(t \ge 0)$ and $\bxt: \En \to \mathbb{R}^3$ is defined by
\begin{equation}
\label{equ:cordbt}
\bxt(\xi)=(x^1(\xi), x^3(\xi), x^4(\xi))=n^{-1}\left(n^S, n^E, n^W\right);
\end{equation}
cf.~\eqref{equ:coordb}, there is no need to keep track of the number of infectives as the model is SI so
$n^I=n-n^S$.  (The labelling of the three coordinates of $\bxt(\xi)$ is chosen to aid connection with the coordinates of $\bx(\xi)$ at~\eqref{equ:coordb}.)  The deterministic approximation to $\bXnt$, obtained by dividing the right-hand sides of~\eqref{equ:dsdt}, \eqref{equ:diEdt} and~\eqref{equ:dwdt} by $\lambda i_E$, setting $\gamma=0$ and
using $i+s=1$,
is $\xt(t)=(\st(t), \iet(t), \wt(t))$ $(t \ge 0)$, where
\begin{align}
\dfrac{d\st}{dt}&=-1, \label{equ:dsdtA}\\
\dfrac{d\iet}{dt}&= -1 +\mu \st - \frac{\iet}{\st} +2 \frac{\wt}{\st}-\frac{\omega}{\lambda} (1-\alpha+\alpha \st),\label{equ:diEdtA}\\
\dfrac{d\wt}{dt}&= \frac{\omega \alpha}{\lambda} \st-2 \frac{\wt}{\st}.\label{equ:dwdtA}
\end{align}

Suppose that $\bXnbar(0) \convp (1-\epsilon, \epsilon, \mu \epsilon(1-\epsilon),0)$ as $n \to \infty$, where $\epsilon \in (0, 1)$, as in part (a) of Theorem~\ref{thm:SIfinalsize}. Then $\bXnt(0) \convp (1-\epsilon, \epsilon, \mu \epsilon(1-\epsilon),0)$ also.  It is shown in Appendix~\ref{app:ttdetSIsoln} that the solution of~\eqref{equ:dsdtA} to~\eqref{equ:dwdtA}, with initial condition $(\st(0), \iet(0), \wt(0))=(1-\epsilon, \mu \epsilon(1-\epsilon),0)$, is
\begin{align}
\st(t)&=1-\epsilon-t,\label{equ:st}\\
\iet(t)&=\st(t)\left\{\left[1+\frac{\omega}{\lambda}(1-\alpha)\right]\log\left(\frac{\st(t)}{1-\epsilon}\right)
+\left(\mu+\frac{\omega\alpha}{\lambda}\right)(1-\st(t))\right.\nonumber\\
&\qquad\qquad\qquad\left.-\frac{\omega\alpha\epsilon}{\lambda}
+\frac{2\omega\alpha}{\lambda}\st(t)\log\left(\frac{\st(t)}{1-\epsilon}\right)\right\},\label{equ:iet}\\
\wt(t)&=-\frac{\omega\alpha}{\lambda}(\st(t))^2 \log\left(\frac{\st(t)}{1-\epsilon}\right).\label{equ:wt}
\end{align}

Let $\Ut=[0,1]\times (0,\infty) \times [0, \infty)$. Let $\bt: \Ut \to \mathbb{R}^3$ be the vector field
given by the right-hand sides of~\eqref{equ:dsdtA} to~\eqref{equ:dwdtA}.
Let $\zetat=\inf\{t>0:\iet(t)=0\}$ and $\tau=1-\st(\zetat)$.  Then it follows from~\eqref{equ:iet} that
$\tau$ is given by the smallest solution in $(\epsilon,1)$ of $F_{\epsilon}(x)=0$, where $F_{\epsilon}$ is defined at~\eqref{equ:Feps}.
(We show in Appendix~\ref{app:ttdeteps>0} that such a solution exists.)
We follow Section 4.2 of~\cite{DN:08}, adapted to $L^2$ estimates rather than exponential estimates as there presented.  For $\epsilon'>0$, let
\begin{align*}
\zetat_{\epsilon'}^-&=\inf\{t \ge 0: B(\xt(t),\epsilon') \cap \Ut^c \neq \emptyset\},\\
\zetat_{\epsilon'}^+&=\inf\{t \ge 0: B(\xt(t),\epsilon') \subseteq \Ut^c\},
\end{align*}
where $B(\xt(t),\epsilon')$ is the closed ball $\{\xt \in \mathbb{R}^3: |\xt-\xt(t)| \le \epsilon'\}$, and let
\[
\rhot(\epsilon')=\sup_{\zetat_{\epsilon'}^- \le t \le \zetat_{\epsilon'}^+}|\xt(t)-\xt(\zetat)|.
\]
We show in Appendix~\ref{app:ttdeteps>0} that, under the conditions of Theorem~\ref{thm:SIfinalsize}(a), $\iet'(\zetat)<0$.  Thus
$\bt(\xt(\zetat))$ is not tangent to the boundary of $\Ut$, so there exists $\Ct < \infty$ such that
$\rhot(\epsilon')< \Ct \epsilon'$ for all sufficiently small $\epsilon'>0$.  Choose such an $\epsilon'$
and then $\ttild_0$ so that $\ttild_0>\zetat_{\epsilon'}^+$.  Let $\st_0=\st(\ttild_0)$ and $\itt^E_0=\min\{\iet(t):0 \le t \le \ttild_0\}$.  Let $\Ut^*=[s_0,1]\times (0,2\mu] \times [0, 2\mu]$.  The vector field $\bt$, now defined on
$[s_0,1]\times (2\itt^E_0,2\mu] \times [0, 2\mu]$, is Lipschitz, with Lipschitz constant $\Kt$ say.

Let $\betatn$ and $\alphatn$ be defined analogously to $\betan$ and $\alphan$ at~\eqref{equ:betan} and~\eqref{equ:alphan}.
Let $\delta'=\epsilon'\re^{-\Kt \ttild_0}/3$ and $\tilde{\Omega}_i$ $(i=0,1,2)$ be corresponding events to
$\Omega_i$ $(i=0,1,2)$ (see the proof of Theorem~\ref{thm:temporalWLLN}) for $\Xnt$.
The components of $\betatn(\xi)$ are given by the corresponding components of $\betan(\xi)$ divided by
$\lambda x^{(3)}(\xi)$.  (Note that $\betan_R$ and $\betan_{RWR}$ are no longer relevant.)  A similar comment applies to the components of $\bt(\bxt(\xi))$.  In an obvious notation, $\alphatn_D(\xi)=
\alphan_D(\xi)/(\lambda x^3(\xi))$, $\alphatn_{RWS}(\xi)=
\alphan_{RWS}(\xi)/(\lambda x^3(\xi))$ and
\begin{equation}
\label{equ:alphantI}
\alphatn_{I}(\xi)=n^{-2}\E[1+(X+Y-Z-1)^2+Z^2];
\end{equation}
cf.~\eqref{equ:alphanI}, the constant $2$ in the expectation in~\eqref{equ:alphanI} is replaced by $1$ in~\eqref{equ:alphantI}, since $\bxt(\xi)$ does not have the component $x^2(\xi)=n^{-1}n^I$ corresponding to the number of infectives.  (Again, $\alphan_R$ and $\alphan_{RWR}$ are no longer relevant.)
The proof of Theorem~\ref{thm:temporalWLLN} is easily modified to show that the events $\tilde{\Omega}_1$ and
$\tilde{\Omega}_2$, with $A_n$ replaced by $\At_n=C_1n^{-1}\log n$ for suitable $C_1 >0$, each occur with high probability.  It then follows using~\cite{DN:08}, Theorem 4.3 (adapted to $L^2$ estimates) that
\begin{equation}
\label{equ:psupA}
\P\left(\left|\bXnt(\zetatn)-\xt(\zetat)\right|>(1+\Ct)\epsilon'\right) \le \frac{4 \At_n}{(\delta')^2}
+\P(\tilde{\Omega}_0^c\cup \tilde{\Omega}_1^c \cup \tilde{\Omega}_2^c).
\end{equation}
Theorem~\ref{thm:SIfinalsize}(a) follows since $\tilde{\Omega}_0$ occurs with high probability, by assumption, and $\epsilon'>0$ can be chosen arbitrarily small.

Turn now to Theorem~\ref{thm:SIfinalsize}(b).  Consider first the non-time-transformed process, let $t_n=\inf\{t>0:\Sn(t) \le n-\log n\}$ and note that $t_n< \infty$ since a major outbreak occurs.  Let $\Xnprime=\{\Xnprime(t):t \ge 0\}$ be defined by $\Xnprime(t)=\Xn(t+t_n)$ $(t \ge 0)$, so $\Xnprime(0)=\Xn(t_n)$ is random.  We apply the same random time-scale transformation as above to the process $\Xnprime$, yielding a process, $\Xntprime$ say, with $\Xntprime(0)=\Xn(t_n)$. We now drop the superfix $'$, from the process $\Xntprime$ to yield a time-transformed process $\Xnt$ with $\Xnt(0)=\Xn(t_n)$.

It follows from the proofs of Theorem~\ref{thm:BPapprox} and Lemma~\ref{lem:majbound} that $\Ine(t_n)$ and $\Wn(t_n)$ are each $O_p(\log n)$, so $\bXnt(0) \convp (1, 0, 0)$ as $n \to \infty$.
The deterministic approximation $\xt(t)=(\st(t), \iet(t), \wt(t))$ $(t \ge 0)$ is now the solution of~\eqref{equ:dsdtA}-\eqref{equ:dwdtA} with initial condition $(\st(0), \iet(0), \wt(0))=(1,0,0)$, which is given by setting $\epsilon=0$ in~\eqref{equ:st}-\eqref{equ:wt}.  (Note that, unlike the corresponding un-time-transformed deterministic approximation, this process does not get stuck at
$(1,0,0)$.)  As before, $\zetat=\inf\{t>0:\iet(t)=0\}$ and $\tau=1-\st(\zetat)$.  Thus $\tau$ satisfies $F_0(\tau)=0$.    We show in Appendix~\ref{app:ttdeteps=0} that, when $R_0>1$, $F_0(\tau)=0$ has a unique solution in $(0,1)$
and $\iet'(\zetat)<0$.

Let $\tau'$ be as in Lemma~\ref{lem:majbound}, $t_1=\tau'/2$ and $\ct=\st(t_1)=1-t_1$.  We now take
\[
\Ut=\left\{\left([1-\ct,1] \times [-1,2\mu)\right)\cup\left([0,1-\ct)\times (0, 2\mu)\right)\right\} \times [0,2\mu).
\]
Note that with this choice of $\Ut$, $\zetat_{\epsilon'}^-$ is close to $\zetat$ for all sufficiently small $\epsilon'>0$.  Note also that $\Inet(t)>0$ for all $t \in [0,t_1]$ with high probability, since if $\Inet(t)=0$ for some $t \in [0,t_1]$ then $\Tnbar \le \tau'/2$ contradicting Lemma~\ref{lem:majbound}. The remainder of the proof now follows that of (a).

\subsection{Proof of Corollary~\ref{cor:SIfinal}}
\label{sec:proofcorSIfinal}
Recall that $\gamma=0$ and we consider the SI model.
By Theorem~\ref{thm:SIfinalsize}(b), $\tau_{SI}(\mu, \lambda, \omega, \alpha)$ is given by the unique solution in $(0,1)$ of $f(x)=0$, where
\begin{equation}
\label{equ:f}
f(x)=\log(1-x)+\frac{(\lambda\mu+\omega\alpha)x}{\lambda+\omega(1-\alpha)+2\omega\alpha(1-x)}.
\end{equation}
Recall from~\eqref{equ:lambdac} that $\lambda_c=\frac{\omega}{\mu-1}$.  The function $f_0$, defined at~\eqref{equ:f0def} is obtained by substituting $\lambda=\lambda_c$ in~\eqref{equ:f}. 

Writing $\theta$ for $\theta(\mu, \alpha)=\frac{2\alpha(\mu-1)}{\mu+\alpha(\mu-1)}$, elementary calculus yields
\[
f_0'(x)=-\frac{1}{1-x}+\frac{1}{(1- \theta x)^2} \qquad\mbox{and}\qquad f_0''(x)=-\frac{1}{(1-x)^2}+\frac{2\theta}{(1- \theta x)^3}.
\]
Thus $f_0'(0)=0$.  If $\theta \in [0, \frac{1}{2}]$ then $f_0'(x)<0$ for all $x \in (0,1)$, whence $f_0$ is strictly decreasing on $[0,1)$ and $0$ is the only root of $f_0$ in $[0,1)$.  If $\theta>\frac{1}{2}$ then $f_0''(0)>0$ and $f_0$ has a unique root in $(0,1)$, since $f_0(1-)=-\infty$ and $f$ has at most two stationary points in $(0,1)$. Part (a) follows since, with $\mu, \omega$ and $\alpha$ held fixed, $\lim_{\lambda \downarrow \lambda_c}f(x)=f_0(x)$ for all $x \in [0,1)$. Part (b) follows on noting that $\theta>\frac{1}{2}$ if and only if $\alpha>\frac{1}{3}$ and $\mu>\frac{3\alpha}{3\alpha-1}$.

Turning to part (c), fix $\mu>1, \omega>0$ and $\alpha>0$, and let
\[
g(x,\lambda)=\log(1-x)+\frac{(\lambda\mu+\omega\alpha)x}{\lambda+\omega(1-\alpha)+2\omega\alpha(1-x)}
\qquad(0 \le x<1,\lambda>\lambda_c),
\]
so $\tau_{SI}(\mu, \lambda, \omega, \alpha)$ is given by the unique $x \in (0,1)$ such that $g(x,\lambda)=0$.  Let $g_{\lambda}$ denote the partial derivative of $g$ with respect to $\lambda$.  Then
\[
g_{\lambda}(x,\lambda)=\frac{\omega x[(1+\alpha-2\alpha x)\mu-\alpha]}
{[\lambda+\omega(1-\alpha)+2\omega\alpha(1-x)]^2},
\]
so
\begin{equation}
\label{equ:g_lambda}
g_{\lambda}(x,\lambda) \begin{cases}
              >0 & \text{ if } x<x_0(\mu,\alpha),\\
              =0 & \text{ if } x=x_0(\mu,\alpha),\\
              <0 & \text{ if } x>x_0(\mu,\alpha),
              \end{cases}
\end{equation}
where (recall~\eqref{equ:x0})
\[
x_0(\mu,\alpha)=\frac{1+\alpha}{2\alpha}-\frac{1}{2 \mu}.
\]
Note that $x_0(\mu,\alpha)$ is independent of both $\lambda$ and $\omega$.

Now $\mu>1$, so $x_0(\mu,\alpha)>\frac{1}{2\alpha}$.  Thus if $\alpha \le \frac{1}{2}$, then
$g_{\lambda}(x,\lambda)>0$ for all $x \in (0,1)$ and it follows that $\tau_{SI}(\mu, \lambda, \omega, \alpha)$ is strictly increasing in $\lambda$ for $\lambda>\lambda_c$.

Now fix $\alpha\in (\frac{1}{2},1)$.  Noting that $\omega(1-\alpha)+2\omega\alpha[1-x_0(\mu,\alpha)]=\frac{\omega\alpha}{\mu}$,
\begin{equation}
\label{equ:htilde}
g(x_0(\mu,\alpha),\lambda)=\tilde{h}(\mu)\qquad(\lambda>\lambda_c),
\end{equation}
where $\tilde{h}(\mu)=h(\mu,\alpha)$ and $h$ is defined at~\eqref{equ:h}.  (Note that $g(x_0(\mu,\alpha),\lambda)$ is independent of $\lambda$.)  If $\mu \ge \frac{\alpha}{1-\alpha}$ then $x_0(\mu,\alpha) \ge 1$ whence, as above, $\tau_{SI}(\mu, \lambda, \omega, \alpha)$ is strictly increasing in $\lambda$.  Thus, consider $1<\mu<\frac{\alpha}{1-\alpha}$.

Now $\tilde{h}(1)=\log\left(1-\frac{1}{2\alpha}\right)+\frac{1}{2\alpha}<0$ and $\tilde{h}(\frac{\alpha}{1-\alpha})=-\infty$.  Further,
\begin{equation}
\label{equ:htildemu}
\tilde{h}'(\mu)=\frac{1+\alpha}{2\alpha}-\frac{\alpha}{\mu[\alpha-(1-\alpha)\mu]}.
\end{equation}
Hence, $\tilde{h}'(1)=\frac{\alpha-1}{2\alpha(2\alpha-1)}<0$, as $\alpha\in (\frac{1}{2},1)$, and
\begin{equation}
\label{equ:h'}
\tilde{h}'(\mu)=0 \mbox{ if and only if } (1-\alpha^2)\mu^2-\alpha(1+\alpha)\mu +2\alpha^2=0.
\end{equation}
The discriminant of the quadratic (in $\mu$) on the right-hand side of~\eqref{equ:h'} is $\Delta=\alpha^2(1+\alpha)(9\alpha-7)$.  Suppose $\alpha<\frac{7}{9}$.  Then $\Delta<0$, so $\tilde{h}$
has no stationary point in $(1, \frac{\alpha}{1-\alpha})$ and $\tilde{h}(\mu)<0$ for all $\mu \in (1,\frac{\alpha}{1-\alpha})$.  Thus, for such
$\mu$, \eqref{equ:htilde} implies that $g(x_0(\mu,\alpha),\lambda)<0$ for all $\lambda >\lambda_c$, whence
$\tau_{SI}(\mu, \lambda, \omega, \alpha)< x_0(\mu,\alpha)$ for all $\lambda >\lambda_c$.  It then follows, using~\eqref{equ:g_lambda} and the definition of $\tau_{SI}(\mu, \lambda, \omega, \alpha)$, that
$\tau_{SI}(\mu, \lambda, \omega, \alpha)$ is strictly increasing in $\lambda$ on $(\lambda_c,\infty)$.
The same conclusion holds if $\alpha=\frac{7}{9}$ since then $\tilde{h}$ has a unique stationary point in $(1, \frac{\alpha}{1-\alpha})$, which must be an inflection as $\tilde{h}'(1)<0$ and $\tilde{h}(\frac{\alpha}{1-\alpha})=-\infty$, so again $\tilde{h}(\mu)<0$ for all $\mu \in (1,\frac{\alpha}{1-\alpha})$.

Suppose now that $\frac{7}{9}<\alpha<1$.  Then $\tilde{h}$ has two stationary points in $(1,\frac{\alpha}{1-\alpha})$, a minimum at $\check{\mu}(\alpha)$ and a maximum at $\hat{\mu}(\alpha)$ say, where $\check{\mu}(\alpha)<\hat{\mu}(\alpha)$ are given by $\frac{\alpha}{2(1-\alpha)}[1 \pm \eta(\alpha)]$, with $\eta(\alpha)=\sqrt{\frac{9\alpha-7}{1+\alpha}}$.  Note that $\eta$ is strictly increasing on $[\frac{7}{9},1]$, $\eta(\frac{7}{9})=0$ and $\eta(1)=1$.  Further, $\eta^{-1}(\theta)=\frac{7+\theta^2}{9-\theta^2}$ $(\theta \in [0,1])$. For $\theta \in [0,1)$, let
$\hat{h}(\theta)=h(\hat{\mu}(\alpha),\alpha)$, where $\alpha=\eta^{-1}(\theta)$.  A little algebra yields
that
\[
\hat{h}(\theta)=2\log(1-\theta)-\log(7+\theta^2)+\frac{3+\theta}{2(1-\theta)}.
\]
Now $\hat{h}(0)=\frac{3}{2}-\log 7<0$, consistent with the above analysis of the case $\alpha=\frac{7}{9}$.  It is shown easily that $\hat{h}$ is strictly increasing and $\hat{h}(1-)=\infty$.
Thus $\hat{h}$ has a unique root $(0,1)$,  which we denote by $\theta^*$.  Numerical calculation yields
$\theta^* \approx 0.4614$, whence $\alpha^*=\eta^{-1}(\theta^*)\approx 0.8209$ and $\hat{\mu}(\alpha^*)\approx 3.3482$.

When $\frac{7}{9}<\alpha<\alpha^*$, then for all $\mu \in (1,\frac{\alpha}{1-\alpha})$ we have $\tilde{h}(\mu)<0$ and the previous analysis shows that $\tau_{SI}(\mu, \lambda, \omega, \alpha)$ is strictly increasing in $\lambda$ on $(\lambda_c,\infty)$. (Thus we have now proved part (c)(i).) The same is true when $\alpha=\alpha^*$
and $\mu \ne \hat{\mu}(\alpha^*)$.  If $\alpha=\alpha^*$ and $\mu =\hat{\mu}(\alpha^*)$ then
$g(x_0(\hat{\mu}(\alpha^*)),\lambda)=0$ for all $\lambda > \lambda_c$, whence
$\tau_{SI}(\hat{\mu}(\alpha^*), \lambda, \omega, \alpha^*)= x_0(\hat{\mu}(\alpha^*))=\tau^*$ for all
$\lambda > \lambda_c$, thus proving part (c)(ii).

Suppose now that $\alpha^*<\alpha<1$.  Then $\tilde{h}(\hat{\mu}(\alpha))>0$, so there exists
$\mu_L^*(\alpha)<\hat{\mu}(\alpha)<\mu_U^*(\alpha)$ such that (recall that $1<\mu<\frac{\alpha}{1-\alpha}$)
\[
\tilde{h}(\mu)
\begin{cases}
              <0 & \text{ if } 1<\mu<\mu_L^*(\alpha) \text{ or } \mu_U^*(\alpha)<\mu<\frac{\alpha}{1-\alpha},\\
              =0 & \text{ if } \mu=\mu_L^*(\alpha) \text{ or } \mu_U^*(\alpha) ,\\
              >0 & \text{ if } \mu_L^*(\alpha)<\mu<\mu_U^*(\alpha).
              \end{cases}
\]
The cases when $\tilde{h}(\mu)<0$ or $\tilde{h}(\mu)=0$ have been analysed previously and the same conclusions follow.  If $\tilde{h}(\mu)>0$ then $g(x_0(\mu,\alpha),\lambda)>0$ for all $\lambda >\lambda_c$.  Thus now
$\tau_{SI}(\mu, \lambda, \omega, \alpha)> x_0(\mu,\alpha)$ for all $\lambda >\lambda_c$, so~\eqref{equ:g_lambda} and the definition of $\tau_{SI}(\mu, \lambda, \omega, \alpha)$ implies that
$\tau_{SI}(\mu, \lambda, \omega, \alpha)$ is strictly decreasing in $\lambda$ on $(\lambda_c,\infty)$,
which completes the proof of part (c)(iii).

Suppose finally that $\alpha^*<\alpha=1$.  Then $\tilde{h}(\mu)=-\log(2\mu)+\mu-\frac{1}{2}$, so
$\tilde{h}'(\mu)=1-\mu^{-1}>0$ for all $\mu>1$.  Further, it is shown
easily that $\tilde{h}$ has a unique root, denoted by $\mu^*(1)$, in $(1,\infty)$. Thus $\tilde{h}(\mu)$
is $<0, =0$ or $>0$ according to $\mu<\mu*(1), =\mu^*(1)$
or $>\mu_*(1)$ and part (c)(iv) follows.


\subsection{Final outcome and discontinuity of the SIR model}
\label{sec:proofthmSIRfinal}
In this section we present arguments in support of conjectures concerning the final size of the SIR epidemic (Conjecture~\ref{conj:SIRfinal}) and its discontinuity at the critical value
$\lambda=\lambda_c$ (see Remark~\ref{rmk:disc}), together with proofs of Theorems~\ref{thm:detdisc} and~\ref{thm:SIRdiscon}.  In Section~\ref{sec:conjSIRsupp}, we first outline the difficulty in extending the proof in Section~\ref{sec:proofthmSIfinal} of Theorem~\ref{thm:SIfinalsize}, concerning the final outcome of the SI epidemic, to the SIR model.  We
then make some remarks on how Conjecture~\ref{conj:SIRfinal} might be proved.  In Section~\ref{sec:detdisc}, we prove Theorem~\ref{thm:detdisc} concerning discontinuity at $\lambda=\lambda_c$ of the final size
of the deterministic SIR model given by~\eqref{equ:dsdt}-\eqref{equ:dwdt} when its initial condition converges to having zero infected in an appropriate fashion.  Finally, in Section~\ref{sec:proofSIRdisc} we prove Theorem~\ref{thm:SIRdiscon}, which gives sufficient conditions and (implicitly) necessary
conditions for the final size of the stochastic SIR model to have a discontinuity at $\lambda=\lambda_c$, and discuss a possible approach to proving that the sufficient conditions are also necessary.

\subsubsection{Final size of SIR epidemic}
\label{sec:conjSIRsupp}
A similar random time-scale transformation to that used in Section~\ref{sec:proofthmSIfinal}, in which the clock is slowed down by a factor given by the total force of infection, leads to the deterministic model
\begin{align}
\dfrac{d\st}{dt}&=-1, \label{equ:dsdtB}\\
\dfrac{d\itt}{dt}&=-\frac{\gamma}{\lambda}\frac{\itt}{\iet}+1,\label{equ:didtB}\\
\dfrac{d\iet}{dt}&= -1-\frac{\gamma}{\lambda}+\mu \st - \frac{\iet}{\st} +2 \frac{\wt}{\st}-\frac{\omega}{\lambda} (1-\alpha+\alpha(1-\itt)),\label{equ:diEdtB}\\
\dfrac{d\wt}{dt}&= \frac{\omega \alpha}{\lambda} \st-2 \frac{\wt}{\st},\label{equ:dwdtB}
\end{align}
which is obtained from~\eqref{equ:dsdt}-\eqref{equ:dwdt} by dividing all right-hand sides by $\lambda i_E$.

Note that the corresponding vector field, $b$ say, is not Lipschitz on, for example,
$U_*=[s_0,1]\times[0,1] \times [0,2\mu]^2$, owing to the term $\frac{\itt}{\iet}$ in~\eqref{equ:didtB}.
Also, the quantity corresponding to $\alphan_R(\xi)$ in~\eqref{equ:alphan1} becomes
\[
\alphatn_R(\xi)=n^{-2}\left(\frac{\gamma}{\lambda}\right)\left(\frac{n^I+n^E_2}{n^E_2}\right),
\]
leading to the bound (cf.~\eqref{equ:alphanb3})
\begin{equation*}
\alphatn_R(\Xnt(t))\le n^{-1}\left(\frac{\gamma}{\lambda}\right)\left(\frac{x^1(\xi)}{x^3(\xi)}+3\log n\right),
\end{equation*}
which causes problems in the application of~\cite{DN:08}, Theorem 4.1, if $x^3(\xi)$ can be arbitrarily small.  These problems disappear if both $\iet$ and $x^3(\xi)$ are bounded away from zero but that excludes the end of an epidemic, and also the start if there are few initial infectives.

Consider Conjecture~\ref{conj:SIRfinal}(a).  For any $\epsilon \in (0, 1- s(\infty))$, it follows immediately from Theorem~\ref{thm:temporalWLLN}, by choosing $t_0$ so that $s(t_0)=s(\infty)+\frac{\epsilon}{2}$,
that
\[
\lim_{n \to \infty}\P(\Tnbar > 1-s(\infty)-\epsilon)=1.
\]
To prove Conjecture~\ref{conj:SIRfinal}(a) it is sufficient to show also that
\begin{equation}
\label{equ:TnbaruB}
\lim_{n \to \infty}\P(\Tnbar < 1-s(\infty)+\epsilon)=1.
\end{equation}
One approach to proving~\eqref{equ:TnbaruB} is as follows (cf.~the proof of~\cite{BB:2007}, Theorem 4.3).  For $\epsilon' \in (0, 1-s(\infty))$, let $t_0(\epsilon')$ be defined by $s(t_0(\epsilon'))=s(\infty)+\epsilon'$. Choose $\epsilon_0'$ sufficiently small so that the process
$\{\bXn(t): t\ge t_0(\epsilon_0')\}$ is subcritical; intuitively this is possible since otherwise the final fraction susceptible would be smaller than
$s(\infty)$.  For $\epsilon' \le \epsilon_0'$, the process of infectives in $\{\bXn(t): t\ge t_0(\epsilon')\}$ can be bounded above by a subcritical branching process, $\B(\epsilon')$ say,
having $2ni(t_0(\epsilon'))$ initial ancestors.  Let $\mu_B(\epsilon')$ be the mean total progeny of $\B(\epsilon')$ if it had a single ancestor.  Note that $\mu_B(\epsilon')$ is increasing in $\epsilon'$ and $\mu_B(\epsilon')<\infty$ for $\epsilon'<\epsilon_0'$.  A simple argument using Markov's inequality shows that, for $\epsilon>0$,
\begin{equation}
\label{equ:TnbaruB1}
\P\left(\Tnbar(\infty)-\Tnbar(t_0(\epsilon')) \ge \epsilon\right) \le \frac{2i(t_0(\epsilon'))\mu_B(\epsilon_0')}{\epsilon}.
\end{equation}
The right-hand side of~\eqref{equ:TnbaruB1} can be made arbitrarily small as $i(t_0(\epsilon'))\to 0$ as $\epsilon' \downarrow 0$, and~\eqref{equ:TnbaruB} follows.  The lack of a closed-form solution to the ordinary differential equations~\eqref{equ:dsdtB}-\eqref{equ:dwdtB} makes it difficult to make the above argument rigorous.

Turning to Conjecture~\ref{conj:SIRfinal}(b), consider the approximating branching process $\B$ defined in Section~\ref{sec:SIR}.  Application of~\cite{Nerman81}, Theorem 5.4 and Corollary 3.2, as at~\eqref{equ:nerman1}, yields
\begin{equation*}
\re^{-r t} I(t) \convas  W_{\infty} \quad \mbox{as } n \to \infty,
\end{equation*}
so, for $\P$-almost all $\omega \in A_{\rm ext}^C$,
\begin{align*}
\lim_{t \to \infty}\frac{I(t)}{I_E(t)}&=\frac{\lambda}{r+\gamma}\\
&=\frac{\lambda}{\lambda(\mu-1)-\omega},
\end{align*}
as $r=\lambda(\mu-1)-\gamma-\omega$.  It is then clear from the proof of Theorem~\ref{thm:BPapprox} that for the epidemic $\EEn$, conditional upon a major outbreak, $\In(t_n)/\Ine(t_n) \convas L$ as $n \to \infty$, where $t_n=\inf\{t \ge 0: \Tn(t)\ge \log n\}$ is the time when the cumulative number of infectives in $\EEn$ first reaches $\log n$. This suggests the initial conditions for the deterministic model used in Conjecture~\ref{conj:SIRfinal}(b).

\subsubsection{Proof of Theorem~\ref{thm:detdisc}}
\label{sec:detdisc}

The lack of a closed-form solution to~\eqref{equ:dsdtB}-\eqref{equ:dwdtB} also means that we do not obtain an equation that is satisfied by
the final size $\tau$.  We can however obtain a condition for the final size of the deterministic model
\eqref{equ:dsdtB}-\eqref{equ:dwdtB}, with initial condition $(\st(0),\itt(0),\iet(0),\wt(0))=(1-\epsilon,\epsilon, L\epsilon,0)$ where $L$ is a suitably chosen constant, to have a discontinuity at the threshold $\lambda=\lambda_c$ in the limit $\epsilon \downarrow 0$ and hence prove Theorem~\ref{thm:detdisc}.

First note that~\eqref{equ:diEdtB}, together with the initial condition, yields
\begin{equation}
\label{equ:iet0}
\iet'(0)=-1-\frac{\gamma+\omega}{\lambda}+\mu > 0 \text{ if and only if } \lambda>\frac{\gamma+\omega}{\mu-1}=\lambda_c.
\end{equation}
Differentiating~\eqref{equ:diEdtB} yields
\begin{align}
\dfrac{d^2\iet}{dt^2}&=\mu \dfrac{d\st}{dt}-\frac{1}{\st}\dfrac{d\iet}{dt}+\frac{1}{\st^2}\iet \dfrac{d\st}{dt}+\frac{2}{\st}\dfrac{d\wt}{dt}-\frac{2}{\st^2}\wt \dfrac{d\st}{dt}+\frac{\omega\alpha}{\lambda}\dfrac{d\itt}{dt}\label{equ:diet2}\\
&=-\mu-\frac{1}{\st}\dfrac{d\iet}{dt}-\frac{1}{\st^2}\iet+\frac{2\omega\alpha}{\lambda}-\frac{2\wt}{\st^2}
+\frac{\omega\alpha}{\lambda}\left(-\frac{\gamma}{\lambda}\frac{\itt}{\iet}+1\right).\nonumber
\end{align}

To determine $\iet''(0)$, we need $L=\lim_{t \downarrow 0}\frac{\itt(t)}{\iet(t)}$.  The value of $L$ is assumed in the theorem, see the definition of $\tau_{SIR}(\mu, \lambda, \gamma, \omega, \alpha)$ just prior to Theorem~\ref{thm:detdisc}, but it can also be obtained using l'Hopital's rule.  Assuming that the limit exists, l'Hopital's rule gives
\begin{equation*}
L=\lim_{t \downarrow 0}\frac{\itt'(t)}{\iet'(t)}
=\frac{-\frac{\gamma}{\lambda}{L}+1}{-1-\frac{\gamma+\omega}{\lambda}+\mu},
\end{equation*}
using~\eqref{equ:didtB} and~\eqref{equ:iet0}, whence
\begin{equation}
\label{equ:L}
L=\frac{\lambda}{\lambda(\mu-1)-\omega}.
\end{equation}

Substituting~\eqref{equ:L} into~\eqref{equ:diet2} yields
\begin{align}
\iet''(0)&=-\mu+\frac{\omega\alpha}{\lambda}\left(3-\frac{\gamma}{\lambda(\mu-1)-\omega}\right)\nonumber\\
&=-\mu+\frac{2\omega\alpha}{\lambda_c}=\frac{\mu[\omega(2\alpha-1)-\gamma]-2\omega\alpha}{\gamma+\omega}\label{equ:diet2a},
\end{align}
when $\lambda=\lambda_c$.
Thus, when $\lambda=\lambda_c$,
\begin{equation*}
\iet''(0)>0 \text{ if and only if } \gamma<\omega(2\alpha-1) \text{ and } \mu > \frac{2\omega\alpha}{\omega(2\alpha-1)-\gamma}
\end{equation*}
and
\begin{equation*}
\iet''(0)<0 \text{ if and only if } \gamma > \omega(2\alpha-1) \text{ or } \mu < \frac{2\omega\alpha}{\omega(2\alpha-1)-\gamma}.
\end{equation*}
Theorem~\ref{thm:detdisc} follows since
\[
\lim_{\lambda \downarrow \lambda_c} \tau_{SIR}(\mu, \lambda, \gamma, \omega, \alpha)
\begin{cases}
	      =0& \text{ if } \iet''(0)<0, \\
	      >0& \text{ if } \iet''(0)>0.
\end{cases}
\]

\subsubsection{Proof of Theorem~\ref{thm:SIRdiscon}}
\label{sec:proofSIRdisc}
Theorem~\ref{thm:SIRdiscon} is proved by considering modifications of the epidemic model $\EEn$ which (i) give upper and lower bounds for the process of infectives in $\EEn$ and (ii) lead to a time-transformed deterministic model whose corresponding vector field is Lipschitz thus enabling proof of an associated law of large numbers.  Part (a) of Theorem~\ref{thm:SIRdiscon} follows by showing that under the given conditions the final size of the deterministic model
for the lower bounding process has a discontinuity at the threshold $\lambda=\lambda_c$; part (b) follows similarly by showing that the final size of
deterministic model
for the upper bounding process is continuous at $\lambda=\lambda_c$.
We present the argument in detail for the lower bounding process.  The proof for the upper bounding process is similar, so only an outline is given.

Suppose first that the model $\EEn$ is modified so that if a susceptible rewires an edge from one infective
to another infective then the edge is dropped.  This clearly leads to a model whose final size is
stochastically smaller than that of the original model $\EEn$.  We run the original model until time
$t_n=\inf\{t>0:\Sn(t) \le n-\log n\}$ and then the modified model, and make the same random time-scale transformation as in Section~\ref{sec:proofthmSIfinal} to the latter.  For the time-transformed modified model, in an obvious notation, let
$\Xnhat=\{(\Snhat(t), \Inhat(t), \Inehat(t), \Wnhat(t)):t \ge 0\}$, where $\Snhat(0)=\Sn(t_n)$ etc.  The corresponding deterministic model is
\begin{align}
\dfrac{d\sh}{dt}&=-1, \label{equ:dsdtBL}\\
\dfrac{d\ih}{dt}&=-\frac{\gamma}{\lambda}\frac{\ih}{\ieh}+1,\label{equ:didtBL}\\
\dfrac{d\ieh}{dt}&= -1-\frac{\gamma}{\lambda}+\mu \sh - \frac{\ieh}{\sh} +2 \frac{\wh}{\sh}-\frac{\omega}{\lambda},\label{equ:diEdtBL}\\
\dfrac{d\wh}{dt}&= \frac{\omega \alpha}{\lambda} \sh-2 \frac{\wh}{\sh},\label{equ:dwdtBL}
\end{align}
with initial condition $(\sh(0),\ih(0),\ieh(0),\wh(0))=(1,0,0,0)$.
Note that the only difference between~\eqref{equ:dsdtB}-\eqref{equ:dwdtB} and~\eqref{equ:dsdtBL}-\eqref{equ:dwdtBL} is that the final term in~\eqref{equ:diEdtB}
is replaced by $\frac{\omega}{\lambda}$ since now
every rewiring necessarily leads to a drop of an infective edge.

For $t \ge 0$, let $\bXnhat(t)=n^{-1}(\Snhat(t), \Inehat(t), \Wnhat(t))$ and $\xh(t)=(\sh(t),\ieh(t),\wh(t))$.
Let $\zetah=\inf\{t>0:\ieh(t)=0\}$.  Observe that~\eqref{equ:dsdtBL}, \eqref{equ:diEdtBL} and~\eqref{equ:dwdtBL} form an autonomous system that is Lipschitz provided that $\sh$ is bounded away from $0$.  Also, in an obvious notation, \eqref{equ:alphan1} becomes
\[
\alphahn_R(\xi)=n^{-2}\left(\frac{\gamma}{\lambda}\right)\frac{n^E_2}{n^E_2},
\]
leading to the bound (cf.~\eqref{equ:alphanb3})
\begin{equation*}
\alphatn_R(\Xnhat(t))\le 2 \frac{\gamma}{\lambda}  n^{-1} \log n.
\end{equation*}
Consequently, it is straightforward to modify the proof of Theorem~\ref{thm:temporalWLLN} to show that, for any $t_0 \in (0,\zetah)$,
\begin{equation}
\label{equ:modifedLLN}
\sup_{0 \le t \le t_0} \left|\bXnhat(t)-\xh(t)\right| \convp 0 \quad \mbox{as } n \to \infty.
\end{equation}

Now $\ieh'(0)=\mu-1-\frac{\gamma+\omega}{\lambda}$, so $\ieh'(0)>0$ if and only if $\lambda>\lambda_c$ and
$\ieh'(0)=0$ if $\lambda=\lambda_c$.  Further, using~\eqref{equ:dsdtBL},
\[
\dfrac{d^2\ieh}{dt^2}=-\mu -\frac{1}{\sh}\dfrac{d\ieh}{dt}-\frac{1}{\sh^2}\ieh +\frac{2}{\sh}\dfrac{d\wh}{dt}+\frac{2\wh}{\sh^2},
\]
so, using~\eqref{equ:dwdtBL}, when $\lambda=\lambda_c$,
\[
\ieh''(0)=-\mu+\frac{2\omega \alpha}{\lambda_c}=\frac{\mu(\omega(2\alpha-1)-\gamma)-2\omega\alpha}{\omega+\gamma}.
\]
Note that $\ieh''(0)>0$ if and only if the condition in Theorem~\ref{thm:SIRdiscon}(a) is satisfied.

Let $\zetahn=\inf\{t \ge 0: \Inehat(t)=0\}$ and $\Tnhat=n-\Snhat(\zetahn)$, so $\Tnhat$ is the size of the
modified epidemic.  Suppose that the condition in Theorem~\ref{thm:SIRdiscon}(a) is satisfied.
Then $\lim_{\lambda \downarrow \lambda_c} (1-\xh(\zetah))>0$ and it follows using~\eqref{equ:modifedLLN}
that there exists $\tau_0>0$ such that
\[
\lim_{\lambda \downarrow \lambda_c} \lim_{n \to \infty} \P(n^{-1}\Tnhat > \tau_0)=1.
\]
Theorem~\ref{thm:SIRdiscon}(a) now follows as $\Tnbar$, conditioned on the occurrence of a major outbreak, is stochastically larger than $n^{-1}\Tnhat$.

Turning to part (b), suppose now that the model $\EEn$ is modified so that if a susceptible rewires an edge from an infective to a recovered individual then
the edge to the infective is also retained.  In the construction this means that when an infective sends a warning down an infective edge and that edge is rewired to a recovered individual then the infective does not lose that infective edge.   This leads to a model whose final size is stochastically larger than that of $\EEn$. As before, we run the model $\EEn$ until time $t_n$ defined above and then the time-transformed modified model. The corresponding time-transformed deterministic model is given by~\eqref{equ:dsdtBL}-\eqref{equ:dwdtBL}, with~\eqref{equ:diEdtBL} replaced by
\begin{equation}
\label{equ:diEdtBU}
\dfrac{d\ieh}{dt}= -1-\frac{\gamma}{\lambda}+\mu \sh - \frac{\ieh}{\sh} +2 \frac{\wh}{\sh}-\frac{\omega}{\lambda}(1-\alpha+\alpha \sh).
\end{equation}
Note that \eqref{equ:dsdtBL}, \eqref{equ:diEdtBU} and~\eqref{equ:dwdtBL} form an autonomous system that is Lipschitz provided that $\sh$ is bounded away from $0$.  The proof parallels that of part (a) in the obvious fashion, the only difference being that now, when $\lambda=\lambda_c$,
\[
\ieh''(0)=-\mu+\frac{3\omega \alpha}{\lambda_c}=\frac{\mu(\omega(3\alpha-1)-\gamma)-3\omega\alpha}{\omega+\gamma}.
\]
Thus $\ieh''(0)<0$ if either $\omega(3\alpha-1) < \gamma$ or $\mu < \frac{3\omega\alpha}{\omega(3\alpha-1)-\gamma}$.  Moreover, a further calculation shows that if $\lambda=\lambda_c$ and $\ieh''(0)=0$ then $\ieh'''(0)=-\frac{2\omega\alpha}{\lambda_c}<0$.  Part (b) of the theorem now follows using a similar argument to that used at the end of the proof of part (a).

The bounding process used in the proof of Theorem~\ref{thm:SIRdiscon}(b) suggests an approach to proving the conditions in part(a) are also necessary for the final size to be discontinuous at $\lambda=\lambda_c$.  Run the epidemic $\EEn$ until time $t_n$ and then make the random time-scale transformation.
For the latter, let
$\bXnt=\{\bXnt(t):t \ge 0\}$, where $\bXnt(t)=\bxt(\Xnt(t))$ $(t \ge 0)$ and $\bxt: \En \to \mathbb{R}^4$ is defined by~\eqref{equ:coordb}.
Let $\Int(t)$ and $\Inet(t)$ denote respectively the numbers of infectives and infective edges in $\bXnt$ at time $t$.
Recall from Section~\ref{sec:detdisc} that for the epidemic $\EEn$, conditional upon a major outbreak, $\In(t_n)/\Ine(t_n) \convas L$ as $n \to \infty$.  This suggests that there exists $L'>L$ and $t_0>0$ such that, conditional upon a major outbreak, $\Int(t) \le L' \Inet(t)$ for all $t \in [0,t_0]$ with high probability.

The term $\frac{\omega}{\lambda} (1-\alpha+\alpha(1-\itt))$ in~\eqref{equ:diEdtB} arises from rewirings to
susceptibles or recovered.  In particular, the probability that a rewiring is to a recovered, and the infective edge is dropped, is $\frac{n-\Xn(t)-\In(t)-2}{n-2}$.  Modifying the model so that this probability is replaced by $\frac{n-\Xn(t)-L'\Ine(t)-2}{n-2}$ leads to an upper bounding process for the true model and hence also for the time-transformed model.
The latter has deterministic limit which satisfies~\eqref{equ:dsdtB}-\eqref{equ:dwdtB} over $[0,t_0]$ but
with~\eqref{equ:diEdtB} replaced by
\begin{equation}
\label{equ:diEdtBU1}
\dfrac{d\iet}{dt} = -1-\frac{\gamma}{\lambda}+\mu \st - \frac{\iet}{\st} +2 \frac{\wt}{\st}-\frac{\omega}{\lambda} (1-\alpha+\alpha(1-L'\iet)).
\end{equation}
Thus~\eqref{equ:dsdtB}, \eqref{equ:diEdtBU1} and~\eqref{equ:dwdtB} form an autonomous system, that is Lipschitz provided that $\sh$ is bounded away from $0$ and hence susceptible to analysis.  Omitting the details, when $\lambda=\lambda_c$ we obtain $\iet''(0)=-\mu+\frac{2\omega\alpha}{\lambda_c}$ and
$\iet'''(0)=-\frac{2\omega\alpha}{\lambda_c}<0$ if $\iet''(0)=0$, which imply that the conditions in Theorem~\ref{thm:SIRdiscon}(b) are also necessary.

It may be possible to make the argument rigorous using branching processes that bound $\EEn$
but that does not seem straightforward since (i) $t_0$ depends on $\lambda$ and $\downarrow 0$ as $\lambda \downarrow \lambda_c$ and (ii) the result is conditional upon a major epidemic, the probability of which also $\downarrow 0$ as $\lambda \downarrow \lambda_c$.

\subsection{Proof of Theorems~\ref{thm:temporalWLLNsusonly} and~\ref{thm:SIfinalsizesusonly}}
\label{sec:proofsusonly}
The proof of Theorem~\ref{thm:temporalWLLNsusonly} is omitted as it parallels that of Theorem~\ref{thm:temporalWLLN} in Section~\ref{sec:proofthmtemporalWLLN}, with only very minor modification.  Turning to Theorem~\ref{thm:SIfinalsizesusonly}, we make the same random time-scale transformation as in Section~\ref{sec:proofthmSIfinal}, in which the clock is slowed down by a factor given by the total force of infection. Recalling~\eqref{equ:dsdt1}-\eqref{equ:dwdt1}, this leads to the deterministic model
given by
\begin{align}
\dfrac{d\st}{dt}&=-1, \label{equ:dsdtC}\\
\dfrac{d\itt}{dt}&=-\frac{\gamma}{\lambda}\frac{\itt}{\iet}  + 1\label{equ:didtC}\\
\dfrac{d\iet}{dt}&= -1-\frac{\gamma}{\lambda}+\mu \st - \frac{\iet}{\st} +2 \frac{\wt}{\st}-\frac{\omega}{\lambda},\label{equ:diEdtC}\\
\dfrac{d\wt}{dt}&=\frac{\omega \alpha}{\lambda} -2 \frac{\wt}{\st},\label{equ:dwdtC}
\end{align}
with initial condition $(\st(0), \itt(0), \iet(0), \wt(0))=(1,0,0,0)$.  The equations for $(\st, \iet, \wt)$ form a closed system having solution
\begin{align}
\st(t)&=1-t,\label{equ:stSO}\\
\iet(t)&=\st(t)\gt(\st(t)),\label{equ:ietSO}\\
\wt(t)&=\frac{\omega\alpha}{\lambda}\st(t)(1-\st(t)),\label{equ:wtSO}
\end{align}
where
\[
\gt(s)=\left(1+\frac{\gamma+\omega(1-2\alpha)}{\lambda}\right)\log \st+\left(\mu-2\frac{\omega\alpha}{\lambda}\right)(1-\st(t)).
\]

The solution is derived as follows.  Equation~\eqref{equ:stSO} follows immediately from~\eqref{equ:dsdtC} and the initial condition.  Equation~\eqref{equ:wtSO}
is obtained by dividing~\eqref{equ:dwdtC} by~\eqref{equ:dsdtC} and
solving the resulting differential equation using the integrating factor $\st^{-2}$, together with the initial condition.  Finally, \eqref{equ:ietSO} is
obtained by dividing~\eqref{equ:diEdtC} by~\eqref{equ:dsdtC}, substituting for $\wt$ from~\eqref{equ:wtSO} and then solving the resulting differential equation using the integrating factor $\st^{-1}$, together with the initial condition.

Let $\zetat=\inf\{t >0:\iet(t)=0\}$ and $\hat{s}=\st(\zetat)$.  Then either $\hat{s}=0$ or
$\gt(\hat{s})=0$.  We thus investigate the solutions of $\gt(s)=0$ in $[0,1]$.  Note that
$\iet'(0)=\mu-1-\frac{\omega+\gamma}{\lambda}>0$, as $\lambda > \lambda_c$, and $\st(1)=0$,
so $\zetat \in (0,1]$, so we investigate the solutions of $\gt(s)=0$ in $[0,1)$.  Now
\[
\gt'(s)=\left(1+\frac{\gamma+\omega(1-2\alpha)}{\lambda}\right)s^{-1}-\mu+2\frac{\omega\alpha}{\lambda},
\]
so $\gt$ has at most one stationary point in $(0,1]$.  Further, $\gt(1) =0$ and
$\gt'(1)=1-\mu+\frac{\omega+\gamma}{\lambda} < 0$, as $\lambda > \lambda_c$, so by considering $\gt(0+)$ it is seen easily that $\gt$ has no root in $(0,1)$ if $\lambda+\gamma+\omega(1-2\alpha)\le 0$ and precisely one root in $(0,1)$ if $\lambda+\gamma+\omega(1-2\alpha)> 0$.  Note that since $\lambda_c=\frac{\gamma+\omega}{\mu-1}$, a necessary (but not sufficient) condition for $\gt$ not to have a root in $(0,1)$ is that
$r(\mu, \gamma, \omega, \alpha)<0$, where $r(\mu, \gamma, \omega, \alpha)$ is defined at~\eqref{equ:rdef}.

Suppose that a major epidemic occurs.  Let $t_n=\inf\{t>0:\Sn(t) \le n-\log n\}<\infty$.
Define the time-transformed process $\bXnt$ analogously to the SI model, using the coordinate functions~\eqref{equ:cordbt}, with $\bXnt(0)=(n^{-1}(\Sn(t_n), \Ine(t_n),\Wn(t_n))$.

Suppose that $r(\mu, \gamma, \omega, \alpha)<0$ and
$\lambda \in (\lambda_c, \omega(2\alpha-1)-\gamma]$.  Then, given any $\epsilon\in (0,1)$, let $t_1=1-\epsilon/2$,
so $\st(t_1)=\epsilon/2$.  Let $\tilde{U}=[\frac{\epsilon}{2},1]\times[-\frac{\epsilon}{2},2\mu]\times [0, 2 \mu]$.  The vector field
$\bt: \Ut \to \mathbb{R}^3$
given by the right-hand sides of~\eqref{equ:dsdtC}, \eqref{equ:diEdtC} and~\eqref{equ:dwdtC} is Lipschitz and the proof of Theorem~\ref{thm:temporalWLLN}
is easily modified to yield
\begin{equation}
\label{equ:susonlywc}
\lim_{n \to \infty} \P\left(\sup_{0 \le t \le t_1} \left|\bXnt(t)-\xt(t)\right|>\frac{\epsilon}{4}\right)=0,
\end{equation}
where $\xt(t)=(\st(t), \iet(t), \wt(t))$.  By a similar argument to that used at the end of the proof of Theorem~\ref{thm:SIfinalsize}(b) in Section~\ref{sec:proofthmSIfinal}, there exists $t_2 \in (0,t_1)$ such that $\Inet(t)>0$ for all $t \in [0, t_2]$ with high probability.
It follows that
\begin{align*}
\lim_{n \to \infty}\P(\Tnbar>1-\epsilon)&\ge \lim_{n \to \infty}\P(\Snt(t_1)<\epsilon)\\
&=1,
\end{align*}
using~\eqref{equ:susonlywc} and the definition of $t_1$.  Theorem~\ref{thm:SIfinalsizesusonly}(b) for the case when
$\lambda \in (\lambda_c, \omega(2\alpha-1)-\gamma]$ now follows as $\epsilon\in (0,1)$ is arbitrary.

The proof of Theorem~\ref{thm:SIfinalsizesusonly} for the other cases, i.e.~when the limit $\taut<1$, follows a similar argument to that of
Theorem~\ref{thm:SIfinalsize}(b) and hence is omitted.

\section{Concluding comments}
\label{sec:conc}

We have presented a construction of an SIR epidemic with preventive rewiring on an Erd\H{o}s-R\'{e}nyi random graph, together with a rigorously justified
deterministic approximation.  For the special case of the SI model, these yielded a detailed analysis of the final outcome of the epidemic, and in particular a proof of a necessary and sufficient condition for the final size to be discontinuous at the phase transition $\lambda=\lambda_c$.
Similar results were also obtained for the SIR model when rewiring is necessarily to a susceptible individual. Moreover in that case the behaviour at the phase transition is very striking in that the final fraction infected by the epidemic jumps from $0$ to $1$.  For the original SIR model only partial results were proved, in that there is a non-negligible gap between the necessary and sufficient conditions for a discontinuity, although stronger results were conjectured together with supporting evidence.  Proving those conjectures is a worthwhile future work.

It seems likely that similar results hold for epidemics on networks constructed using the configuration model.  Indeed~\cite{YD:20} provides
such results although, as described in Section~\ref{sec:YaoDurrett}, we believe that their results are for a model which is different from the original
SIR model and does not provide an exact asymptotic approximation to the original model in the limit as $n \to \infty$.  However it seems hard to extend our construction to epidemics on configuration-model networks since the construction exploits symmetries and independence properties of a Erd\H{o}s-R\'{e}nyi random graph that are not present in configuration-model networks.  Although not considered here, it seems likely that the construction can be extended to
models that are more akin to the Erd\H{o}s-R\'{e}nyi random graph, such as the stochastic block model.

\appendix

\section{Time-transformed deterministic SI model.}
\label{app:ttdetSI}
\subsection{Solution of time-transformed deterministic model~\eqref{equ:dsdtA}-\eqref{equ:dwdtA}}
\label{app:ttdetSIsoln}
We solve~\eqref{equ:dsdtA}-\eqref{equ:dwdtA} with the initial condition $(\st(0), \iet(0),\wt(0))=(1-\epsilon, \mu \epsilon(1-\epsilon),0)$.
The solution~\eqref{equ:st} for $\st(t)$ follows immediately from~\eqref{equ:dsdtA} as $\xt(0)=1-\epsilon$.  Dividing~\eqref{equ:dwdtA} by~\eqref{equ:dwdtA} yields
\[
\dfrac{d\wt}{d\st}= -\frac{\omega \alpha}{\lambda} \st+2 \frac{\wt}{\st},
\]
which with the above initial condition yields, using the integrating factor $\st^{-2}$,
\begin{equation}
\label{equ:equ:wtapp}
\wt(\st)=\frac{\omega \alpha}{\lambda}\st^2\log\left(\frac{\st}{1-\epsilon}\right),
\end{equation}
and~\eqref{equ:wt} follows.  Dividing~\eqref{equ:diEdtA} by\eqref{equ:dwdtA} and substituting from~\eqref{equ:equ:wtapp} gives
\[
\dfrac{d\iet}{d\st}=1 -\mu \st + \frac{\iet}{\st} +2\frac{\omega\alpha}{\lambda}\st \frac{\wt}{\st}
\log\left(\frac{\st}{1-\epsilon}\right)
+\frac{\omega}{\lambda} (1-\alpha+\alpha \st),
\]
which, with the above initial condition can be solved using the integrating factor $\st^{-1}$ to yield~\eqref{equ:iet}.

\subsection{Final size and $\iet'(\zetat)$ when $\epsilon>0$}
\label{app:ttdeteps>0}
Define $\fe: [\epsilon,1) \to \mathbb{R}$ by
\[
\fe(x)=\log\left(\frac{1-x}{1-\epsilon}\right)+\frac{(\lambda\mu+\omega\alpha)x-\omega\alpha\epsilon}
{\lambda+\omega(1-\alpha)+2\omega\alpha(1-x)},
\]
so $\tau$ satisfies $F_{\epsilon}(\tau)=0$ if and only if $\fe(\tau)=0$.
From~\eqref{equ:iet}, we have $\iet(t)=H(\st(t))$, where
\[
H(\st)=\frac{\st}{\lambda}G(\st), \text{ with } G(\st)=(\lambda+\omega(1-\alpha)+2\omega\alpha \st)\fe(1-\st).
\]
Now $\fe(\epsilon)>0$ and $\fe(1-)=-\infty$, so $\fe$ has at least one root in $[\epsilon,1)$, and
since $\st(t)=1-\epsilon-t$ $(t \ge 0)$, $\zetat=\inf\{t>0:\iet(t)=0\}$ is finite and $\st(\zetat)=1-\tau$, where $\tau$ is the smallest solution of $F_{\epsilon}(x)=0$ in $(\epsilon,1)$.
Further, since $\st'(\zetat)=-1$,
\begin{align*}
\iet'(\zetat)&=-H'(\st(\zetat))\\
&=-\frac{1}{\lambda}\left[\st(\zetat)G'(\st(\zetat))+G(\st(\zetat))\right]\\
&=-\frac{\st(\zetat)}{\lambda}G'(\st(\zetat)),
\end{align*}
as $G(\st(\zetat))=0$.  Now,
\[
G'(\st(\zetat))=-(\lambda+\omega(1-\alpha)+2\omega\alpha \st(\zetat))\fe'(1-\st(\zetat)),
\]
since $\fe(1-\st(\zetat))=0$, so $\iet'(\zetat)<0$ if and only if $\fe'(\tau)<0$.  It is shown easily that
$\fe'(\tau)<0$ if and only if $F_{\epsilon}'(\tau)<0$.

\subsection{Final size and $\iet'(\zetat)$ when $\epsilon=0$}
\label{app:ttdeteps=0}

Setting $\epsilon=0$ in~\eqref{equ:iet} shows that $\tau=1-\st(\zetat)$ is given by the smallest solution in $(0,1)$ of $F_0(x)=0$, provided there is at least one solution, otherwise $\tau=1$.  We show now that when $R_0>1$ there is a unique solution in $(0,1)$.  Let $f$ be the function defined at~\eqref{equ:f},
so $\tau$ satisfies $F_0(\tau)=0$ if and only if $f(\tau)=0$.  Now
\[
f'(x)=-\frac{1}{1-x}+\frac{(\lambda\mu+\omega\alpha)[\lambda+\omega(1+\alpha)]}
{[\lambda+\omega(1-\alpha)+2\omega\alpha(1-x)]^2},
\]
so, for $x \in (0,1)$, $f'(x)=0$ if and only if $g_1(1-x)=g_2(1-x)$, where
$g_1(x)=(\lambda+\omega(1-\alpha)+2\omega\alpha x)^2$ and $g_2(x)=(\lambda\mu+\omega\alpha)[\lambda+\omega(1+\alpha)]x$.
Now $g_1(0)>g_2(0)$, $\lim_{x \to \infty}g_1(x)-g_2(x)=\infty$ and a simple calculation shows that
$g_1(1)<g_2(1)$ if and only if $\lambda>\frac{\omega}{\mu-1}=\lambda_c$.  Suppose $R_0>1$, so
$\lambda>\lambda_c$.  Then $g_1-g_2$ has a root in
$(0,1)$ and a root in $(1, \infty)$.  Moreover, $g_1-g_2$ has precisely one root in each of these intervals, as it convex on $\mathbb{R}$.  Hence $f'$ has a unique root in $(0,1)$.  Further, $f(0)=0$,
$f'(0)=-1+(\lambda\mu+\omega\alpha)/[\lambda+\omega(1+\alpha)]>0$, as $\lambda>\lambda_c$, and
$f(1-)=-\infty$, so the unique stationary point of $f'$ in $(0,1)$ must be a maximum.  It follows that
$f(x)=0$, and hence also $F_0(x)=0$ has a unique solution $\tau$ in $(0,1)$ and $f'(\tau)<0$.
Letting $\epsilon=0$ in the argument in Appendix~\eqref{app:ttdeteps>0} shows that $\iet'(\zetat)<0$.

\section*{Acknowledgements}

Tom Britton is grateful to the Swedish Research Council (grant 2015-05015) for financial support.


\begin{thebibliography}{1}


\bibitem[Altmann (1998)]{Altmann98}
Altmann, M. (1998)
The deterministic limit of of infectious disease models with dynamic partners.
{\sl Math.~Biosci.} {\bf 150}, 153--175.




\bibitem[Andersson and Britton(2000)]{AB00}
Andersson, H. and Britton, T.~(2000)
{\it Stochastic Epidemic Models and Their Statistical Analysis}.
Springer, New York.







\bibitem[Ball and Britton(2007)]{BB:2007}
Ball, F. and Britton, T.~(2007) An epidemic model with infector-dependent severity.
{\sl Adv.~Appl.~Prob.} {\bf 39}, 949--972.




\bibitem[Ball and Neal(2008)]{BN:2008}
Ball, F. and Neal, P.~(2008) Network epidemic models with two levels of mixing.
{\sl Math.~Biosci.} {\bf 212}, 69--87.





\bibitem[Barbour and Reinert(2013)]{BR:2013}
Barbour, A. and Reinert, G.~(2013) Approximating the epidemic curve.
{\sl Electron.~J.~Probab.} {\bf 18}(54), 1--30.






\bibitem[Bollob\'as(1980)]{Bollobas:1980}
Bollob\'as, B.~(1980) A probabilistic proof of an asymptotic formula for the number of labelled regular graphs.
{\sl European J.~Combin.} {\bf 1}, 311-316.




\bibitem[Britton and O'Neill(2002)]{BO'N:2002}
Britton, T. and O'Neill, P.D.~(2002)
Bayesian inference for stochasric epidemics in populations with random social structure.
{\sl Scand.~J.~Stat.} {\bf 29}, 375--390.



\bibitem[Britton and Trapman(2012)]{BT:2012}
Britton, T. and Trapman, P.~(2012)
Maximizing the size of the giant.
{\sl J.~Appl.~Prob.} {\bf 49}, 1156--1165.

\bibitem[Chinazzi et al.(2020)]{C20}
Chinazzi, M., Davis, J.T., Ajelli, M., Gioannini, C., Litvinova, M., Merler, S., Pastore y Piontti, A., Mu, K., Rossi, L., Sun, K., Viboud, C., Xiong, X., Yu, H., Halloran, M.E., Longini, I.M. and Vespignani, A. (2020) The effect of travel restrictions on the spread of the 2019 novel coronavirus (COVID-19) outbreak.
{\sl Science} {\bf 368}, 395--400.


\bibitem[Darling and Norris(2008)]{DN:08}
Darling, R.W.R. and Norris, J.R.~(2008)
Differential equation approximations for Markov chains.
{\sl Probability Surveys} {\bf 5}, 37--79.

\bibitem[Durrett(2007)]{Durrett:2007}
Durrett, R. (2007)
{\it Random Graph Dynamics}.
Cambridge University Press, Cambridge.





\bibitem[Ethier and Kurtz(1986)]{EK86}
Ethier, S.N. and Kurtz, T.G.~(1986)
{\it Markov Processes: Characterization and Convergence}.
Wiley, New York.

\bibitem[Ferguson et al.(2006)]{F06}
Ferguson, N.M., Cummings, D.A.T., Fraser, C., Cajka, J.C., PC Cooley P.C.\ and Burke, D.S.~(2006) Strategies for mitigating an influenza pandemic.
{\sl Nature} {\bf 442}, 448-452.


\bibitem[Funk et al.(2010)]{FSJ10} Funk, S., Salathe, M., and Jansen, V.A.A.~(2010)  Modelling the influence of human behaviour onthe spread of infectious diseases:  a review. {\sl J.~R.~Soc.~Interface}, {\bf 7}:1247–1256, doi:10.1098/rsif.2010.0142.

\bibitem[Gross et al.(2006)]{gross06}
Gross, T., Dommar D'Lima, C.~J.~and Blasius, B.~(2006)
Epidemic dynamics on an adaptive network.
{\sl Phys.~Rev.~Letters} {\bf 96}, 208701.

\bibitem[van der Hofstad(2016)]{vdHofstad:2016}
van der Hofstad, R.~(2016)
{\it Random Graphs and Complex Networks Volume 1}.
Cambridge University Press, Cambridge.


\bibitem[Jacobsen et al.(2018)]{JBTR:2018}
Jacobsen, K.A., Burch, M.~G., Tien, J.~H.~and Rempala, G.~A.~(2018)
The large graph limit of a stochastic epidemic model on a dynamic
multilayer network.
{sl J.~Biol.~Dyn.} {\bf 12}, 746--788.


\bibitem[Janson(2009)]{Janson:2009}
Janson, S.~(2009)
The probability that a random multigraph is simple.
{\sl Combin.~Prob.~Comput.} {\bf 18}, 205--225.





\bibitem[Janson et al.(2014)]{JLW:2014}
Janson, S., Luczak, M. and Windridge, P.~(2014)
Law of large numbers for the {SIR} epidemic on a random graph with given degrees.
{\sl Random Structures Algorithms} {\bf 45}, 726--763.

\bibitem[Jiang et al.(2019)]{JKYJD:2019}
Jiang, Y., Kassem, R., York, G., Junge, M.~and Durrett, D.~(2019)
SIR epidemics on evolving graphs.
{\sl arXiv}:1901.06568v1.

\bibitem[Kimber(1983)]{Kimber83}
Kimber, A.C.~(1983)
A note on Poisson maxima.
{\sl Probability Theory and Related Fields} {\bf 63}, 551–-552.






\bibitem[Leung et al.(2018)]{LBSB18} Leung, K., Ball, F., Sirl, D.\ and Britton, T. (2018).
Individual preventive social distancing during an epidemic may have negative population-level outcomes.
{\sl J.~R.~Soc.~Interface} {\bf 15}:20180296.

\bibitem[Martin-L{\"o}f(1986)]{ML:1986}
Martin-L{\"o}f, A.~(1986) Symmetric sampling procedures, general epidemic processes and their threshold limit theorems.
{\sl J.~Appl.~Prob.} {\bf 23}, 265--282.



\bibitem[Molloy and Reed(1995)]{MR:1995}
Molloy, M. and Reed, B.~(1995) A critical point for random graphs with a given degree sequence.
{\sl Random Structures Algorithms} {\bf 6}, 161--179.

\bibitem[Neal(2003)]{Neal:2003}
Neal, P.~(2003) SIR epidemics on a Bernoulli random graph.
{\sl J.~Appl.~Prob.} {\bf 23}, 265--282.

\bibitem[Nerman(1981)]{Nerman81}
Nerman, O.~(1981) On the convergence of supercritical general (C-M-J)
branching processes.
{\sl Z. Wahrscheinlichkeitsth} {\bf 57}, 365--395.


\bibitem[Newman et al.(2001)]{NSW:2001}
Newman, M.E.J., Strogratz, S.H. and Watts, D.J.~(2001)
Random graphs with arbitrary degree distributions and their applications.
{\sl Phys.~Rev.~E} {\bf 64}, 026118.






\bibitem[Watson(1980)]{Watson80}
Watson, R.~(1980) A useful random time-scale transformation for the standard epidemic model.
{\sl J.~Appl.~Prob.} {\bf 17}, 324--332.


\bibitem[Yao and Durrett(2020)]{YD:20}
Yao, D. and Durrett, R.~(2020) Epidemics on evolving graphs.
{\sl arXiv}:2003.08534v1.
\end{thebibliography}
\end{document}